\documentclass[a4paper,12pt,amsart,frenchb]{article}

\usepackage{amsmath,amsbsy,amsfonts,amssymb}
\usepackage[french]{babel}
\newcommand{\ass}[2]{\vskip0.3cm\noindent
{\bf {#1}}. { \sl {#2}}\vskip0.3cm\noindent
}
  
\begin{document}

 \title{  La conjecture locale de Gross-Prasad pour les  groupes sp\'eciaux orthogonaux: le cas g\'en\'eral}
\author{C. Moeglin, J.-L. Waldspurger}
\date{31 d\'ecembre 2009}
\maketitle

{\bf Introduction}

Soit $F$ un corps local non archim\'edien de caract\'eristique nulle.  On appelle espace quadratique un couple $(V,q)$ form\'e d'un espace vectoriel $V$ sur $F$ de dimension finie et d'une forme bilin\'eaire $q$ sur $V$, sym\'etrique et non d\'eg\'en\'er\'ee. Soient $(V,q)$ et $(V',q')$ deux espaces quadratiques. On note $d$ et $d'$ leurs dimensions et $G$ et $G'$ leurs groupes sp\'eciaux orthogonaux. On suppose $d>d'$ et $d$ et $d'$ de parit\'es distinctes. On suppose donn\'ee une d\'ecomposition orthogonale  de $(V,q)$ en la somme de $(V',q')$ et d'un espace quadratique qui est lui-m\^eme somme orthogonale de plans hyperboliques et d'une droite quadratique $(D,q_{D})$. Le groupe $G'$  est alors un sous-groupe de $G$.  On d\'efinit $\nu_{0}\in F^{\times}/F^{\times,2}$ en fixant un \'el\'ement non nul $v_{0}\in D$ et en posant
$$\nu_{0}= (-1)^dq_{D}(v_{0},v_{0})/2 .$$
Soient $\sigma$, resp. $\sigma'$, une repr\'esentation admissible irr\'eductible de $G(F)$, resp. $G'(F)$. Gross et Prasad ont d\'efini une multiplicit\'e $m(\sigma,\sigma')$. Par exemple, dans le cas o\`u $d=d'+1$, introduisons des espaces $E_{\sigma}$ et $E_{\sigma'}$ dans lesquels se r\'ealisent $\sigma$ et $\sigma'$. Alors $m(\sigma,\sigma')$ est la dimension de l'espace complexe des applications lin\'eaires $l:E_{\sigma}\to E_{\sigma'}$ telles que $l\circ \sigma(g')=\sigma'(g')\circ l$ pour tout $g'\in G'(F)$. La d\'efinition g\'en\'erale est rappel\'ee en 1.2. D'apr\`es [AGRS] et [GGP] corollaire 15.2, on a toujours $m(\sigma,\sigma')\leq1$.

Rappelons la classification conjecturale des repr\'esentations admissibles irr\'eductibles de  $G(F)$ dans le cas $d$ impair.  Pour tout entier naturel pair $N$, on fixe une forme symplectique sur ${\mathbb C}^N$ et on note $Sp(N,{\mathbb C})$ son groupe symplectique. Notons $W_{F}$ le groupe de Weil absolu de $F$ et $W_{DF}=W_{F}\times SL(2,{\mathbb C})$ le groupe de Weil-Deligne. Notons $\Phi(G)$ l'ensemble des classes de conjugaison par $Sp(d-1,{\mathbb C})$ d'homomorphismes  continus $\varphi:W_{DF}\to Sp(d-1,{\mathbb C})$ qui sont semi-simples et dont la restriction \`a $SL(2,{\mathbb C})$ est alg\'ebrique. On conjecture que l'ensemble des classes de repr\'esentations admissibles irr\'eductibles de $G(F)$ est union disjointe de $L$-paquets $\Pi^G(\varphi)$ index\'es par les $\varphi\in \Phi(G)$. Cette classification se ram\`ene \`a celle des repr\'esentations temp\'er\'ees de la fa\c{c}on suivante.  Consid\'erons les donn\'ees suivantes:

$\bullet$ un L\'evi $\hat{L}=GL(d_{1},{\mathbb C})\times ...\times GL(d_{t},{\mathbb C})\times Sp(d_{0}-1,{\mathbb C})$ de $Sp(d-1,{\mathbb C})$;

$\bullet$ des homomorphismes $\varphi_{0}:W_{DF}\to Sp(d_{0}-1,{\mathbb C})$ et $\varphi_{j}:W_{DF}\to GL(d_{j},{\mathbb C})$ pour $j=1,...,t$, v\'erifiant les m\^emes conditions que $\varphi$ et qui sont temp\'er\'es, c'est-\`a-dire que les images de $W_{F}$ par ces homomorphismes sont relativement compactes;

$ \bullet$ des r\'eels $b_{1}>b_{2}>...>b_{t}>0$.

Notons $\vert .\vert _{F}$ la valeur absolue usuelle de $F^{\times}$, que l'on identifie par la th\'eorie du corps de classes \`a un caract\`ere de $W_{F}$, puis de $W_{DF}$. Introduisons l'homomorphisme 
$$\varphi^{\hat{L}}=(\varphi_{1}\otimes \vert .\vert _{F}^{b_{1}})\otimes  ...\otimes(\varphi_{t}\otimes\vert .\vert _{F}^{b_{t}})\otimes \varphi_{0}$$
de $W_{DF}$ dans $\hat{L}$. En le poussant par l'inclusion de $\hat{L}$ dans $Sp(d-1,{\mathbb C})$, il devient un \'el\'ement de $\Phi(G)$. Inversement, soit $\varphi\in\Phi(G)$. Alors il existe des donn\'ees comme ci-dessus, uniques \`a conjugaison pr\`es, de sorte que $\varphi=\varphi^{\hat{L}}$. Il peut ne correspondre \`a $\hat{L}$ aucun L\'evi de $G$ (c'est le cas si $G$ n'est pas d\'eploy\'e et $d_{0}=1$). Dans ce cas, on pose $\Pi^G(\varphi)=\emptyset$ et, pour unifier les notations, $S(\varphi)/S(\varphi)^0=\{1\}$, ${\cal E}^{G}(\varphi)=\emptyset$. Supposons qu'\`a $\hat{L}$ corresponde un L\'evi $L$ de $G$. On a 
$$L=GL(d_{1})\times ...\times GL(d_{t})\times G_{0},$$
o\`u $G_{0}$ est un groupe de m\^eme type que $G$. Pour tout $j=1,...,t$, notons $\pi(\varphi_{j})$ la repr\'esentation admissible irr\'eductible de $GL(d_{j},F)$ d\'etermin\'ee par $\varphi_{j}$ par la correspondance de Langlands, prouv\'ee par Harris-Taylor et Henniart. Admettons la conjecture pour les repr\'esentations temp\'er\'ees. A $\varphi_{0}$ correspond un $L$-paquet $\Pi^{G_{0}}(\varphi_{0})$ de repr\'esentations temp\'er\'ees de $G_{0}(F)$. Soit $P$ le sous-groupe parabolique de $G$, de composante de L\'evi $L$, pour lequel la suite $(b_{1},...,b_{t})$ est positive dans un sens usuel. Pour $\sigma_{0}\in \Pi^{G_{0}}(\varphi_{0})$, notons $\sigma$ le quotient de Langlands de la repr\'esentation induite
$$Ind_{P}^{G}((\pi(\varphi_{1})\otimes\vert det\vert _{F}^{b_{1}})\otimes...\otimes(\pi(\varphi_{t})\otimes \vert det\vert _{F}^{b_{t}})\otimes \sigma_{0}),$$
c'est-\`a-dire son unique quotient irr\'eductible. On note $\Pi^G(\varphi)$ l'ensemble de ces repr\'esentations $\sigma$ quand $\sigma_{0}$ d\'ecrit $\Pi^{G_{0}}(\varphi_{0})$. Rappelons que le paquet $\Pi^{G_{0}}(\varphi_{0})$ est param\'etr\'e (conjecturalement) par un ensemble ${\cal E}^{G_{0}}(\varphi_{0})$ de caract\`eres d'un groupe $S(\varphi_{0})/S(\varphi_{0})^0$ isomorphe \`a un produit de facteurs ${\mathbb Z}/2{\mathbb Z}$. On a rappel\'e ce param\'etrage en [W1] 4.2, on le note $\epsilon\mapsto \sigma(\varphi_{0},\epsilon)$. On pose $S(\varphi)/S(\varphi)^0=S(\varphi_{0})/S(\varphi_{0})^0$,  ${\cal E}^{G}(\varphi)={\cal E}^{G_{0}}(\varphi_{0})$ et, pour un \'el\'ement $\epsilon \in {\cal E}^{G}(\varphi)$, on note $\sigma(\varphi,\epsilon)$ la repr\'esentation d\'eduite de $\sigma_{0}=\sigma(\varphi_{0},\epsilon)$.

Les repr\'esentations de $G(F)$ se param\`etrent de fa\c{c}on similaire dans le cas o\`u $d$ est pair. Pour tout entier naturel $N$, on fixe une forme bilin\'eaire sym\'etrique non d\'eg\'en\'er\'ee sur ${\mathbb C}^N$ et on note $O(N,{\mathbb C})$, resp. $SO(N,{\mathbb C})$, son groupe orthogonal, resp. sp\'ecial orthogonal. Consid\'erons un homomorphisme $\varphi:W_{DF}\to O(N,{\mathbb C})$ v\'erifiant les m\^emes conditions que pr\'ec\'edemment. Par composition avec le d\'eterminant, on obtient un caract\`ere quadratique de $W_{DF}$, forc\'ement trivial sur $SL(2,{\mathbb C})$, donc un caract\`ere quadratique de $W_{F}$. Par la th\'eorie du corps de classes, il d\'etermine un \'el\'ement $\delta(\varphi)\in F^{\times}/F^{\times,2}$. D\'efinissons un autre \'el\'ement de ce groupe par $\delta(q)=(-1)^{d/2}det(q)$. On note $\Phi(G)$ l'ensemble des classes de conjugaison par $SO(d,{\mathbb C})$ d'homomorphismes $\varphi:W_{DF}\to O(d,{\mathbb C})$ v\'erifiant les conditions pr\'ec\'edentes et tels que $\delta(\varphi)=\delta(q)$. Pour $\varphi\in \Phi(G)$, on d\'efinit comme ci-dessus le $L$-paquet $\Pi^G(\varphi)$: ou bien il est vide, ou bien c'est l'ensemble des quotients de Langlands issus de repr\'esentations $\pi(\varphi_{j})\otimes \vert det\vert _{F}^{b_{j}}$ de groupes $GL(d_{j},F)$ et des \'el\'ements $\sigma_{0}$ d'un $L$-paquet $\Pi^{G_{0}}(\varphi_{0})$. Le $L$-paquet est param\'etr\'e par un ensemble ${\cal E}^G(\varphi)={\cal E}^{G_{0}}(\varphi_{0})$ de caract\`eres d'un groupe $S(\varphi)/S(\varphi)^0=S(\varphi_{0})/S(\varphi_{0})^0$ produit de facteurs ${\mathbb Z}/2{\mathbb Z}$. Signalons toutefois que, pour associer un L\'evi $L$ de $G$ \`a un L\'evi $\hat{L}$ de $O(d,{\mathbb C})$, il faut consid\'erer $O(d,{\mathbb C})$ comme le $L$-groupe de $G$, ce qui sous-entend que l'on fixe des donn\'ees suppl\'ementaires cach\'ees (pr\'ecis\'ement, on identifie un sous-tore maximal de $SO(d,{\mathbb C})$ au $L$-groupe d'un sous-tore maximal de $G$).

  On admet la validit\'e des conjectures pour les repr\'esentations temp\'er\'ees, telles qu'on les a formul\'ees en [W1] 4.2 et 4.3, compl\'et\'ees comme en 2.1 ci-dessous dans le cas $d$ impair. Dans [W1]  4.8, on a pr\'ecis\'e les param\'etrages des $L$-paquets temp\'er\'es et ce sont ces param\'etrages que nous utilisons. Les param\'etrages des $L$-paquets de $G(F)$  ne d\'ependent d'aucune donn\'ee auxiliaire dans le cas o\`u $d$ est impair.  Ils d\'ependent par contre de  l'\'el\'ement $\nu_{0}$ d\'efini ci-dessus dans le cas o\`u $d$ est pair. Celui-ci sert \`a normaliser les facteurs de transfert.  

Soit $\varphi\in \Phi(G)$. Si $G$ est  quasi-d\'eploy\'e,  on dit que $\varphi$ est g\'en\'erique s'il existe $\sigma\in \Pi^{G}(\varphi)$   et un type de mod\`eles de Whittaker de sorte que $\sigma$ admette un mod\`ele de ce type. Rappelons que si $d$ est impair, il n'y a qu'un seul type de mod\`eles de Whittaker tandis que, si $d$ est pair,  il y a plusieurs types de tels mod\`eles, autant que  d'orbites unipotentes r\'eguli\`eres de $G(F)$. Si $G$ n'est pas quasi-d\'eploy\'e, on introduit un espace quadratique $(\underline{V},\underline{q})$ de m\^emes dimension et d\'eterminant que $(V,q)$, mais dont le groupe sp\'ecial orthogonal $\underline{G}$ est quasi-d\'eploy\'e. A $\varphi$ est associ\'e un $L$-paquet $\Pi^{\underline{G}}(\varphi)$ de repr\'esentations de $\underline{G}(F)$. On dit que $\varphi$ est g\'en\'erique s'il existe un \'el\'ement $\underline{\sigma}$ de ce paquet et un type de mod\`eles de Whittaker de sorte que $\underline{\sigma}$ admette un mod\`ele de ce type. 

On peut \'evidemment remplacer $G$ par $G'$ dans les consid\'erations ci-dessus. Soient $\varphi\in \Phi(G)$ et $\varphi'\in \Phi(G')$. En suivant Gross et Prasad, on a d\'efini en [W1]  4.9 un signe $E(\varphi,\varphi')$ et des caract\`eres $\boldsymbol{\epsilon}$ de $S(\varphi)/S(\varphi)^0$ et $\boldsymbol{\epsilon}'$ de $S(\varphi')/S(\varphi')^0$. On pose
$$\mu(G,G')=\left\lbrace\begin{array}{cc}1,&\text{ si }G\text{ et } G' \text{ sont quasi-d\'eploy\'es,}\\ -1,&\text{ sinon.}\\ \end{array}\right.$$
On sait que si $E(\varphi,\varphi')=\mu(G,G')$, $\boldsymbol{\epsilon}$ appartient \`a ${\cal E}^G(\varphi)$ et $\boldsymbol{\epsilon}'$ appartient \`a ${\cal E}^{G'}(\varphi')$. Si au contraire $E(\varphi,\varphi')=-\mu(G,G')$, $\boldsymbol{\epsilon}$ n'appartient pas \`a ${\cal E}^G(\varphi)$ et $\boldsymbol{\epsilon}'$ n'appartient pas \`a ${\cal E}^{G'}(\varphi')$.

En plus des conjectures de classification, on doit admettre certains r\'esultats issus de la formule des traces locale tordue, afin d'assurer la validit\'e du th\'eor\`eme  de [W2].

Notre principal r\'esultat est le th\'eor\`eme suivant, qui sera d\'emontr\'e dans la section 3.

\ass{Th\'eor\`eme}{Soient $\varphi\in \Phi(G)$ et $\varphi'\in \Phi(G')$. On suppose $\varphi$ et $\varphi'$ g\'en\'eriques. Alors:

(i) toutes les repr\'esentations induites dont les \'el\'ements de $\Pi^G(\varphi)$ et de $\Pi^{G'}(\varphi')$ sont les quotients de Langlands sont irr\'eductibles;

(ii) si $E(\varphi,\varphi')=-\mu(G,G')$, on a $m(\sigma,\sigma')=0$ pour tous $\sigma\in \Pi^G(\varphi)$, $\sigma'\in \Pi^{G'}(\varphi')$;

(iii) si $E(\varphi,\varphi')=\mu(G,G')$, on a
$$m(\sigma(\varphi,\boldsymbol{\epsilon}),\sigma'(\varphi',\boldsymbol{\epsilon}'))=1$$
et $m(\sigma,\sigma')=0$ pour tous $\sigma\in \Pi^G(\varphi)$, $\sigma'\in \Pi^{G'}(\varphi')$ tels que $(\sigma,\sigma')\not=(\sigma(\varphi,\boldsymbol{\epsilon}), \sigma'(\varphi',\boldsymbol{\epsilon}'))$.}

 Pour $d$ impair, les conjectures de classification des repr\'esentations de $G(F)$ font  partie des r\'esultats annonc\'es par Arthur, au moins dans le cas d'un groupe quasi-d\'eploy\'e ([A1] th\'eor\`eme 30.1). Dans le cas $d$ pair, Arthur annonce un r\'esultat un peu plus faible, o\`u on ne distingue pas deux repr\'esentations conjugu\'ees par le groupe orthogonal tout entier, cf. [W1] 4.4. En admettant seulement ces conjectures plus faibles, on peut \'enoncer un th\'eor\`eme similaire \`a celui ci-dessus, o\`u on regroupe les repr\'esentations conjugu\'ees par le groupe orthogonal. 

Il y a trois ingr\'edients dans la preuve du th\'eor\`eme. En pr\'ecisant un raisonnement d\^u \`a Gan, Gross et Prasad, on prouve dans la premi\`ere section que les multiplicit\'es sont compatibles \`a l'induction, sous des hypoth\`eses convenables (proposition 1.3). Cela ram\`ene les assertions (ii) et (iii) du th\'eor\`eme au cas temp\'er\'e, pourvu que toutes les induites intervenant soient irr\'eductibles, c'est-\`a-dire pourvu que l'assertion (i) soit v\'erifi\'ee. Dans la deuxi\`eme section, on \'etablit un crit\`ere g\'en\'eral d'irr\'eductibilit\'e pour ces induites (th\'eor\`eme 2.13). Signalons qu'il est \'egalement vrai pour les groupes symplectiques. La cons\'equence  de ce crit\`ere est que, parmi ces induites, il y en a  qui sont "les plus r\'eductibles" qu'il est possible. C'est-\`a-dire que si l'une des induites est r\'eductible, celles-ci le sont aussi. Ce sont celles pour lesquelles la repr\'esentation que l'on induit admet un mod\`ele de Whittaker (il y en a plusieurs s'il y a plusieurs types de tels mod\`eles). Pour obtenir le (i) du th\'eor\`eme, il reste \`a utiliser une conjecture de Shahidi d\'emontr\'ee par Muic qui affirme que si un  quotient de Langlands est g\'en\'erique (ce qui fait partie de nos hypoth\`eses: les paquets sont g\'en\'eriques), alors l'induite dont il est quotient est irr\'eductible.

\section{Induction et multiplicit\'es}

\bigskip

\subsection{Notations}

Selon le cas, on d\'esigne une repr\'esentation d'un groupe soit par la repr\'esentation elle-m\^eme, disons $\pi$, soit par un espace $E_{\pi}$ dans lequel elle se r\'ealise. On introduira parfois une repr\'esentation $\sigma$ d'un groupe sp\'ecial orthogonal sans pr\'eciser au d\'epart de quel groupe il s'agit. Dans ce cas, on notera $d_{\sigma}$ la dimension de l'espace quadratique $(V,q)$ du groupe sp\'ecial orthogonal duquel $\sigma$ est une repr\'esentation. De m\^eme, on introduira parfois une repr\'esentation $\pi$ d'un groupe lin\'eaire  (on entend par l\`a un groupe $GL(d,F)$) sans pr\'eciser ce groupe. On notera $d_{\pi}$ l'entier tel que $\pi$ est une repr\'esentation de $GL(d_{\pi},F)$.  On note $\check{\sigma}$ la contragr\'ediente d'une repr\'esentation $\sigma$.

Soient $(V,q)$ un espace quadratique de groupe sp\'ecial orthogonal $G$ et $\sigma$ une repr\'esentation lisse de $G(F)$. La notation
$$\sigma=\pi_{1}\times...\times \pi_{t}\times \sigma_{0},$$
ou encore
$$\sigma=(\times_{i=1,...,t}\pi_{i})\times \sigma_{0}$$
signifie ce qui suit.  On fixe une d\'ecomposition
$$V=X_{1}\oplus...\oplus X_{t}\oplus V_{0}\oplus Y_{t}\oplus...\oplus Y_{1}.$$
 On suppose que les espaces $X_{i}$ et $Y_{i}$ sont  totalement isotropes et non nuls. En posant  $V_{i}=X_{i}\oplus Y_{i}$ pour $i=1,...,t$, on suppose que les espaces $V_{i}$, pour $i=0,...,t$ sont orthogonaux.  On note $d_{i}$ la dimension de $X_{i}$ (ou encore celle de $Y_{i}$) et $d_{0}$ celle de $V_{0}$. On note $P$ le sous-groupe parabolique de $G$ qui conserve le drapeau
$$X_{1}\subset X_{1}\oplus X_{2}\subset...\subset X_{1}\oplus...\oplus X_{t},$$
$U$ son radical unipotent et $M$ sa composante de L\'evi form\'ee des \'el\'ements qui conservent  tous les espaces $X_{i}$ et $Y_{i}$. On note $G_{0}$ le groupe sp\'ecial orthogonal de $V_{0}$ et on fixe une base de $X_{i}$ qui permet de noter   $GL(d_{i})$ le groupe lin\'eaire de l'espace $X_{i}$. On a l'\'egalit\'e
$$M=GL(d_{1})\times...\times GL(d_{t})\times G_{0}.$$
Le terme $\sigma_{0}$ est une repr\'esentation lisse de $G_{0}(F)$ et, pour tout $i=1,...,t$, $\pi_{i}$ est une repr\'esentation lisse de $GL(d_{i},F)$. Alors $\sigma$ est la repr\'esentation
$$Ind_{P}^G(\pi_{1} \otimes...\otimes \pi_{t}\otimes   \sigma_{0})$$
de $G(F)$.

Dans le cas o\`u la dimension $d$ de $V$ est paire, cette notation est ambig\" ue car cette repr\'esentation induite peut d\'ependre du choix de la d\'ecomposition de $V$. En pratique, l'important sera que ces d\'ecompositions soient choisies de fa\c{c}on coh\'erentes.

On utilise une notation similaire
$$\pi=\pi_{1}\times...\times\pi_{t}=\times_{i=1,...,t}\pi_{i}$$
pour des repr\'esentations de groupes lin\'eaires.

Soit $\pi$ une repr\'esentation admissible irr\'eductible d'un groupe lin\'eaire. Elle poss\`ede un caract\`ere central $\omega_{\pi}$, qui est un caract\`ere de $F^{\times}$. Il existe un unique r\'eel $e_{\pi}$ tel que $\vert \omega_{\pi}(\lambda)\vert =\vert \lambda\vert _{F}^{e_{\pi}d_{\pi}}$ pour tout $\lambda\in F^{\times}$. On appelle $e_{\pi}$ l'exposant de $\pi$. Soit $b\in {\mathbb R}$. On note $\pi\vert .\vert _{F}^b$ la repr\'esentation $\pi\otimes \vert det\vert _{F}^b$. Si $\pi$ est unitaire, on a $e_{\pi\vert .\vert _{F}^b}=b$. On appelle support cuspidal de $\pi$ l'ensemble avec multiplicit\'es $\{\rho_{1},..,\rho_{k}\}$ form\'e de repr\'esentations irr\'eductibles et cuspidales tel que $\pi$ soit un sous-quotient de l'induite $\times_{i=1,...,k}\rho_{i}$.  On utilise ici, \`a la suite de Zelevinsky, la notion d'ensembles avec multiplicit\'es. Ils sont toujours finis et peuvent \^etre consid\'er\'es comme des familles finies prises \`a l'ordre pr\`es.

Soient $e$ et $f$ deux r\'eels tels que $e-f\in {\mathbb N}$. Par abus de notation, on d\'esigne par $[e,f]$ l'ensemble des r\'eels $x$ tels que $e\geq x\geq f$ et $e-x\in {\mathbb N}$. On ordonne cet ensemble par l'ordre  oppos\'e  de l'ordre ordinaire. Soit $\rho$ une  repr\'esentation admissible irr\'eductible et cuspidale d'un groupe lin\'eaire. On note $<e,f>_{\rho}$ l'unique sous-module irr\'eductible de la repr\'esentation induite $\times_{x\in [e,f]}\pi\vert .\vert _{F}^x$, le produit \'etant pris dans l'ordre que l'on vient de pr\'eciser. Dans le cas particulier o\`u $\rho$ est unitaire et o\`u $e=(a-1)/2$ et $f=(1-a)/2$ pour un entier $a\geq1$, on note aussi $St(\rho,a)$ cette repr\'esentation (une repr\'esentation de Steinberg g\'en\'eralis\'ee). En g\'en\'eral, $<e,f>_{\rho}$ est une telle repr\'esentation de Steinberg g\'en\'eralis\'ee tordue par le caract\`ere $\vert det\vert _{F}^{(e+f)/2}$. Pour simplifier la notation, on utilisera aussi la notation $[e,f]$ dans le cas o\`u $f=e+1$. Par convention, $[e,f]=\emptyset$ dans ce cas.

Soit $N\geq1$ un entier.  Pour tout $k=1,...,N$ notons $Q_{N-k,k}$ le sous-groupe parabolique triangulaire sup\'erieur de $GL_{N}$ de composante de L\'evi 
$$L_{N-k,k}=GL_{N-k}\times GL(1)\times...\times GL(1),$$
avec $k$ termes $GL(1)$.  Notons $U_{N-k,k}$ son radical unipotent et  $P_{N-k,k}$ le sous-groupe de $Q_{N-k,k}$ form\'e des \'el\'ements dont les projections dans les $k$ derniers facteurs $GL(1)$ de $L_{N-k,k}$ valent $1$. En particulier,  $P_{N-1,1}$ est le groupe appel\'e "mirabolique". On fixe un caract\`ere continu et non trivial $\psi_{F}$ de $F$ et on d\'efinit un caract\`ere $\psi_{k}$  de $P_{N-k,k}(F)$  par
$$\psi_{k}(p)=\psi_{F}(\sum_{i=N-k+1,...,N-1}p_{i,i+1}),$$
o\`u les $p_{i,i+1}$ sont les coefficients  de $p$.

\bigskip

\subsection{D\'efinition des multiplicit\'es}

On fixe pour la suite de la section  deux espaces quadratiques  $(V,q)$ et $(V',q')$ comme dans l'introduction.  On pose $r=(d-d'-1)/2$ et on introduit un espace $Z$ de dimension $2r+1$ sur $F$, muni d'une base $(v_{k})_{k=-r,...,r}$. On d\'efinit une forme bilin\'eaire sym\'etrique $q_{Z}$ sur $Z$ par
$$q_{Z}(\sum_{k=-r,...,r}\lambda_{k}v_{k},\sum_{k=-r,...,r}\lambda'_{k}v_{k})=2(-1)^d\nu_{0}\lambda_{0}\lambda'_{0}+\sum_{k=1,...,r}(\lambda_{-k}\lambda'_{k}+\lambda_{k}\lambda'_{-k}).$$
On fixe un isomorphisme de $(V,q)$ sur la somme orhogonale $(V',q')\oplus (Z,q_{Z})$. On note $Q$ le sous-groupe parabolique de $G$ qui conserve le drapeau de sous-espaces
$$Fv_{r}\subset Fv_{r}\oplus Fv_{r-1}\subset...\subset Fv_{r}\oplus...\oplus Fv_{1}.$$
On note $N$ son radical unipotent. On  d\'efinit un caract\`ere $\psi_{N}$ de $N(F)$ par la formule
$$\psi_{N}(n)=\psi_{F}(\sum_{k=0,...,r-1}q(nv_{k},v_{-k-1})).$$
Le groupe $G'$ est inclus dans $Q$ et $\psi_{N}$ est invariant par conjugaison par $G'(F)$. Soient $\sigma$ et $\sigma'$ des repr\'esentations lisses de $G(F)$ et $G'(F)$. On note $Hom_{G'(F),\psi_{F}}(\sigma,\sigma')$ l'espace des applications lin\'eaires $l:E_{\sigma'}\to E_{\sigma}$ telles que 
$$l(\sigma(ng')e)=\psi_{N}(n)\sigma(g')l(e)$$
pour tous $n\in N(F)$, $g'\in G'(F)$, $e\in E_{\sigma}$. On note $m(\sigma,\sigma')$ ou $m(\sigma',\sigma)$ la dimension de cet espace $Hom_{G'(F),\psi_{F}}(\sigma',\sigma)$. Cette dimension ne d\'epend pas des choix effectu\'es (l'espace $Z$, sa base, le caract\`ere $\psi_{F}$). Si $\sigma$ et $\sigma'$ sont irr\'eductibles, on a $m(\sigma,\sigma')\leq1$. 

\bigskip

\subsection{Induction}

 Consid\'erons une repr\'esentation induite
 $$\sigma=\pi_{1}\vert .\vert _{F}^{b_{1}}\times...\times \pi_{t}\vert .\vert _{F}^{b_{t}}\times \sigma_{0}$$
 de $G(F)$, o\`u
 
 $\bullet$ pour $i=1,...,t$, $\pi_{i}$ est une repr\'esentation admissible irr\'eductible et temp\'er\'ee d'un groupe lin\'eaire $GL(N_{i},F)$;
 
 $\bullet$ $\sigma_{0}$ est une repr\'esentation admissible irr\'eductible et temp\'er\'ee de $G_{0}(F)$, o\`u $G_{0}$ est le groupe sp\'ecial orthogonal d'un certain sous-espace $V_{0}$ de $V$, cf. 1.1;
 
 $\bullet$ les $b_{i}$ sont des nombres r\'eels tels que $b_{1}\geq b_{2}\geq...\geq b_{t}\geq0$.

 On consid\`ere une repr\'esentation
  $$\sigma'=\pi'_{1}\vert .\vert _{F}^{b'_{1}}\times...\times \pi'_{t'}\vert .\vert _{F}^{b'_{t'}}\times \sigma'_{0}$$
 de $G'(F)$ v\'erifiant des hypoth\`eses similaires.
 
  On a d\'efini la multiplicit\'e $m(\sigma,\check{\sigma}')$ dans le paragraphe pr\'ec\'edent . On remarque qu'\`a conjugaison pr\`es, on peut supposer que le plus petit des deux espaces $V_{0}$ et $V'_{0}$ est inclus dans le plus grand et qu'alors, ces deux espaces verifient les m\^emes hypoth\`eses que $V$ et $V'$, avec la m\^eme constante $\nu_{0}$. On peut donc d\'efinir la multiplicit\'e $m(\sigma_{0},\sigma'_{0})$. La section est consacr\'ee \`a la preuve de la proposition suivante.

\ass{Proposition}{Sous ces hypoth\`eses, on a l'\'egalit\'e
$$m(\sigma,\check{\sigma}')=m(\sigma_{0},\sigma'_{0}).$$} 

{\bf Remarque.} On pourrait remplacer $\sigma_{0}$ ou $\sigma'_{0}$ par leurs contragr\'edientes dans le second membre de cette \'egalit\'e, cf. [W3] 7.1(1).

 {\bf Convention.} Dans la suite interviendront des in\'egalit\'es portant sur $b_{1}$ et $b'_{1}$. Dans le cas o\`u $t=0$, auquel cas $b_{1}$ n'est pas d\'efini, on consid\`ere que $b_{1}=0$ dans ces in\'egalit\'es. De m\^eme pour $b'_{1}$.

 \bigskip

\subsection{Une premi\`ere in\'egalit\'e}

  On suppose  $ d=d'+1$.  On consid\`ere les donn\'ees du paragraphe pr\'ec\'edent en ce qui concerne le groupe $G'$.  On consid\`ere une repr\'esentation induite
  $$\sigma=\pi\vert .\vert _{F}^b\times \sigma_{0}$$
  de $G(F)$,  o\`u $\pi$ est une repr\'esentation admissible irr\'eductible d'un groupe lin\'eaire $GL(N,F)$, unitaire et de la s\'erie discr\`ete, $b$ est un r\'eel et $\sigma_{0}$ est une repr\'esentation admissible, pas forc\'ement irr\'eductible, d'un groupe $G_{0}(F)$.

   \ass{Lemme}{Supposons $b\geq b'_{1}$. Alors $m(\sigma,\check{\sigma}')\leq m(\sigma',\check{\sigma}_{0})$.}
  
 Conform\'ement \`a notre convention, l'in\'egalit\'e $b\geq b'_{1}$ signifie $b\geq0$ si $t'=0$. La d\'emonstration ci-dessous est inspir\'ee de celle du th\'eor\`eme 15.1 de [GGP].

  Preuve.   R\'ealisons $\sigma$ comme en 1.1. On fixe une d\'ecomposition 
  $$V=X\oplus V_{0}\oplus Y$$
  comme dans ce paragraphe, dont on d\'eduit des groupes $P$ et 
  $$M=GL(N)\times G_{0}.$$
  On a
  $$\sigma=Ind_{P}^G(\pi\vert .\vert _{F}^b\otimes \sigma_{0}).$$
  A tout $g\in P(F)\backslash G(F)$, associons le sous-espace $g^{-1}(X)$ de $V$. Il est totalement isotrope et de dimension $N$. Notons ${\cal U}$ le sous-ensemble des $g\in P(F)\backslash G(F)$ tels que $dim(g^{-1}(X)\cap V')=N-1$ et ${\cal X}$ celui des $g$ tels que $dim(g^{-1}(X)\cap V')=N$. Alors ${\cal U}$ est un ouvert de compl\'ementaire le ferm\'e ${\cal X}$.  On v\'erifie que ${\cal U}$   est form\'e d'une seule orbite sous l'action de $G'(F)$. L'ensemble ${\cal X}$ est en g\'en\'eral form\'e d'une seule orbite, sauf dans le cas o\`u $d$ est impair et $N=d'/2$, o\`u il y en a deux.

  Notons $E_{\sigma}$ l'espace  de la repr\'esentation induite $\sigma$ form\'e de fonctions sur $G(F)$ \`a valeurs dans $E_{\pi}\otimes_{{\mathbb C}}E_{\sigma_{0}}$ v\'erifiant la condition habituelle de transformation \`a gauche par $P(F)$. Notons $E_{\sigma,{\cal U}}$ le sous-espace form\'e des fonctions \`a support dans ${\cal U}$ et $E_{\sigma,{\cal X}}$ le quotient $E_{\sigma}/E_{\sigma,{\cal U}}$. Ces espaces sont stables par $G'(F)$ et on a une suite exacte
  $$ 0\to Hom_{G'(F)}(E_{\sigma,{\cal X}},E_{\check{\sigma}'})\to Hom_{G'(F)}(E_{\sigma},E_{\check{\sigma}'})\to Hom_{G'(F)}(E_{\sigma,{\cal U}},E_{\check{\sigma}'}).\eqno(1)$$
  
  Montrons que
  $$ Hom_{G'(F)}(E_{\sigma,{\cal X}},E_{\check{\sigma}'})=\{0\}.\eqno(2)$$
  
 Supposons que ${\cal X}$ soit form\'e d'une seule orbite pour l'action de $G'(F)$. Quitte \`a effectuer une conjugaison par un \'el\'ement de $G(F)$, on peut supposer $X\subset V'$ et $Y\subset V'$. En posant $V'_{0}=V_{0}\cap V'$, on a 
  $$V'=X\oplus V'_{0}\oplus Y$$
  Notons $P'$  et $M'$ le sous-groupe parabolique de $G'$ et sa composante de L\'evi associ\'es \`a cette d\'ecomposition. On a $P'=P\cap G'$ et l'application $P'(F)g'\mapsto P(F)g'$ est un hom\'eomorphisme de $P'(F)\backslash G'(F)$ sur ${\cal X}$. On en d\'eduit que la repr\'esentation de $G'(F)$ dans $E_{\sigma,{\cal X}}$ est une repr\'esentation induite \`a partir d'une repr\'esentation de $M'(F)$. Pour d\'ecrire celle-ci, on doit remarquer que, pour $m'=(x,g'_{0})\in GL(N,F)\times G'_{0}(F)$, on a $\delta_{P}(m')^{1/2}=\vert det(x)\vert _{F}^{(d_{0}+N-1)/2}$ tandis que $\delta_{P'}(m')^{1/2}=\vert det(x)\vert _{F}^{(d_{0}+N-2)/2}$. On en d\'eduit qu'en notant $\sigma_{0\vert G'}$ la restriction de $\sigma_{0}$ \`a $G'(F)$, la repr\'esentation de $G'(F)$ dans $E_{\sigma,{\cal X}}$ est \'egale \`a
  $$Ind_{P'}^{G'}((\pi\otimes \vert det\vert_{F} ^{b+1/2})\otimes \sigma_{0\vert G'}),$$
  c'est-\`a-dire \`a
  $$\pi\vert .\vert_{F} ^{b+1/2}\times \sigma_{0\vert G'},$$
  avec les notations de 1.1.  Par le th\'eor\`eme de seconde r\'eciprocit\'e de Bernstein, on a
  $$Hom_{G'(F)}(E_{\sigma,{\cal X}},E_{\check{\sigma}'})=Hom_{M'(F)}((\pi\otimes \vert det\vert ^{b+1/2})\otimes \sigma_{0\vert G'},(\check{\sigma}')_{\bar{P}'}),$$
  le  terme $(\check{\sigma}')_{\bar{P}'}$ d\'esignant le module de Jacquet normalis\'e. La repr\'esentation $\pi$ \'etant irr\'eductible et unitaire,  l'exposant de $\pi\vert .\vert ^{b+1/2}$ est \'egal \`a $b+1/2$. La repr\'esentation $(\check{\sigma}')_{\bar{P}'}$ est de longueur finie. Pour d\'emontrer (2), il suffit de prouver que, pour tout sous-quotient irr\'eductible cette repr\'esentation, de la forme $\pi_{\sharp}\otimes \sigma'_{\sharp}$, l'exposant $e_{\sharp}$ de $\pi_{\sharp}$ ne peut pas \^etre \'egal \`a $b+1/2$. Par le calcul habituel du module de Jacquet d'une induite, on sait qu'une telle repr\'esentation $\pi_{\sharp}$ intervient dans une induite 
  $$\pi_{1,\flat}\times ...\times \pi_{t',\flat}\times \pi_{0,\sharp}\times \pi_{t',\natural}\times...\times \pi_{1,\natural}$$
  o\`u:
  
  $\bullet$ pour $i=1,...,t'$, il existe une repr\'esentation irr\'eductible $\tau_{i,\flat}$ d'un groupe lin\'eaire de sorte que $\pi_{i,\flat}\otimes \tau_{i,\flat}$ intervienne comme sous-quotient d'un module de Jacquet de $\check{\pi}'_{i}\vert .\vert ^{-b'_{i}}$ relatif \`a un sous-groupe parabolique triangulaire inf\'erieur (parce qu'on calcule le module de Jacquet de $\sigma'$ relativement \`a $\bar{P}'$);

  $\bullet$ pour $i=1,...,t'$, il existe une repr\'esentation irr\'eductible $\tau_{i,\natural}$ d'un groupe lin\'eaire de sorte que $\pi_{i,\natural}\otimes \tau_{i,\natural}$ intervienne comme sous-quotient d'un module de Jacquet de $\pi'_{i}\vert .\vert ^{b'_{i}}$ relatif \`a un sous-groupe parabolique triangulaire inf\'erieur;

$\bullet$ il existe une repr\'esentation irr\'eductible $\sigma'_{0,\flat}$ d'un groupe sp\'ecial orthogonal de sorte que $\pi_{0,\sharp}\times \sigma'_{0,\sharp}$ intervienne dans un module de Jacquet de $\check{\sigma}'_{0}$ relativement \`a un sous-groupe parabolique "triangulaire inf\'erieur".

  Notons $N_{i,\flat}$, $N_{i,\natural}$ et $N_{0,\sharp}$ les rangs des groupes lin\'eaires dont  $\pi_{i,\flat}$, $\pi_{i,\natural}$ et $\pi_{0,\sharp}$ sont des repr\'esentations et $e_{i,\flat}$, $e_{i,\natural}$ et $e_{0,\sharp}$ les  exposants de ces repr\'esentations. Alors  on a les \'egalit\'es
  $$N=(\sum_{i=1,...,t'}N_{i,\flat })+N_{0,\sharp}+(\sum_{i=1,...,t'}N_{i,\natural}),$$
$$Ne_{\sharp}=(\sum_{i=1,...,t'}N_{i,\flat }(e_{i,\flat}-b'_{i}))+N_{0,\sharp}e_{0,\sharp}+(\sum_{i=1,...,t'}N_{i,\natural}(e_{i,\natural}+b'_{i})).$$
Parce que les repr\'esentations $\pi'_{i}$ et $\check{\sigma}'_{0}$ sont temp\'er\'ees et qu'on consid\`ere leurs modules de Jacquet relatifs \`a des sous-groupes paraboliques inf\'erieurs, les nombres $e_{i,\flat}$, $e_{i,\natural}$ et $e_{0,\sharp}$ sont n\'egatifs ou nuls. On en d\'eduit l'in\'egalit\'e $e_{\sharp}\leq b'_{1}$. Sous l'hypoth\`ese de l'\'enonc\'e, on a donc $e_{\sharp}\leq b$, a fortiori $e_{\sharp}\not=b+1/2$, et cela prouve (2) dans le cas o\`u ${\cal X}$ est form\'e d'une seule orbite. Quand il y a deux orbites, $E_{\sigma,{\cal X}}$ est somme de deux sous-repr\'esentations auxquelles on peut appliquer le m\^eme raisonnement. D'o\`u (2).

Etudions la repr\'esentation $E_{\sigma,{\cal U}}$ de $G'(F)$. Quitte \`a effectuer une conjugaison par un \'el\'ement de $G(F)$, on se ram\`ene \`a la situation suivante. On a une d\'ecomposition
$$V'=X'\oplus D'\oplus V_{0}\oplus Y',$$
o\`u $X'$ et $Y'$ sont totalement isotropes de dimension $N-1$ et les espaces $X'\oplus Y'$, $D'$ et $V_{0}$ sont orthogonaux. L'espace $D'$ est une droite port\'ee par un vecteur $v'_{0}$ tel que $q'(v'_{0},v'_{0})=2(-1)^{d'}\nu_{0}$. Rappelons que $V=V'\oplus D$, o\`u $D$ est une droite port\'ee par un vecteur $v_{0}$ tel que $q(v_{0};v_{0})=2(-1)^d\nu_{0}$. On a les \'egalit\'es $X=X'\oplus F(v_{0}+v'_{0})$, $Y=Y'\oplus F(v_{0}-v'_{0})$. Posons $V'_{\infty}=D'\oplus V_{0}$, notons $G'_{\infty}$ son groupe sp\'ecial orthogonal, $P'_{\infty}$ le sous-groupe parabolique de $G'$ form\'e des \'el\'ements qui conservent $X'$, $P'_{\infty}=M'_{\infty}U'_{\infty}$ sa d\'ecomposition habituelle, o\`u $M'_{\infty}=GL(N-1)\times G'_{\infty}$. On a $G_{0}\subset G'_{\infty}$ et on v\'erifie que $P\cap G'=(GL(N-1)\times G_{0})U'_{\infty}$. La repr\'esentation $E_{\sigma,{\cal U}}$ de $G'(F)$ est une induite compacte \`a partir de ce groupe $P\cap G'$. Ici, il n'y a pas de d\'ecalage sur les modules: la restriction de $\delta_{P}$ \`a $P\cap G'$ est le module usuel de ce groupe. Par contre, la repr\'esentation que l'on induit n'est pas triviale sur  $U'_{\infty}(F)$.  Ce groupe $U'_{\infty}(F)$ admet une filtration \`a deux crans
$$0\to \wedge^2(X')\to U'_{\infty}(F)\to (V'_{\infty}\otimes_{F} X')\to 0.$$
Par projection orthogonale de $V'_{\infty}$ sur $D'$, on obtient une projection de $U'_{\infty}(F)$ sur $D'\otimes_{F}X'$, dont on note $U'_{\flat}(F)$ le noyau, puis une projection de $(GL(N-1,F)\times G_{0}(F))U'_{\infty}(F)$ sur le produit semi-direct $GL(N-1,F)\ltimes (D'\otimes_{F}X')$. En choisissant des bases convenables, on peut identifier  ce dernier groupe au sous-groupe mirabolique $P_{N-1,1}(F)$ de $GL(N,F)$. On a aussi une projection de $(GL(N-1,F)\times G_{0}(F))U'_{\infty}(F)$ sur $G_{0}(F)$, d'o\`u une projection
$$(GL(N-1,F)\times G_{0}(F))U'_{\infty}(F)\to P_{N-1,1}(F)\times G_{0}(F).$$
Notons $(\pi\vert .\vert _{F}^b)_{\vert P_{N-1,1}}$ la restriction de $\pi\vert .\vert _{F}^b$ en une repr\'esentation de $P_{N-1,1}(F)$. Par la projection ci-dessus, on peut consid\'erer $(\pi\vert .\vert _{F}^b)_{\vert P_{N-1,1}}\otimes \sigma_{0}$ comme une repr\'esentation de  $(GL(N-1,F)\times G_{0}(F))U'_{\infty}(F)$. On a alors
$$E_{\sigma,{\cal U}}=ind_{(GL(N-1)\times G_{0})U'_{\infty}}^{G'}((\pi\vert .\vert _{F}^b)_{\vert P_{N-1,1}}\otimes \sigma_{0}),$$
o\`u $ind$ d\'esigne l'induction \`a supports compacts. D'apr\`es [BZ] 3.5, la repr\'esentation $(\pi\vert .\vert _{F}^b)_{\vert P_{N-1,1}}$ poss\`ede une filtration
$$\{0\}=\tau_{N+1}\subset \tau_{N}\subset \tau_{N-1}\subset...\subset \tau_{1}=(\pi\vert .\vert _{F}^b)_{\vert P_{N-1,1}}.$$
Les quotients de cette filtration  v\'erifient
$$\tau_{k}/\tau_{k+1}\simeq ind_{P_{N-k,k}}^{P_{N-1,1}}(\Delta^k(\pi\vert .\vert _{F}^b)\otimes \psi_{k})$$
o\`u  $\Delta^k(\pi\vert .\vert _{F}^b)$ est la $k$-i\`eme d\'eriv\'ee de $\pi\vert .\vert _{F}^b$, cf. [BZ] 4.3. L'induction \`a supports compacts \'etant un foncteur exact, on en d\'eduit une filtration
$$\{0\}=\mu_{N+1}\subset \mu_{N}\subset \mu_{N-1}\subset...\subset \mu_{1}=E_{\sigma,{\cal U}}$$
dont les quotients v\'erifient
$$ \mu_{k}/\mu_{k+1}\simeq ind_{P_{N-k,k}G_{0}U'_{\flat}}^{G'}((\Delta^k(\pi\vert .\vert _{F}^b)\otimes \psi_{k}\otimes \sigma_{0}).\eqno(3)$$
Pour tout $k=1,...,N$, on a une suite exacte
$$ 1\to Hom_{G'(F)}(\mu_{k}/\mu_{k+1},\check{\sigma}')\to Hom_{G'(F)}(\mu_{k},\check{\sigma}')\to Hom_{G'(F)}(\mu_{k+1},\check{\sigma}'). \eqno(4)$$ 
Montrons que

(5) pour $k=1,...,N-1$, $Hom_{G'(F)}(\mu_{k}/\mu_{k+1},\check{\sigma}')=\{0\}$.

Introduisons le sous-groupe parabolique $Q'_{N-k}$ de $G'$ form\'e des \'el\'ements qui conservent le sous-espace engendr\'e par les $N-k$ premiers vecteurs de base de $X'$.   Sa composante de L\'evi $L'_{N-k}$ s'identifie \`a $GL(N-k)\times G'_{k}$, o\`u $G'_{k}$ est un groupe sp\'ecial orthogonal de m\^eme type que $G'$.   On a $P_{N-k,k}G_{0}U'_{\flat}\subset Q'_{N-k}$ et, dans la formule (3), la repr\'esentation que l'on induit est triviale  sur le radical unipotent de $Q'_{N-k}(F)$. Introduisons la repr\'esentation $\mu'_{k}$ de $G'_{k}(F)$ d\'efinie par
$$\mu'_{k}=ind_{(U_{N-k,k}G_{0}U'_{\flat})\cap G'_{k}}^{G'_{k}}(\psi_{k}\otimes \sigma_{0}).$$
Alors 
$$\mu_{k}/\mu_{k+1}\simeq Ind_{Q'_{N-k}}^{G'} (\Delta^k(\pi\vert .\vert _{F}^b)\otimes \mu'_{k})=\Delta^k(\pi\vert .\vert _{F}^b)\times \mu'_{k}.$$
Par le th\'eor\`eme de seconde adjonction de Bernstein, on a
$$Hom_{G'(F)}(\mu_{k}/\mu_{k+1},\sigma')=Hom_{L'_{N-k}(F)}(\Delta^k(\pi\vert .\vert _{F}^b)\otimes \mu'_{k}, \sigma'_{\bar{Q}'_{N-k}}).$$
Puisque $\pi$ est unitaire et de la s\'erie discr\`ete, elle est de la forme $St(\rho,a)$, pour une repr\'esentation irr\'eductible et cuspidale $\rho$ et un entier $a\geq1$. D'apr\`es [Z] proposition 9.6, $\Delta^k(\pi\vert .\vert _{F}^b)$ est nulle si $k$ n'est pas un multiple de $d_{\rho}$. Si $k=ld_{\rho}$ pour un entier $l\geq1$, on a $\Delta^k(\pi\vert .\vert _{F}^b)=<(a-1)/2+b,(1-a)/2+b+l>_{\rho}$. L'exposant de cette repr\'esentation est strictement sup\'erieur \`a $b$. La preuve de (2) montre que, pour tout sous-quotient irr\'eductible $\pi_{\sharp}\otimes \sigma'_{\sharp}$ de $(\check{\sigma}')_{\bar{Q}'_{N-k}}$, l'exposant de $\pi_{\sharp}$ est inf\'erieur ou \'egal \`a $b$. L'espace d'homomorphismes ci-dessus est donc nul, ce qui prouve (5).

En utilisant (4) et (5) successivement pour $k=1,...,N-1$, on obtient une injection
$$ Hom_{G'(F)}(E_{\sigma,{\cal U}},E_{\check{\sigma}'})\to Hom_{G'(F)}(\mu_{N},\check{\sigma}').\eqno(6)$$
Pour $k=N$, la repr\'esentation $\Delta^k(\pi\vert .\vert_{F} ^b)$ dispara\^{\i}t de la d\'efinition de $\mu_{N}$. La contragr\'ediente de l'induite \`a supports compacts d'une repr\'esentation admissible \'etant l'induite ordinaire de la contragr\'ediente, on a
$$Hom_{G'(F)}(\mu_{N},\check{\sigma}')=Hom_{G'(F)}(\sigma',Ind_{P_{0,N}G_{0}U'_{\flat}}^{G'}(\psi_{N}^{-1}\otimes \check{\sigma}_{0}))$$
$$=Hom_{P_{0,N}(F)G_{0}(F)U'_{\flat}(F)}(\sigma',\psi_{N}^{-1}\otimes \check{\sigma}_{0}).$$
 Remarquons que $P_{0,N}U'_{\flat}$ n'est autre que le radical unipotent  d'un sous-groupe parabolique de $G'$ contenu dans $P'_{\infty}$, de composante de L\'evi $GL(1)^{N-1}\times G'_{\infty}$. Alors le dernier espace ci-dessus n'est autre que celui not\'e $Hom_{G_{0}(F),\psi'_{F}}(\sigma',\check{\sigma}_{0})$ en 1.2, pour un choix convenable de caract\`ere $\psi'_{F}$. Sa dimension est $m(\sigma',\check{\sigma}_{0})$. En utilisant (1), (2) et (6), on obtient le lemme. $\square$

 \bigskip
 
 \subsection{Rappel d'un r\'esultat de Gan, Gross et Prasad}
 
 On conserve la situation du paragraphe pr\'ec\'edent. On suppose $\pi$ cuspidale et on fait l'hypoth\`ese suivante:
 
 (H) soit $\pi_{\sharp}$ une repr\'esentation irr\'eductible d'un groupe lin\'eaire telle qu'il existe une repr\'esentation irr\'eductible $\tau_{\sharp}$ d'un groupe lin\'eaire ou d'un groupe sp\'ecial orthogonal  de sorte que $\pi_{\sharp}\otimes \tau_{\sharp}$ intervienne comme sous-quotient d'un module de Jacquet de l'une des repr\'esentations $\pi'_{i}$, $\check{\pi}'_{i}$ ou $\check{\sigma}'_{0}$; soit de plus $\beta\in {\mathbb R}$; alors aucun \'el\'ement du support cuspidal de $\pi_{\sharp}$ n'est \'egal \`a $\pi\vert .\vert_{F}^{\beta}$.
 
 \ass{Lemme}{Sous ces hypoth\`eses, on a l'\'egalit\'e
 $$m(\sigma,\check{\sigma}')=m(\sigma',\check{\sigma}_{0}).$$}

 Cf. [GGP] th\'eor\`eme 15.1, dont on reproduit la d\'emonstration.
 
 Preuve. Bernstein a d\'ecompos\'e la cat\'egorie des repr\'esentations lisses de $G'(F)$ en somme directe de sous-cat\'egories  index\'ees par des classes d'inertie de supports cuspidaux. Notons $proj$ la projection sur la sous-cat\'egorie contenant $\check{\sigma}'$. On a    la suite exacte
 $$0\to proj(E_{\sigma,{\cal U}})\to proj(E_{\sigma})\to proj(E_{\sigma,{\cal X}})\to 0.$$
D'apr\`es la description de $E_{\sigma,{\cal X}}$ donn\'ee en 1.4, l'hypoth\`ese (H) implique que $proj(E_{\sigma,{\cal X}})=\{0\}$. Donc 
$$ proj(E_{\sigma,{\cal U}})= proj(E_{\sigma}).$$
Mais, pour toute repr\'esentation lisse $\tau$ de $G'(F)$, on a
$$Hom_{G'(F)}(\tau,\check{\sigma}')=Hom_{G'(F)}(proj(\tau),\check{\sigma}').$$
Donc
$$Hom_{G'(F)}(E_{\sigma},E_{\check{\sigma}'})=Hom_{G'(F)}(E_{\sigma,{\cal U}},E_{\check{\sigma}'}).$$
On d\'ecrit $E_{\sigma,{\cal U}}$ comme en 1.4. Puisque $\pi$ est cuspidale, la filtration de $(\pi\vert .\vert _{F}^b)_{P_{N-1,1}}$ n'a qu'un quotient non nul: on a $(\pi\vert .\vert _{F}^b)_{P_{N-1,1}}=\tau_{N}$. L'injection (6) de 1.4 devient une bijection et le r\'esultat s'ensuit. $\square$

\bigskip

\subsection{L'in\'egalit\'e $m(\sigma,\check{\sigma}')\leq m(\sigma_{0},\sigma'_{0})$}

 Dans la situation de 1.3, on va d\'emontrer l'in\'egalit\'e
 $$ m(\sigma,\check{\sigma}')\leq m(\sigma_{0},\sigma'_{0}).\eqno(1)$$
 Remarquons d'abord que, puisque  les repr\'esentations $\pi_{i}$ et $\pi'_{i}$ sont irr\'eductibles et temp\'er\'ees, ce sont des induites de s\'eries discr\`etes. En les \'ecrivant comme de telles induites, on obtient de nouvelles expressions de $\sigma$ et $\sigma'$, similaires aux expressions primitives, mais o\`u les nouvelles repr\'esentations $\pi_{i}$ et $\pi'_{i}$ sont de la s\'erie discr\`ete. En oubliant cette construction, on peut supposer que les $\pi_{i}$ et $\pi'_{i}$ sont de la s\'erie discr\`ete. 
 
 On utilisera plusieurs fois la construction suivante. Choisissons une repr\'esentation admissible irr\'eductible $\rho$ de $GL((d+1-d')/2,F)$, unitaire et cuspidale, v\'erifiant les analogues de l'hypoth\`ese (H) du paragraphe pr\'ec\'edent o\`u $\sigma'$ est remplac\'e par $\check{\sigma}'$ ou $\sigma$. Une telle repr\'esentation existe. Consid\'erons la repr\'esentation $\rho\times \sigma'$ d'un groupe  de m\^eme type que $G'$.   On  a
 $$\rho\times \sigma'=\rho\times\pi'_{1}\vert .\vert _{F}^{b'_{1}}\times...\times \pi'_{t'}\vert .\vert _{F}^{b'_{t'}}\times \sigma'_{0}.$$
 Mais les hypoth\`eses sur $\rho$ nous permettent de permuter les premiers facteurs:
$$\rho\times \sigma'=\pi'_{1}\vert .\vert _{F}^{b'_{1}}\times...\times \pi'_{t'}\vert .\vert _{F}^{b'_{t'}}\times \rho\times\sigma'_{0}.$$
C'est une expression du m\^eme type que celle de $\sigma'$: $t'$ est devenu $t'+1$, le r\'eel suppl\'ementaire $b'_{t'+1}$ est nul et on n'a pas chang\'e $\sigma'_{0}$. Remarquons que $d_{\rho\times \sigma'}=d+1$. On a l'\'egalit\'e

 $$m(\sigma,\check{\sigma}')=m(\rho\times \sigma',\check{\sigma}).\eqno(2)$$

En effet, cela r\'esulte de 1.5, o\`u l'on remplace respectivement $\pi$, $\sigma_{0}$ et $\sigma'$ par $\rho$, $\sigma'$ et $\sigma$.

 Revenons \`a notre probl\`eme et supposons d'abord que tous les $b_{i}$ et $b'_{i}$ sont nuls. Dans ce cas, prouvons:
 
 (3) on a l'\'egalit\'e $m(\sigma,\check{\sigma}')= m(\sigma_{0},\sigma'_{0})$.
 
 Il s'agit ici d'induites unitaires, donc $\sigma$ et $\sigma'$ sont semi-simples. D\'ecomposons $\sigma'$ en composantes irr\'eductibles: $\sigma'=\oplus_{j=1,...,k}\sigma'_{j}$. On a
 $$m(\sigma,\check{\sigma}')=\sum_{j=1,...,k}m(\sigma,\check{\sigma}'_{j}).$$
 Il s'agit ici de repr\'esentations temp\'er\'ees. En utilisant la proposition 5.7 et les lemmes 5.3 et 5.4 de [W3], on a, pour tout $j=1,...,k$,
 $$m(\sigma,\check{\sigma}'_{j})=m(\sigma_{0},\check{\sigma}'_{j}).$$
 Si $d_{\sigma_{0}}<d'$, on a
 $$\sum_{j=1,...,k}m(\sigma_{0},\check{\sigma}'_{j})=\sum_{j=1,...,k}m(\check{\sigma}'_{j},\sigma_{0})=m(\check{\sigma}',\sigma_{0}).$$
 Les m\^emes r\'esultats de [W3] entra\^{\i}nent que cette derni\`ere multiplicit\'e vaut $m(\check{\sigma}'_{0},\sigma_{0})$. Mais, pour les repr\'esentations irr\'eductibles, le changement d'une repr\'esentation en sa contragr\'ediente ne modifie pas la multiplicit\'e ([W3] 7.1(1)). Le dernier nombre est donc \'egal \`a $m(\sigma_{0},\sigma'_{0})$, d'o\`u (3) dans ce cas. Supposons maintenant $d_{\sigma_{0}}>d'$. On choisit $\rho$ comme ci-dessus, relativement aux repr\'esentations $\sigma_{0}$ et $\sigma'$. Remarquons que $\rho$ v\'erifie alors pour tout $j$ les m\^emes hypoth\`eses relativement aux repr\'esentations $\sigma_{0}$ et $\sigma'_{j}$. Donc, d'apr\`es (2), 
 $m(\sigma_{0},\check{\sigma}'_{j})=m(\rho\times\sigma'_{j},\check{\sigma}_{0})$, puis
 $$\sum_{j=1,...,k}m(\sigma_{0},\check{\sigma}'_{j})=m(\rho\times \sigma',\check{\sigma}_{0}).$$
 On peut maintenant appliquer les m\^emes r\'esultats de [W3]. Ils entra\^{\i}nent que ce dernier terme vaut $m(\sigma_{0},\sigma'_{0})$. Cela prouve (3).
 
 Revenons au cas g\'en\'eral. On d\'efinit un invariant
$$N(\sigma,\sigma')=(\sum_{i=1,...,t; b_{i}\not=0}N_{i})+(\sum_{i=1,...,t'; b'_{i}\not=0}N'_{i}).$$ 
Soit $N$ un entier naturel. On va d\'emontrer par r\'ecurrence sur $N$ que l'in\'egalit\'e (1)
  est v\'erifi\'ee pour toutes donn\'ees $G$, $G'$, $\sigma$, $\sigma'$ v\'erifiant les conditions requises et telles que $N(\sigma,\sigma')=N$.  Le cas $N=0$ est couvert par (3).  On  suppose maintenant $N>0$ et on fixe des donn\'ees avec $N(\sigma,\sigma')=N$.

$1^{er}$ cas. On suppose $d=d'+1$, $t\geq1$ et $b_{1}\geq b'_{1}$. Posons
$$\sigma_{1}= \pi_{2}\vert .\vert_{F} ^{b_{2}}\times...\times\pi_{t}\vert .\vert _{F}^{b_{t}}\times \sigma_{0}.$$
C'est une repr\'esentation d'un groupe  $G_{1}$ de m\^eme type que $G$ et on a $\sigma=\pi_{1}\vert .\vert _{F}^{b_{1}}\times\sigma_{1}$. La situation permet d'appliquer le lemme 1.4: on a
$$m(\sigma,\check{\sigma}')\leq m(\sigma', \check{\sigma}_{1}).$$
On applique l'hypoth\`ese de r\'ecurrence aux groupes $G'$ et $G_{1}$ et \`a leurs repr\'esentations $\sigma'$ et $\sigma_{1}$. C'est loisible puisque $N(\sigma', \sigma_{1})=N(\sigma,\sigma')- d_{\pi_{1}}<N$. On en d\'eduit $m(\sigma', \check{\sigma}_{1})\leq m(\sigma_{0},\sigma'_{0})$, puis (1).

$2^{\grave{e}me}$ cas. On suppose $t'\geq1$ et $b'_{1}\geq b_{1}$. Choisissons une repr\'esentation $\rho$ v\'erifiant les hypoth\`eses permettant d'appliquer (2). Gr\^ace \`a cette relation, on a $m(\sigma,\check{\sigma}')=m(\rho\times \sigma',\check{\sigma})$. Mais les repr\'esentations $\rho\times \sigma'$ et $\sigma$ v\'erifient les hypoth\`eses du premier cas. Donc $m(\rho\times\sigma',\check{\sigma})\leq m(\sigma_{0},\sigma'_{0})$.

 $3^{\grave{e}me}$ cas. On suppose $t\geq1$ et $b_{1}\geq b'_{1}$. On raisonne comme dans le deuxi\`eme cas, \`a ceci pr\`es que les repr\'esentations $\rho\times \sigma'$ et $\sigma$ v\'erifient maintenant les hypoth\`eses du deuxi\`eme cas.
 
 Pour $N>0$, on est forc\'ement dans l'un des deuxi\`eme ou troisi\`eme cas et cela prouve (1). $\square$
 
 \bigskip
 
 \subsection{Produit multilin\'eaire}
 
  On suppose $d=d'+1$. Fixons un sous-tore d\'eploy\'e maximal $A'$ de $G'$ et un sous-tore d\'eploy\'e maximal $A$ de $G$ qui contient $A'$. Fixons des sous-groupes paraboliques minimaux $P'_{min}$ de $G'$ et $P_{min}$ de $G$, contenant respectivement $A'$ et $A$. On peut les choisir de sorte que:
  
  $\bullet$ soient $n'\leq n$ les entiers tels que $A\simeq GL(1)^n$ et $A'\simeq GL(1)^{n'}$; alors le plongement de $A'$ dans $A$  s'identifie \`a $(a_{1},...,a_{n'})\mapsto (a_{1},...,a_{n'},1,...,1)$;
  
  $\bullet$ l'ensemble des racines simples de $a'$ relatif \`a $P'_{min}$ est form\'e  des  caract\`eres $(a_{1},...,a_{n'})\mapsto a_{i}a_{i+1}^{-1}$ pour $i=1,...,n'-1$ et de $(a_{1},...,a_{n'})\mapsto a_{n'}$, sauf dans le cas o\`u $d'=2n'$, auquel cas la derni\`ere racine est remplac\'ee par $(a_{1},...,a_{n'})\mapsto a_{n'-1}a_{n'}$; de m\^eme pour l'ensemble des racines simples de $T$ relatif \`a $P_{min}$.
  
   On peut supposer que les sous-groupes paraboliques servant \`a d\'efinir les repr\'esentations $\sigma$ et $\sigma'$ de 1.3 contiennent les sous-groupes paraboliques $P_{min}$, resp. $P'_{min}$.  On fixe des sous-groupes compacts sp\'eciaux $K'$ de $G'(F)$ et $K$ de $G(F)$, en bonne position relativement \`a $P'_{min}$ et $P_{min}$.
   
   En 1.3, on a suppos\'e que les $b_{i}$ et $b'_{i}$ \'etaient r\'eels et v\'erifiaient certaines in\'egalit\'es. Oublions ces conditions en prenant pour $b_{i}$ et $b'_{i}$ des nombres complexes quelconques. On introduit les param\`etres ${\bf z}=(z_{1},...,z_{t})$ et  ${\bf z}'=(z'_{1},...,z'_{t'})$, avec $z_{i}=q^{-b_{i}}$ et   $z'_{i}=q^{-b'_{i}}$, o\`u $q$ est le nombre d'\'el\'ements du corps r\'esiduel de $F$. On note plut\^ot nos repr\'esentations $\sigma_{{\bf z}}$ et $\sigma'_{{\bf z}'}$. Suivant Bernstein, on peut consid\'erer $\sigma_{{\bf z}}$ et $\sigma'_{{\bf z}'}$ comme les sp\'ecialisations pour ces valeurs des param\`etres de repr\'esentations \`a valeurs dans l'alg\`ebre 
   $${\cal R}={\mathbb C}[z_{1}^{\pm 1},...,z_{t}^{\pm 1},(z'_{1})^{\pm 1},...,(z'_{t'})^{\pm 1}].$$
   La repr\'esentation $\sigma_{{\bf z}}$ se r\'ealise dans un espace ${\cal E}$ de fonctions
   $$e:K\to E_{\pi_{1}}\otimes_{{\mathbb C}}...\otimes_{{\mathbb C}}E_{\pi_{t}}\otimes E_{\sigma_{0}}$$
   et cet espace est ind\'ependant de ${\bf z}$. On note $\check{\sigma}_{{\bf z}}$ la contragr\'ediente de $\sigma_{{\bf z}}$. Elle  se r\'ealise de m\^eme dans un espace $\check{\cal E}$ de fonctions
   $$\check{e}:K\to  E_{\check{\pi}_{1}}\otimes_{{\mathbb C}}...\otimes_{{\mathbb C}}E_{\check{\pi}_{t}}\otimes E_{\check{\sigma}_{0}}.$$
  On introduit le produit bilin\'eaire naturel sur ${\cal E}\times \check{\cal E}$:
  $$<e,\check{e}>=\int_{K}<e(k),\check{e}(k)>dk,$$
  o\`u le produit int\'erieur est l'accouplement naturel sur
$$(E_{\pi_{1}}\otimes_{{\mathbb C}}...\otimes_{{\mathbb C}}E_{\pi_{t}}\otimes E_{\sigma_{0}})\times (E_{\check{\pi}_{1}}\otimes_{{\mathbb C}}...\otimes_{{\mathbb C}}E_{\check{\pi}_{t}}\otimes E_{\check{\sigma}_{0}}).$$
Les m\^emes consid\'erations valent pour $\sigma'_{{\bf z}'}$ et on introduit des espaces ${\cal E}'$ et $\check{\cal E}'$, munis d'un produit bilin\'eaire. Pour $e\in {\cal E}$, $\check{e}\in \check{\cal E}$, $e'\in {\cal E}'$ et $\check{e}'\in \check{\cal E}'$, posons
$${\cal L}({\bf z},{\bf z}';e',\check{e}',e,\check{e})=\int_{G'(F)}<\sigma'_{{\bf z}'}(x)e', \check{e}'><\sigma_{{\bf z}}(x)e,\check{e}>dx.$$
Notons ${\cal D}$ le domaine de $({\mathbb C}^{\times})^n\times ({\mathbb C}^{\times})^{n'}$ d\'efini par les relations
$$q^{-1/2}<\, \vert z_{i}\vert ,\,\vert z_{j}\vert ,\,\vert z_{i}z'_{j}\vert ,\,\vert z_{i}z_{j}^{_{'}-1}\vert\, < q^{1/2},$$
les $i$, $j$ parcourant tous les entiers possibles.

\ass{Lemme}{(i) L'int\'egrale ${\cal L}({\bf z},{\bf z}';e',\check{e}',e,\check{e})$ est absolument convergente pour $({\bf z},{\bf z}')\in {\cal D}$.

(ii) Il existe un polyn\^ome non nul $D\in {\cal R}$ et,   pour tous $e$, $\check{e}$, $e'$, $\check{e}'$, il existe un polyn\^ome $L({\bf z},{\bf z}';e',\check{e}',e,\check{e})\in {\cal R}$ de sorte que
$$D({\bf z},{\bf z}'){\cal L}({\bf z},{\bf z}';e',\check{e}',e,\check{e})=L({\bf z},{\bf z}';e',\check{e}',e,\check{e})$$
pour tout $({\bf z},{\bf z}')\in {\cal D}$.}

Preuve. Pour simplifier, on suppose $d\not=2n$ et $d'\not=2n'$. On indiquera plus loin comment adapter la preuve si l'une de ces conditions n'est pas v\'erifi\'ee. Posons ${\cal A}=X_{*}(A')$ et ${\cal A}_{{\mathbb R}}=X_{*}(A')\otimes_{{\mathbb Z}}{\mathbb R}^{n'}$.  L'ensemble des racines simples de $A'$ relatif \`a $P'_{min}$ s'identifie \`a un ensemble $\Delta'$ de formes lin\'eaires sur ${\cal A}_{{\mathbb R}}$. On introduit l'ensemble des  poids $\{\varpi_{\alpha}; \alpha\in \Delta'\}$.
Fixons une uniformisante $\varpi_{F}$ de $F$ et identifions ${\cal A}$ \`a un sous-groupe de $A'(F)$ par
$${\bf m}=(m_{1},...,m_{n'})\mapsto (\varpi_{F}^{m_{1}},...,\varpi_{F}^{m_{t'}}).$$
Introduisons le sous-ensemble ${\cal A}^+$ form\'e des ${\bf m}\in {\cal A}$ tels que $m_{1}\geq m_{2}\geq...\geq m_{t'}\geq0$. D'apr\`es la d\'ecomposition de Cartan, il existe un sous-ensemble fini $\Gamma'$ de $G'(F)$, contenu dans le commutant de $A'$, de sorte que $G'(F)$ soit union disjointe des ensembles $K'{\bf m}\gamma'K'$ pour $ {\bf m}\in {\cal A}^+$ et $\gamma'\in \Gamma'$. On a
$${\cal L}({\bf z},{\bf z}';e',\check{e}',e,\check{e})= \int_{K'\times K'}\sum_{\gamma'\in \Gamma'}\sum_{{\bf m}\in {\cal A}^+}mes({\bf m}\gamma')<\sigma'_{{\bf z}'}(k_{1}{\bf m}\gamma'k_{2})e', \check{e}'>$$
$$<\sigma_{{\bf z}}(k_{1}{\bf m}\gamma'k_{2})e,\check{e}> dk_{1}dk_{2},$$
o\`u $mes({\bf m}\gamma')=mes(K'{\bf m}\gamma'K')mes(K')^{-2}$. Nos repr\'esentations \'etant lisses, cela nous permet de fixer $\gamma'\in \Gamma'$ et de remplacer l'int\'egrale  ${\cal L}({\bf z},{\bf z}';e',\check{e}',e,\check{e})$ par la s\'erie
$$S({\bf z},{\bf z}';e',\check{e}',e,\check{e})=\sum_{{\bf m}\in {\cal A}^+}mes({\bf m}\gamma')<\sigma'_{{\bf z}'}({\bf m})e', \check{e}'><\sigma_{{\bf z}}({\bf m})e,\check{e}>$$
  Consid\'erons un \'el\'ement $ {\bf T}=(T_{1},...,T_{t'})\in {\cal A}$ tel que $T_{1}>...>T_{t'}>0$. Un tel \'el\'ement permet de d\'ecomposer ${\cal A}^+$ en union disjointe de sous-ensembles ${\cal A}^+(Q')$, o\`u $Q'$ parcourt les sous-groupes paraboliques de $G'$ qui contiennent $P'_{min}$, cf. [A2] 3.9. A $Q'$ sont associ\'es un sous-ensemble $\Delta'(Q')\subset \Delta'$ et une d\'ecomposition  ${\cal A}_{{\mathbb R}}={\cal A}_{Q',{\mathbb R}}\oplus {\cal A}_{{\mathbb R}}^{Q'}$. Le premier sous-espace est l'intersection des annulateurs des $\alpha\in \Delta'(Q')$ et le second est engendr\'e par les coracines $\check{\alpha}$ pour $\alpha\in \Delta'(Q')$. On \'ecrit tout ${\bf m}\in {\cal A}_{{\mathbb R}}$ sous la forme ${\bf m}={\bf m}_{Q'}+{\bf m}^{Q'}$ conform\'ement \`a cette d\'ecomposition. L'ensemble ${\cal A}^+_{Q'}$ est celui des ${\bf m}\in {\cal A}^+$ tels que

$\bullet$ $\alpha({\bf m}_{Q'}-{\bf T}_{Q'})>0$ pour $\alpha\in \Delta\setminus \Delta'(Q')$;

$\bullet$ $\varpi_{\alpha}({\bf m}^{Q'}-{\bf T}^{Q'})\leq0$ pour $\alpha\in \Delta'(Q')$.

On peut fixer $Q'$ et remplacer $S({\bf z},{\bf z}';e',\check{e}',e,\check{e})$ par la s\'erie $S(Q';{\bf z},{\bf z}';e',\check{e}',e,\check{e})$ o\`u l'on restreint la somme aux ${\bf m}\in {\cal A}^+(Q')$. Notons $M'$ la composante de L\'evi de $Q'$ qui contient $A'$. On \'ecrit  $M'=GL(N_{1})\times...\times GL(N_{k})\times G'_{1}$, o\`u $G'_{1}$ est un groupe de m\^eme type que $G'$. Les \'el\'ements   $e'$ et $\check{e}'$ \'etant fix\'es, les r\'esultats de Casselman nous disent que, si $\alpha({\bf T})$ est assez grand pour tout $\alpha\in \Delta$, la propri\'et\'e suivante  est v\'erifi\'ee. Notons $p_{Q'}:E_{\sigma'_{{\bf z}'}}\to E_{\sigma'_{{\bf z}',Q'}}$ et $\check{p}_{\bar{Q}'}:E_{\check{\sigma}'_{{\bf z}'}}\to E_{(\check{\sigma}'_{{\bf z}'})_{\bar{Q}'}}$ les projections sur les modules de Jacquet. Alors il existe un produit bilin\'eaire $M'(F)$ invariant sur ces modules de sorte que, pour tout ${\bf m}\in {\cal A}^+(Q')$, on ait l'\'egalit\'e
$$<\sigma'_{{\bf z}'}({\bf m})e', \check{e}'> =\delta_{Q'}({\bf m})^{1/2}<\sigma'_{{\bf z}',Q'}({\bf m})p_{Q'}(e'), \check{p}_{\bar{Q}'}(\check{e}')>.$$
Remarquons qu'\`a $Q'$ est naturellement associ\'e un sous-groupe parabolique $Q$ de $G$ contenant $P_{min}$:   le L\'evi $M$ de $Q$ est $GL(N_{1})\times...\times GL(N_{k})\times G_{1}$, o\`u $G_{1}$ est de m\^eme type que $G$. Un \'el\'ement de ${\cal A}^+(Q')$ v\'erifie les m\^emes conditions relativement \`a $Q$ que celles indiqu\'ees ci-dessus et on peut appliquer le r\'esultat de Casselman. Les \'el\'ements $e$ et $\check{e}$ \'etant fix\'es,  si $\alpha({\bf T})$ est assez grand pour tout $\alpha\in \Delta$, on a une \'egalit\'e
$$<\sigma_{{\bf z}}({\bf m})e,\check{e}>=\delta_{Q}({\bf m})^{1/2}<\sigma_{{\bf z},Q}({\bf m})p_{Q}(e),p_{\bar{Q}}(\check{e})>$$
pour tout ${\bf m}\in {\cal A}^+(Q')$. Notons ${\cal C}(M')$ l'intersection de ${\cal A}$ avec le centre de $M'(F)$.  On a ${\cal C}(M')={\mathbb Z}^k$. Remarquons que ${\cal C}(M')$ est aussi contenu dans le centre de $M(F)$. La repr\'esentation $\sigma_{{\bf z},Q}$ est de longueur finie, la longueur \'etant born\'ee ind\'ependamment de ${\bf z}$. D'apr\`es le m\^eme calcul qu'en 1.4(2), les restrictions \`a ${\cal C}(M')={\mathbb Z}^k$ des caract\`eres centraux de ses sous-quotients irr\'eductibles sont de la forme
$${\bf c}=(c_{1},...,c_{k})\mapsto \chi_{{\bf z}}({\bf c})=\chi({\bf c}) \prod_{i=1,...,t;j=1,...,k}z_{i}^{(f_{i,j}-f_{-i,j})c_{j}},\eqno(1)$$
o\`u:

$\bullet$ les $f_{i,j}$ et $f_{-i,j}$ sont des entiers naturels v\'erifiant
$$\sum_{i=1,...,t}(f_{i,j}+f_{-i,j})=N_{j};$$

$\bullet$ $\chi$ est le caract\`ere central d'un sous-quotient irr\'eductible de $\sigma_{{\bf 1},Q}$, o\`u ${\bf 1}=(1,...,1)\in ({\mathbb C}^{\times})^n$.

Cela entra\^{\i}ne qu'il existe un entier $l\geq 0$ de sorte que, pour tous $e$, $\check{e}$ et tout ${\bf m}\in {\cal A}$, la fonction
$$ {\bf c}\mapsto <\sigma_{{\bf z},Q}({\bf c}{\bf m})p_{Q}(e),p_{\bar{Q}}(\check{e})>\eqno(2)$$
sur ${\cal C}(M')$ est combinaison lin\'eaire de fonctions
$$ {\bf c}\mapsto \chi_{{\bf z}}({\bf c})\prod_{j=1,...,k}c_{j}^{l_{j}},\eqno(3)$$
o\`u les $l_{j}$ sont des entiers naturels inf\'erieurs ou \'egaux \`a $l$. Un r\'esultat analogue vaut pour la fonction
$${\bf c}\mapsto <\sigma'_{{\bf z}',Q'}({\bf m})p_{Q'}(e'), \check{p}_{\bar{Q}'}(\check{e}')>  .$$

Le groupe ${\cal C}(M')$ agit par translations sur ${\cal A}$. On v\'erifie que l'ensemble ${\cal A}^+(Q')$ est contenu dans la r\'eunion d'un nombre fini d'orbites. On peut d\'ecomposer notre s\'erie  $S(Q';{\bf z},{\bf z}';e',\check{e}',e,\check{e})$ en somme finie de s\'eries o\`u l'on restreint l'ensemble de sommation \`a une intersection $({\bf m}{\cal C}(M'))\cap {\cal A}^+(Q')$. Consid\'erons une telle s\'erie.   L'application ${\bf c}\mapsto {\bf m}{\bf c}$ identifie l'intersection pr\'ec\'edente \`a un c\^one ${\cal C}(M')^+$ d\'efini par des in\'egalit\'es $c_{j}-c_{j+1}\geq C_{j}$ pour $j=1,...,k-1
$ et $c_{k}\geq C_{k}$, o\`u les $C_{j}$ sont certaines constantes. On v\'erifie qu'il y a une constante $m>0$ telle que  $mes(K'{\bf m}{\bf c}\gamma'K')=m\delta_{Q'}({\bf c})^{-1}$ pour tout  ${\bf c}$ dans ce c\^one. Alors notre s\'erie est born\'ee par une somme de s\'eries de la forme
$$ \sum_{{\bf c}\in {\cal C}(M')^+}\delta_{Q'}({\bf c})^{-1/2}\delta_{Q}({\bf c})^{1/2}\chi'_{{\bf z}'}({\bf c})\chi_{{\bf z}}({\bf c})\prod_{j=1,...,k}\vert c_{j}\vert ^{l_{j}+l'_{j}}.\eqno(4)$$
On calcule
$$\delta_{Q'}({\bf c})^{-1/2}\delta_{Q}({\bf c})^{1/2}=q^{-\sum_{j=1,...,k}N_{j}c_{j}/2}.$$
Consid\'erons l'expression (1). La repr\'esentation $\sigma_{{\bf 1}}$ est temp\'er\'ee, donc $\chi$ est born\'ee sur ${\cal C}(M')^+$. Si on note $\underline{z}$ le plus grand des nombres $\vert z_{i}\vert ^{\pm 1}$,  $\chi_{{\bf z}}({\bf c})$ est essentiellement born\'e sur ${\cal C}(M')^+$ par $\prod_{j=1,...,k}\underline{z}^{\sum_{j=1,...,k}N_{j}c_{j}}$. De m\^eme,  $\chi'_{{\bf z}'}({\bf c})$ est essentiellement born\'e   par $\prod_{j=1,...,k}(\underline{z}')^{\sum_{j=1,...,k}N_{j}c_{j}}$. La s\'erie (4) est donc essentiellement born\'ee par
$$\sum_{{\bf c}\in {\cal C}(M')^+}\prod_{j=1,...,t}(\underline{z}\underline{z}'q^{-1/2})^{N_{j}c_{j}}\vert c_{j}\vert ^{l_{j}+l'_{j}}.$$
Cette derni\`ere s\'erie est convergente sous les hypoth\`eses du (i) de l'\'enonc\'e et cela d\'emontre cette assertion.

On v\'erifie que la somme de la s\'erie (4) est une fraction rationnelle en ${\bf z}$, ${\bf z}'$. Remarquons que les termes qui y interviennent parcourent des ensembles finis ind\'ependants des vecteurs $e$, $\check{e}$, $e'$, $\check{e}'$. En reprenant le calcul ci-dessus, on voit que, pour d\'emontrer le (ii) de l'\'enonc\'e, il suffit de prouver l'assertion suivante. Consid\'erons les diff\'erentes fonctions de la forme (3) qui peuvent intervenir. Elles sont d\'etermin\'ees par un caract\`ere $\chi$ et des familles d'entiers $f_{i,j}-f_{-i,j}$ et $l_{j}$, ces donn\'ees parcourant des ensembles finis. Notons $(f_{h,{\bf z}})_{h\in H}$ cette famille de fonctions, l'ensemble d'indices $H$ \'etant donc fini.
 Notons $  f_{{\bf z}}$ la fonction (2). Ecrivons 
$$f_{{\bf z}}=\sum_{h\in H}C_{h}({\bf z})f_{h,{\bf z}}.$$
On doit prouver que les diff\'erents coefficients $C_{h}({\bf z})$ sont des fractions rationnelles en ${\bf z}$, de d\'enominateur born\'e ind\'ependamment de $e$ et $\check{e}$. On doit aussi prouver l'assertion similaire relative au groupe $G'$, mais elle  se prouve \'evidemment de la m\^eme fa\c{c}on.  Pour ${\bf z}$ en position g\'en\'erale, la famille de fonctions  $(f_{h,{\bf z}})_{h\in H}$ est lin\'eairement ind\'ependante. On peut donc fixer une famille $({\bf c}_{h})_{h\in H}$ d'\'el\'ements de ${\cal C}(M')$ telle que le d\'eterminant de la matrice $(f_{h,{\bf z}}({\bf c}_{h'}))_{h,h'\in H}$ soit non nul pour au moins une valeur de ${\bf z}$. Ce d\'eterminant est donc un \'el\'ement non nul de ${\cal R}$. Les coefficients $C_{h}({\bf z})$ sont d\'etermin\'es par le syst\`eme d'\'equations
$$f_{{\bf z}}({\bf c}_{h'})=\sum_{h\in H}C_{h}({\bf z})f_{h,{\bf z}}({\bf c}_{h'})$$
pour tout $h'\in H$. Le membre de gauche de cette \'equation appartient \`a ${\cal R}$. Il en r\'esulte que les $C_{h}({\bf z})$ sont des fractions rationnelles, de d\'enominateur divisant le d\'eterminant ci-dessus, lequel ne d\'epend pas des vecteurs $e$ et $\check{e}$. Cela d\'emontre l'assertion requise et  le lemme, sous les hypoth\`eses $d\not=2n$, $d'\not=2n'$.

Supposons $d=2n$. Cela implique que $G$ est d\'eploy\'e et l'hypoth\`ese sur le plongement de $V'$ dans $V$ implique que $G'$ est lui-aussi d\'eploy\'e.  Donc $n'=n-1$. La seule chose qui change dans le raisonnement ci-dessus est que, si $N_{1}+...+N_{k}=n'$, le L\'evi du parabolique $Q$ associ\'e \`a $Q'$ est $GL(N_{1})\times...\times GL(N_{k})\times GL(1)$. Cela n'a aucune incidence. Supposons maintenant $d'=2n'$. Les deux groupes $G$ et $G'$ sont d\'eploy\'es et 
on a $n=n'$. Dans la d\'efinition de ${\cal A}^+$, la derni\`ere in\'egalit\'e $m_{n'}\geq0$ doit a priori \^etre remplac\'ee par $m_{n'-1}+m_{n'}\geq0$. Mais on peut d\'ecomposer cet ensemble en deux, l'un sur lequel $m_{n'}\geq0$ et l'autre sur lequel $m_{n'}<0$. En changeant l'identification de $A'$ avec $GL(1)^{n'}$ en inversant la derni\`ere coordonn\'ee, et en changeant en cons\'equence le groupe $P_{min}$, le deuxi\`eme ensemble se ram\`ene \`a un ensemble du premier type (la condition sur $m_{n'}$ devient $m_{n'}>0$ au lieu de $m_{n'}\geq0$, mais c'est sans importance). On peut donc consid\'erer que ${\cal A}^+$ est d\'efini par les m\^emes in\'egalit\'es que pr\'ec\'edemment. Pour un sous-groupe parabolique $Q'$ tel que $N_{1}+...+N_{k}<n'$, rien ne change. 

Soit $Q'$  tel que $N_{1}+...+N_{k}=n'$ et $N_{k}\geq2$. Si $\Delta'(Q')$ contient la racine ${\bf m}\mapsto {\bf m}_{n'-1}-{\bf m}_{n'}$, de nouveau rien ne change. Supposons que $\Delta'(Q')$ contienne la racine ${\bf m}\mapsto {\bf m}_{n'-1}+{\bf m}_{n'}$. On d\'efinit $Q$ de sorte que sa composante de L\'evi $M$ soit \'egale \`a $GL(N_{1})\times...\times GL(N_{k-1})\times G_{1}$.  On v\'erifie les deux propri\'et\'es suivantes:

(5) notons $\Delta$ l'ensemble des racines simples de $A$ associ\'e \`a $P_{min}$ et $\Delta(Q)$ le sous-ensemble associ\'e \`a $Q$; pour ${\bf m}\in {\cal A}^+(Q')$ et $\alpha\in \Delta\setminus \Delta(Q)$, on a $\alpha({\bf m})>\alpha({\bf T})$;

(6) il existe $C\in {\mathbb N}$ tel que $m_{i}\leq C$ pour tout ${\bf m}\in {\cal A}^+(Q')$ et tout $i=N_{1}+...+N_{k-1}+1,...,n'$.

L'assertion (5) vient de l'inclusion $\Delta\setminus \Delta(Q)\subset \Delta'\setminus \Delta'(Q')$. D\'emontrons (6). Posons $e=N_{1}+...+N_{k-1}$. Soit ${\bf m}\in {\cal A}^+(Q')$.  Cette hypoth\`ese implique que l'on peut \'ecrire
$$(m_{e+1}-t_{e+1},...,m_{n'}-t_{n'})=(z,...,z,-z)+(p_{1},..., p_{N_{k}-1},-p_{N_{k}}),$$
avec $z>0$, $p_{1}+...+p_{f}\leq0$ pour tout $f=1,...,N_{k-1}$, $p_{1}+...+p_{N_{k}}=0$. Cela entra\^{\i}ne $p_{N_{k}}\geq0$. Alors  $z=t_{n'}-m_{n'}-p_{N_{k}}$ est major\'e (puisque l'on a suppos\'e $m_{n'}\geq0$). On a aussi $p_{1}\leq0$, donc $m_{e+1}=t_{e+1}+z+p_{1}$ est major\'e. Puisque $m_{e+1}\geq m_{e+2}\geq...\geq m_{n'}\geq0$, (6) s'ensuit. 

La relation (5) nous permet d'appliquer les r\'esultats de Casselman aux termes provenant du groupe $G$, pour le sous-groupe parabolique $Q$. On remplace dans les raisonnements ci-dessus le groupe ${\cal C}(M')$ par son analogue ${\cal C}(M)$ pour le groupe $M$. Gr\^ace \`a la relation (6), ${\cal A}^+(Q')$ est inclus  dans un nombre fini d'orbites pour l'action de ce groupe. La preuve se poursuit alors comme pr\'ec\'edemment. 

Soit enfin $Q'$ tel que $N_{1}+...+N_{k}=n'$ et $N_{k}=1$. On d\'ecompose notre s\'erie en une somme sur les ${\bf m}$ tels que $m_{n'}> t_{n'}$ et d'une somme sur les $m_{n'}\leq t_{n'}$. La premi\`ere se traite comme dans le cas g\'en\'eral, en prenant $Q$ de composante de L\'evi $GL(N_{1})\times...\times GL(N_{k})$. La seconde se traite comme dans le cas particulier ci-dessus en prenant $Q$ de composante de L\'evi $GL(N_{1})\times...\times GL(N_{k-1})\times G_{1}$. Cela ach\`eve la preuve. $\square$

\bigskip

\subsection{Preuve de l'in\'egalit\'e $m(\sigma_{0},\sigma'_{0})\leq m(\sigma,\check{\sigma}')$}

On consid\`ere la situation de 1.3 et on suppose d'abord $d=d'+1$. Il n'y a rien \`a prouver si $m(\sigma_{0},\sigma'_{0})=0$. On suppose donc $m(\sigma_{0},\sigma'_{0})=1$.  Introduisons un polyn\^ome $D({\bf z},{\bf z}')$ v\'erifiant les conditions du (ii) du lemme 1.7. Fixons  d'abord des familles  ${\bf z}$ et ${\bf z}'$, telles que $D({\bf z},{\bf z}')\not=0$  et dont toutes les coordonn\'ees sont de module $1$. Les repr\'esentations $\sigma_{{\bf z}}$ et $\sigma'_{{\bf z}'}$ sont temp\'er\'ees. La proposition 5.7 et le lemme 5.3 de [W3] nous disent qu'il existe $e\in {\cal E}$, $\check{e}\in \check{{\cal E}}$, $e'\in {\cal E}'$ et $\check{e}'\in \check{{\cal E}}'$ tels que ${\cal L}({\bf z},{\bf z}';e',\check{e}',e,\check{e})\not=0$. 

{\bf Remarque.} Dans [W3], on a consid\'er\'e des produits hermitiens plut\^ot que des produits bilin\'eaires, mais la traduction est facile puisque les repr\'esentations sont unitaires.

On a donc aussi $L({\bf z},{\bf z}';e',\check{e}',e,\check{e})\not=0$. Ecrivons les coordonn\'ees 
$z_{i}$ et $z'_{i}$ de nos familles sous la forme $z_{i}=q^{-\beta_{i}}$, $z'_{i}=q^{-\beta'_{i}}$. Pour $s\in {\mathbb C}$, introduisons les familles ${\bf z}(s)$ et ${\bf z}'(s)$ de coordonn\'ees $z_{i}(s)=q^{-s\beta_{i}+(s-1)b_{i}}$, $z'_{i}(s)=q^{-s\beta'_{i}+(s-1)b'_{i}}$. Pour tous $e$, $\check{e}$, $e'$, $\check{e}'$, l'application $s\mapsto L({\bf z}(s),{\bf z}'(s);e',\check{e}',e,\check{e})$ est holomorphe. Elle est non nulle pour au moins un choix de $e$, $\check{e}$, $e'$, $\check{e}'$. Notons $N$ le plus grand entier tel que, pour tous $e$, $\check{e}$, $e'$, $\check{e}'$, cette application soit divisible par $s^N$. Posons
$$L_{\sigma,\sigma'}(e',\check{e}',e,\check{e})=lim_{s\to 0}s^{-N} L({\bf z}(s),{\bf z}'(s);e',\check{e}',e,\check{e}).$$
Il y a au moins un choix de vecteurs $e$, $\check{e}$, $e'$, $\check{e}'$ tel que $L_{\sigma,\sigma'}(e',\check{e}',e,\check{e})\not=0$. Dans le domaine ${\cal D}$, on v\'erifie ais\'ement l'\'egalit\'e
$${\cal L}({\bf z},{\bf z}';\sigma'_{{\bf z}'}(g')e', \check{e}',\sigma_{{\bf z}}(g')e,\check{e})={\cal L}({\bf z},{\bf z}';e',\check{e}',e,\check{e})$$
pour tout $g'\in G'(F)$. On a une m\^eme relation pour l'application $L({\bf z},{\bf z}';e',\check{e}',e,\check{e})$ et cette \'egalit\'e se prolonge alg\'ebriquement \`a toutes familles ${\bf z}$, ${\bf z}'$. Il en r\'esulte que l'application $L_{\sigma,\sigma'}$ ci-dessus v\'erifie la relation
$$L_{\sigma,\sigma'}(\sigma'(g')e', \check{e}',\sigma(g')e,\check{e})=L_{\sigma,\sigma'}(e',\check{e}',e,\check{e}).$$
 Fixons alors $\check{e}$ et $\check{e}'$ et d\'efinissons une application lin\'eaire $l:E_{\sigma}\to E_{\check{\sigma}'}$ par l'\'egalit\'e
$$<e',l(e)>=L_{\sigma,\sigma'}(e',\check{e}',e,\check{e}).$$
Elle appartient \`a $Hom_{G'}(\sigma,\check{\sigma}')$. Pour un bon choix de $\check{e}$, $\check{e}'$, elle est non nulle. Donc $m(\sigma,\check{\sigma}')\geq1$.
  
 Supposons maintenant $d>d'+1$. On choisit $\rho$ comme en 1.6(2). On a $m(\sigma,\check{\sigma}')=m(\rho\times\sigma',\check{\sigma})$. D'apr\`es ce que l'on vient de prouver, cette derni\`ere multiplicit\'e est sup\'erieure ou \'egale \`a $1$.
 
 L'in\'egalit\'e que l'on vient de prouver et celle de 1.6 d\'emontrent la proposition 1.3. 
 
 \bigskip

\section{Irr\'eductibilit\'e et repr\'esentations g\'en\'eriques}

 \bigskip
 
 \subsection{Rappels sur les param\'etrages}
 
 Dans cette section,  $G$ est un groupe sp\'ecial orthogonal ou symplectique d\'efini sur $F$. Pr\'ecis\'ement, cela signifie que l'on fixe un espace vectoriel $V$ sur $F$ de dimension finie et une forme bilin\'eaire $q$ non d\'eg\'en\'er\'ee sur $V$ qui est soit sym\'etrique, soit antisym\'etrique. Alors $G$ est la composante neutre du groupe d'automorphismes de $(V,q)$. On note $d_{G}$ la dimension de $V$. Un groupe  similaire $G'$ est dit de m\^eme type que $G$ s'il est associ\'e \`a un couple $(V',q')$ satisfaisant les conditions suivantes: $q'$ v\'erifie la m\^eme condition de sym\'etrie que $q$; le plus grand des espaces $(V,q)$ et $(V',q')$ est isomorphe \`a la somme orthogonale du plus petit et de plans hyperboliques.
 
 On consid\`ere que le $L$-groupe $\hat{G}$ de $G$ est:
 
 $\bullet$ le groupe sp\'ecial orthogonal complexe $SO(\hat{d}_{G},{\mathbb C})$, o\`u $\hat{d}_{G}=d_{G}+1$, si $G$ est symplectique;
 
 $\bullet$ le groupe symplectique complexe $Sp(\hat{d}_{G},{\mathbb C})$, o\`u $\hat{d}_{G}=d_{G}-1$, si $G$ est sp\'ecial orthogonal "impair" (c'est-\`a-dire que $d_{G}$ est impair);
 
 $\bullet$ le groupe orthogonal complexe $O(\hat{d}_{G},{\mathbb C})$, o\`u $\hat{d}_{G}=d_{G}$, si $G$ est sp\'ecial orthogonal "pair" (c'est-\`a-dire que $d_{G}$ est pair).
 
 Toutes les repr\'esentations de groupes r\'eductifs que l'on consid\'erera seront suppos\'ees admissibles et de longueur finie. Dans les cas symplectique ou sp\'ecial orthogonal impair, on conjecture (et nous admettons cette conjecture) que l'ensemble des classes de conjugaison de repr\'esentations irr\'eductibles temp\'er\'ees de $G(F)$ se d\'ecompose en union disjointe de $L$-paquets $\Pi^G(\varphi)$, o\`u $\varphi$ parcourt les classes de conjugaison par $\hat{G}$ d'homomorphismes $\varphi:W_{DF}\to \hat{G}$ qui v\'erifient quelques conditions usuelles de continuit\'e et de semi-simplicit\'e et qui sont temp\'er\'es, c'est-\`a-dire que l'image de $W_{F}$ par $\varphi$ est relativement compacte. Un tel homomorphisme $\varphi$ se pousse en un homomorphisme de $W_{DF}$ dans $GL(\hat{d}_{G},{\mathbb C})$. Via la correspondance de Langlands (th\'eor\`eme de Harris-Taylor et Henniart), on associe \`a $\varphi$ une repr\'esentation irr\'eductible $\pi(\varphi)$ de $GL(\hat{d}_{G},F)$. Elle est temp\'er\'ee et autoduale. On peut la prolonger en une repr\'esentation du produit semi-direct $GL(\hat{d}_{G},F)\rtimes \{1,\theta\}$, o\`u $\theta$ est l'automorphisme ext\'erieur habituel d\'efini par $\theta(g)={^tg}^{-1}$. Notons $\tilde{\pi}(\varphi)$ sa restriction \`a la composante connexe non neutre $GL(\hat{d}_{G},F)\theta$. D'autre part, pour toute repr\'esentation $\sigma$, notons $\Theta_{\sigma}$ son caract\`ere. Les propri\'et\'es essentielles du $L$-paquet $\Pi^G(\varphi)$ sont les suivantes:
 
 (1) la somme $\Theta_{\Pi^G(\varphi)}=\sum_{\pi\in \Pi^G(\varphi)}\Theta_{\pi}$ est une distribution stable;
 
 (2) si $G$ est quasi-d\'eploy\'e, il existe $c\in {\mathbb C}^{\times}$ tel que la distribution $\Theta_{\tilde{\pi}(\varphi)}$  soit le transfert, par endoscopie tordue, de $c\Theta_{\Pi^G(\varphi)}$;
 
 (3) si $G$ n'est pas quasi-d\'eploy\'e, introduisons sa forme quasi-d\'eploy\'ee $\underline{G}$ et le $L$-paquet $\Pi^{\underline{G}}(\varphi)$ de repr\'esentations de $\underline{G}(F)$; alors  il existe $c\in {\mathbb C}^{\times}$ tel que $\Theta_{\Pi^G(\varphi)}$ soit le transfert, par endoscopie ordinaire, de $c\Theta_{\Pi^{\underline{G}}(\varphi)}$.
 
 Ici, les correspondances endoscopiques entre classes de conjugaison stable sont injectives. Les conditions (2) et (3) se traduisent concr\`etement de la fa\c{c}on suivante. Soit $g\in G(F)$ un \'el\'ement semi-simple fortement r\'egulier. Alors on a une \'egalit\'e
 $$\Theta_{\Pi^G(\varphi)}(g)=c\sum_{\tilde{x}}\Delta(g,\tilde{x})\Theta_{\tilde{\pi}(\varphi)}(\tilde{x}),\eqno(4)$$
 o\`u $\tilde{x}$ parcourt un certain sous-ensemble fini de $GL(\hat{d}_{G},F)\theta$ associ\'e \`a $g$ et $\Delta(g,\tilde{x})$ est un facteur de transfert (c'est l'inverse du facteur de Kottwitz-Shelstad). Par ind\'ependance lin\'eaire des caract\`eres, ces \'egalit\'es d\'eterminent uniquement le paquet $\Pi^G(\varphi)$. 
 
 {\bf Remarque.} Dans le cas d'un groupe orthogonal impair, les conjectures pos\'ees en [W1] 4.2 faisaient intervenir l'endoscopie tordue entre $G$ et $GL(\hat{d}_{G}+1)\theta$. Ici, on utilise l'endoscopie tordue plus habituelle entre $G$ et $GL(\hat{d}_{G})\theta$. Pour la validit\'e des r\'esultats de cette section, celle-ci suffit. Mais, pour le reste de l'article, on doit admettre aussi la validit\'e des conjectures telles qu'on les a formul\'ees en [W1].

  Dans le cas d'un groupe sp\'ecial orthogonal pair, il faut imposer \`a $\varphi$ une condition portant sur le d\'eterminant $det\circ\varphi$ (celui-ci doit correspondre au discriminant de la forme $q$, cf. [M1] paragraphe 2.1 et l'introduction ci-dessus), et on consid\`ere les classes de conjugaison de tels $\varphi$ par $SO(\hat{d}_{G},{\mathbb C})$ et non par $O(\hat{d}_{G},{\mathbb C})$. Le paquet $\Pi^G(\varphi)$ v\'erifie les m\^emes conditions que ci-dessus. Mais les correspondances entre classes de conjugaison stable ne sont plus injectives: deux classes conjugu\'ees par le groupe orthogonal tout entier ne sont pas discernables par endoscopie tordue. Notons $G^+$ ce groupe orthogonal et fixons un \'el\'ement $w\in G^+(F)\setminus G(F)$. La relation (4) devient
 $$\Theta_{\Pi^G(\varphi)}(g)+\Theta_{\Pi^G(\varphi)}(wgw^{-1})=c\sum_{\tilde{x}}\Delta(g,\tilde{x})\Theta_{\tilde{\pi}(\varphi)}(\tilde{x}).$$ 
 Ces relations ne d\'eterminent plus $\Pi^G(\varphi)$. Toutefois, pour toute repr\'esentation $\sigma$ de $G(F)$, notons $\sigma^w$ la repr\'esentation $g\mapsto \sigma(wgw^{-1})$ et posons $\Pi^G(\varphi)^w=\{\pi^w; \pi\in \Pi^G(\varphi)\}$. Alors l'ensemble avec multiplicit\'es $\bar{\Pi}(\varphi)=\Pi^G(\varphi)\sqcup \Pi^G(\varphi)^w$ est uniquement d\'etermin\'e.  
 
 Comme on l'a expliqu\'e dans l'introduction, la d\'efinition des $L$-paquets s'\'etend au cas non temp\'er\'e en utilisant la classification de Langlands.    
  
 {\bf Notation.} Soit $\pi$ une repr\'esentation irr\'eductible de $G(F)$. On pose $\pi^{GL}=\pi(\varphi)$, o\`u $\varphi$ est le param\`etre de Langlands tel que $\pi\in \Pi^G(\varphi)$. C'est une repr\'esentation de $GL(\hat{d}_{G},F)$. 
 
 \bigskip
 
 \subsection{Induction et modules de Jacquet}  
 
 Soit $\pi$ une repr\'esentation irr\'eductible temp\'er\'ee d'un groupe lin\'eaire $GL(d_{\pi},F)$. Le support cuspidal de $\pi$ a la forme suivante: il existe un ensemble fini avec multiplicit\'es, que l'on note $Jord(\pi)$, de couples $(\rho,a)$, o\`u $\rho$ est une repr\'esentation  irr\'eductible cuspidale et  unitaire d'un groupe lin\'eaire et $a \geq1$ est un entier,  de telle sorte que le support cuspidal de $\pi$ est exactement 
$$\cup_{(\rho,a)\in Jord(\pi)}\cup_{x\in [(a-1)/2,-(a-1)/2]}\{\rho\vert.\vert_{F}^{x}\}.
$$

On a donc d\'efini $Jord(\pi)$ pour toute repr\'esentation irr\'eductible temp\'er\'ee $\pi$ d'un groupe lin\'eaire. On transpose cela en une d\'efinition de $Jord(\pi)$ pour toute repr\'esentation irr\'eductible temp\'er\'ee $\pi$ de $G(F)$, on posant $Jord(\pi)=Jord(\pi^{GL})$, avec la notation introduite en 2.1.  Remarquons que $Jord(\pi)$ ne d\'epend que du $L$-paquet $\Pi$ contenant $\pi$, ce qui permet de d\'efinir $Jord(\Pi)$ pour un tel $L$-paquet.

  Soit $\pi$ une repr\'esentation irr\'eductible temp\'er\'ee de $G(F)$. L'ensemble $Jord(\pi)$ a les propri\'et\'es suivantes. Pour tout $(\rho,a)\in Jord(\pi)$, ou bien la mutliplicit\'e de $(\rho,a)$ dans $Jord(\pi)$ est paire, ou bien la repr\'esentation de $W_{DF}$ associ\'ee \`a la repr\'esentation de Steinberg $St(\rho,a)$  est \`a valeurs dans un groupe classique de m\^eme type que le groupe dual de $G$. Dans ce dernier cas, on dira que $(\rho,a)$ a bonne parit\'e. Dans  les autres cas, on dira que $(\rho,a)$ n'a pas bonne parit\'e. Ce dernier cas  couvre \`a la fois celui o\`u $\rho$ n'est pas autodual et celui o\`u le param\`etre de Langlands de $St(\rho,a)$, bien qu'autodual,  ne se factorise pas par le bon groupe classique.

Le facteur de transfert de Kottwitz et Shelstad est compatible \`a l'induction. Cela signifie la chose suivante. Supposons $G$ quasi-d\'eploy\'e et fixons un tel facteur $\Delta$ pour l'endoscopie tordue entre $G$ et $GL(\hat{d}_{G})\theta$. Soit $L$ un L\'evi de $G$. Il lui correspond un L\'evi $\theta$-stable ${\bf L}$ de $GL(\hat{d}_{G})$.  Soient $g\in L(F)$ et $\tilde{x}\in {\bf L}(F)\theta$ des \'el\'ements suffisamment r\'eguliers. D\'efinissons un facteur $\Delta_{L,{\bf L}\theta}(g,\tilde{x})$ par

$\Delta_{L,{\bf L}\theta}(g,\tilde{x})= \Delta(g,\tilde{x})$ si les classes de conjugaison stable  de $g$ dans $L(F)$ et de $\tilde{x}$ dans ${\bf L}(F)\theta$ se correspondent;

$\Delta_{L,{\bf L}\theta}(g,\tilde{x})=0$ sinon.

Alors $\Delta_{L,{\bf L}\theta}$ est un facteur de transfert pour le couple $(L,{\bf L}\theta)$. Cela entra\^{\i}ne que le transfert est compatible \`a l'induction et cela a la cons\'equence suivante. Soit  $\Pi$ un paquet de repr\'esentations temp\'er\'ees de $G(F)$. Supposons donn\'ee  une d\'ecomposition en union disjointe (au sens des ensembles avec multiplicit\'es):
$$Jord(\Pi)={\cal E}\sqcup {\cal F}\sqcup \check{{\cal E}},$$
o\`u $\check{{\cal E}}=\{(\check{\rho},a); (\rho,a)\in {\cal E}\}$. Il existe un $L$-paquet temp\'er\'e $\Pi'$ d'un groupe de m\^eme type que $G$ tel que ${\cal F}=Jord(\Pi')$. Alors

(1) si $G$ est symplectique ou sp\'ecial orthogonal impair, $\Pi$ est exactement form\'e des  composantes irr\'eductibles de toutes les induites 
$$(\times_{(\rho,a)\in {\cal E}}St(\rho,a))\times \pi'$$
quand $\pi'$ d\'ecrit $\Pi'$.

Dans le cas o\`u $G$ est sp\'ecial orthogonal pair, on a une assertion analogue. On doit remplacer $\Pi$ et $\Pi'$ par les paquets $\bar{\Pi}$ et $\bar{\Pi}'$ d\'efinis en 2.1. Dans le cas o\`u le sous-groupe parabolique $P$ servant \`a d\'efinir l'induite ci-dessus n'est pas semblable \`a $wPw^{-1}$, il faut sommer sur les deux induites possibles.

Notons $\varphi$ et $\varphi'$ les param\`etres de Langlands de $\Pi$ et $\Pi'$. L'assertion (1)  r\'esulte simplement du fait que la somme des caract\`eres des composantes en question est stable et se transf\`ere en le caract\`ere de l'induite tordue 
$$(\times_{(\rho,a)\in {\cal E}}St(\rho,a))\times \tilde{\pi}(\varphi').$$
Or cette derni\`ere n'est autre que $\tilde{\pi}(\varphi)$. $\square$

En particulier, un paquet temp\'er\'e $\Pi$ est form\'e de repr\'esentations de la s\'erie discr\`ete si et seulement si tous les \'el\'ements de $Jord(\Pi)$  sont de bonne parit\'e et interviennent avec multiplicit\'e $1$.

Une autre cons\'equence concerne les modules de Jacquet.  Fixons un sous-groupe parabolique $P$ de $G$ de L\'evi $GL(d_{1})\times...\times GL(d_{m})\times G'$.  Pour $i=1,...,m$, soient $\rho_{i}$ une repr\'esentation irr\'eductible cuspidale  de $GL(d_{i},F)$, pas forc\'ement unitaire. Pour  une repr\'esentation irr\'eductible $\pi$ de $G(F)$, notons $Jac_{\rho_{1},..., \rho_{m}}(\pi)$ la repr\'esentation semi-simple de $G'(F)$ telle que le semi-simplifi\'e du module de Jacquet $\pi_{P}$ soit la somme de $\rho_{1}\otimes...\otimes \rho_{m}\otimes Jac_{\rho_{1},..., \rho_{m}}(\pi)$ et de repr\'esentations dont les premi\`eres composantes ne sont pas \'egales \`a $\rho_{1}\otimes...\otimes \rho_{m}$.  Soit $\pi$ une repr\'esentation irr\'eductible temp\'er\'ee de $G(F)$. On a

(2) supposons $Jac_{\rho_{1},...,\rho_{m}}(\pi)\not=\{0\}$; alors pour tout $(\rho,a)\in Jord(\pi)$, il existe deux \'el\'ements $b_{\rho,a} $ et $\check{b}_{\rho,a}$ de $ [(a+1)/2,-(a+1)/2]$, avec $b_{\rho,a}>\check{b}_{\rho,a}$ tels que:

$\bullet$ la famille $(\rho_{i})_{i=1,...,m}$ s'obtienne en m\'elangeant les familles  $(\rho\vert .\vert _{F}^x)_{x\in [(a-1)/2,b_{{\rho,a}}]}$, pour $(\rho,a)\in Jord(\pi)$,  sans changer l'ordre dans chacune d'elles;

$\bullet$ la famille $(\check{\rho}_{i})_{i=m,...,1}$ s'obtienne en m\'elangeant les familles  $(\rho\vert .\vert _{F}^x)_{x\in [\check{b}_{\rho,a},-(a-1)/2]}$, pour $(\rho,a)\in Jord(\pi)$,  sans changer l'ordre dans chacune d'elles.

En effet, supposons $G$ symplectique ou sp\'ecial orthogonal impair. Soit $\Pi$ le $L$-paquet contenant $\pi$. D'apr\`es un lemme analogue \`a [MW] lemme 4.2.1, la somme des caract\`eres des $Jac_{\rho_{1},...,\rho_{m}}(\pi_{1})$, pour $\pi_{1}\in\Pi$, est stable et a pour transfert un certain module de Jacquet "tordu" $Jac^{\theta}_{\rho_{1},...,\rho_{m}}(\pi^{GL})$. En tant que repr\'esentation d'un groupe lin\'eaire non tordu, celui-ci se construit de fa\c{c}on analogue \`a ci-dessus. On note ${\bf P}$ un sous-groupe parabolique de $GL(\hat{d}_{G}) $ de L\'evi 
$$GL(d_{1})\times...\times GL(d_{m})\times GL(d_{0})\times GL(d_{m})\times...\times GL(d_{1}).$$
Alors $Jac^{\theta}_{\rho_{1},...,\rho_{m}}(\pi^{GL})$ est la repr\'esentation semi-simple de $GL(d_{0},F)$ telle que le module de Jacquet $(\pi^{GL})_{{\bf P}}$ soit la somme de
$$\rho_{1}\otimes...\otimes \rho_{m}\otimes Jac^{\theta}_{\rho_{1},...,\rho_{m}}(\pi^{GL}))\otimes \check{\rho}_{m}\otimes...\otimes \check{\rho}_{1},$$
et de repr\'esentations dont les composantes non centrales ne sont pas celles ci-dessus.  L'hypoth\`ese que $Jac_{\rho_{1},...,\rho_{m}}(\pi)\not=\{0\}$ impose que ce module de Jacquet tordu n'est pas nul. Mais on sait bien calculer $(\pi^{GL})_{{\bf P}}$ et la conclusion s'ensuit. Dans le cas o\`u $G$ est sp\'ecial orthogonal pair, la d\'emonstration s'adapte en consid\'erant le paquet $\bar{\Pi}$. $\square$

On a un r\'esultat plus pr\'ecis pour des familles particuli\`eres. Par exemple, fixons une repr\'esentation irr\'eductible $\rho$ cuspidale et unitaire d'un groupe lin\'eaire  et un r\'eel $x>0$ et supposons $\rho_{i}=\rho\vert .\vert _{F}^x$ pour tout $i=1,...,m$. On a

(3) supposons $Jac_{\rho_{1},...,\rho_{m}}(\pi)\not=\{0\}$; alors $x$ est un demi-entier et $(\rho,2x+1)$ intervient avec multiplicit\'e au moins $m$ dans $Jord(\pi)$.

Cela r\'esulte imm\'ediatement de (2). $\square$

Soit $ \Pi$ un paquet de repr\'esentations temp\'er\'ees de $G(F)$. Soit $(\rho,a)\in Jord(\Pi)$ avec $a\geq2$, notons $m$ sa multiplicit\'e.  Notons $\Pi^-$ le paquet temp\'er\'e tel que $Jord(\Pi^-)$ se d\'eduise de $Jord(\Pi)$ en rempla\c{c}ant les $m$ copies de $(\rho,a)$ par $m$ copies de $(\rho,a-2)$. Prenons $\rho_{i}=\rho\vert .\vert _{F}^{(a-1)/2}$ pour $i=1,...,m$. Soit $\pi\in \Pi$. Alors

(4)   toutes les composantes irr\'eductibles de $Jac_{\rho_{1},...,\rho_{m}}(\pi)$ appartiennent \`a $\Pi^-$ dans le cas symplectique ou orthogonal impair, \`a $\bar{\Pi}^-$ dans le cas orthogonal pair.

En notant $\varphi$ et $\varphi^-$ les param\`etres de Langlands des paquets $\Pi$ et $\Pi^-$, on calcule facilement 
$$Jac^{\theta}_{\rho_{1},...,\rho_{m}}(\pi(\varphi))=\pi(\varphi^-).$$
Alors la m\^eme preuve que celle de (2) entra\^{\i}ne (4). $\square$

\bigskip

\subsection{Support cuspidal \'etendu}

Dans la suite de la section, on suppose $G$ symplectique ou sp\'ecial orthogonal impair. On indiquera dans le dernier paragraphe 2.15 comment adapter les arguments au cas d'un groupe sp\'ecial orthogonal pair.

   Soit $\pi$ une repr\'esentation irr\'eductible de $G(F)$.  On appelle support cuspidal \'etendu de $\pi$ le support cuspidal ordinaire de $\pi^{GL}$.  C'est un ensemble avec multiplicit\'es de repr\'esentations  irr\'eductibles cuspidales non n\'ecessairement unitaires de groupes lin\'eaires. On remarque que par sa construction, cet ensemble est stable par passage \`a la  contragr\'ediente.

\ass{Lemme}{Soit $\pi$ une repr\'esentation de $G(F)$ de la forme $\pi=\sigma\times \pi'$, o\`u $\sigma$ est une repr\'esentation irr\'eductible d'un groupe lin\'eaire et $\pi'$ une repr\'esentation irr\'eductible d'un groupe de m\^eme type que $G$. Alors tout sous-quotient irr\'eductible de  $\pi$  a  pour support cuspidal \'etendu l'union disjointe du support cuspidal \'etendu de $\pi'$ et des supports cuspidaux ordinaires de  $\sigma$ et $\check{\sigma}$.}

Ceci permet de parler du support cuspidal \'etendu d'une repr\'esentation induite m\^eme si celle-ci n'est pas irr\'eductible: c'est le support cuspidal \'etendu de n'importe lequel de ses sous-quotients irr\'eductibles.

Preuve. Soit $\pi$  une repr\'esentation irr\'eductible de $G(F)$. On peut la r\'ealiser comme  sous-quotient  d'une induite
$$(\times_{i=1,.. .,v}\rho_{i})\times \pi_{cusp},$$
o\`u $\pi_{cusp }$ est une repr\'esentation irr\'eductible cuspidale d'un groupe de m\^eme type que $G$ et, pour tout $i=1,...,v$, $\rho_{i}$ est une repr\'esentation irr\'eductible cuspidale d'un groupe lin\'eaire, pas forc\'ement unitaire. Appelons support cuspidal ordinaire de $\pi$ l'ensemble avec multiplicit\'es $\{\rho_{i};i=1,...,v\}\sqcup\{\check{\rho}_{i};i=1,...,v\}\sqcup\{\pi_{cusp}\}$. On sait qu'il ne d\'epend pas de l'induite (1) choisie. Il est form\'e de repr\'esentations de groupes lin\'eaires et d'au plus une repr\'esentation $\pi_{cusp}$ d'un groupe de m\^eme type que $G$ (cette repr\'esentation dispara\^{\i}t ici et dans la suite quand ce groupe est r\'eduit \`a $1$). On va montrer

(1)  le support cuspidal \'etendu de $\pi$ est l'union disjointe de celui de $\pi_{cusp}$ et du compl\'ementaire de $\{\pi_{cusp}\}$ dans le support cuspidal ordinaire de $\pi$.

Cette assertion entra\^{\i}ne le lemme car, dans sous les hypoth\`eses de  l'\'enonc\'e, le support cuspidal ordinaire de tout sous-quotient irr\'eductible de $\sigma\times\pi'$ est l'union disjointe de celui de $\pi'$ et de ceux de $\sigma$ et $\check{\sigma}$.

Pour prouver (1), on utilise la remarque suivante:

(2) supposons que $\pi$ apparaisse comme sous-quotient d'une induite $\sigma_{1}\times...\times \sigma_{v}\times \pi'$, o\`u $\pi'$ et les $\sigma_{i}$ sont irr\'eductibles et o\`u $v\geq1$, et que l'on sache que le support cuspidal \'etendu soit r\'eunion de celui de $\pi'$ et des supports cuspidaux ordinaires des $\sigma_{i}$ et des $\check{\sigma}_{i}$ pour $i=1,...,v$; alors (1) est vrai pour $\pi$.

 En effet, le support cuspidal ordinaire de $\pi$ est forc\'ement r\'eunion de celui de $\pi'$ et de ceux des $\sigma_{i}$ et $\check{\sigma}_{i}$ pour $i=1,...,v$. En raisonnant par r\'ecurrence sur $d_{G}$, on peut supposer que le support cuspidal \'etendu de $\pi'$ se d\'eduit de son support cuspidal \'etendu de la fa\c{c}on prescrite par (1) et on en d\'eduit que le support cuspidal \'etendu de $\pi$ se d\'eduit de la m\^eme fa\c{c}on de son support cuspidal ordinaire.
 
 Si $\pi$ n'est pas temp\'er\'ee, on r\'ealise $\pi$ comme quotient de Langlands d'une induite comme en (2), o\`u $\pi'$ est temp\'er\'ee et les $\sigma_{i}$ sont des repr\'esentations temp\'er\'ees tordues par un caract\`ere. Par d\'efinition du support cuspidal \'etendu de $\pi$, les hypoth\`eses  de (2) sont satisfaites. Donc (1) est v\'erifi\'ee pour $\pi$.
 
 Supposons que $\pi$ soit temp\'er\'ee et qu'il existe un \'el\'ement $(\rho,a)$ de $Jord(\pi)$ tel que, ou bien $(\rho,a)$ ne soit pas la bonne parit\'e, ou bien $(\rho,a)$ soit de bonne parit\'e mais intervienne avec multiplicit\'e au moins $2$. Alors 2.2(1) permet de r\'ealiser $\pi$ comme sous-quotient d'une induite  de sorte que les hypoth\`eses de (2) soient satisfaites. D'o\`u la conclusion dans ce cas.
 
 Supposons que $\pi$ soit temp\'er\'ee, que $Jord(\pi)$ soit form\'e de couples $(\rho,a)$ de bonne parit\'e intervenant avec multiplicit\'e $1$, et que $\pi$ ne soit pas cuspidale. Cette derni\`ere hypoth\`ese entra\^{\i}ne que l'on peut r\'ealiser $\pi$ comme sous-module d'un induite $\rho'\times \pi'$, o\`u $\pi'$ est irr\'eductible et $\rho'$ est irr\'eductible et cuspidale, pas forc\'ement unitaire.  D'apr\`es l'hypoth\`ese sur $Jord(\pi)$, 2.2(2) entra\^{\i}ne qu'il existe $(\rho,a)\in Jord(\pi)$ tel que $a\geq2$ et $\rho'=\rho\vert .\vert _{F}^{(a-1)/2}$. La relation 2.2(3) entra\^{\i}ne alors que le support cuspidal \'etendu de $\pi$ est r\'eunion de celui de $\pi'$ et de $\{\rho',\check{\rho}'\}$. Autrement dit, les hypoth\`eses de (2) sont v\'erifi\'ees et on conclut.
 
 On est ramen\'e au cas o\`u $\pi$ est cuspidale. Mais alors (1) est tautologique. $\square$

{\bf Remarque.} Soient $\pi$ une repr\'esentation irr\'eductible temp\'er\'ee de $G(F)$, $\rho$ une repr\'esentation  irr\'eductible cuspidale et unitaire d'un groupe lin\'eaire et $s\in {\mathbb R}$. Supposons que l'induite $\rho\vert.\vert_{F}^{s}\times \pi$ ait un sous-quotient irr\'eductible ayant m\^eme support cuspidal qu'une repr\'esentation temp\'er\'ee. Alors l'une des conditions suivantes est v\'erifi\'ee:

$\bullet$ $s=0$;

$\bullet$  $s$ est un demi-entier avec $\vert s\vert \geq1$, $\rho\simeq \check{\rho}$ et $Jord(\pi)$ contient le couple $(\rho,2\vert s\vert -1)$ 

$\bullet$  $s= \pm 1/2$ et  $(\rho,2)$ a bonne parit\'e.

Cela r\'esulte du lemme ci-dessus et des propri\'et\'es du support cuspidal \'etendu d'une repr\'esentation temp\'er\'ee.

\bigskip

\subsection{Support cuspidal \'etendu et induction dans le cas temp\'er\'e}
 
 \ass{Lemme}{Soit $\pi$ une repr\'esentation irr\'eductible temp\'er\'ee de $G(F)$. On suppose qu'elle est sous-quotient d'une induite 
 $$<e,f>_{\rho}\times \pi'$$
 o\`u $[e,f]$ est un segment tel que $e\geq0\geq f$, $\rho$ est une rep\'esentation irr\'eductible cuspidale et unitaire  d'un groupe lin\'eaire  et $\pi'$ est une repr\'esentation irr\'eductible temp\'er\'ee d'un groupe de m\^eme type que $G$ .
Alors  l'une des propri\'et\'es suivantes est v\'erifi\'ee:

$\bullet$ $\rho\simeq \check{\rho}$, $e$ et $f$ sont des demi-entiers et $Jord(\pi)=Jord(\pi')\cup\{(\rho,2e+1), (\rho,-2f+1)\}$;

$\bullet$  $\rho\not\simeq \check{\rho}$,  $e=-f$ et $Jord(\pi)=Jord(\pi')\cup\{(\rho,2e+1), (\check{\rho},2e+1)\}$.}

Preuve. Le support cuspidal \'etendu de $\pi$ est, d'apr\`es le lemme pr\'ec\'edent, l'union de celui de $\pi'$ avec $\cup_{x\in [e,f]}\{\rho\vert.\vert_{F}^{x}\}\cup_{x\in [-f,-e]}\{\check{\rho}\vert.\vert_{F}^{x}\}$. Les conditions sur $\rho,e,f$ r\'esultent des propri\'et\'es rappel\'ees en 2.2 du support cuspidal \'etendu des repr\'esentations temp\'er\'ees. $\square$

\bigskip

\subsection{ Un lemme technique}

\ass{Lemme}{Consid\'erons une induite
$$\sigma\vert .\vert ^s\times \pi',$$
o\`u  $\pi'$  est une repr\'esentation irr\'eductible temp\'er\'ee d'un groupe de m\^eme type que $G$, $\sigma=St(\rho,a)$ est une repr\'esentation de Steinberg g\'en\'eralis\'ee temp\'er\'ee d'un groupe lin\'eaire et $s >0$ est un r\'eel.  Alors le quotient de Langlands de  cette induite  est l'unique sous-quotient irr\'eductible   qui poss\`ede une inclusion dans une induite de la forme 
$$
(\times_{x\in [(a-1)/2,-(a-1)/2]}\check{\rho}\vert.\vert_{F}^{x-s})\times \tau,
$$
o\`u $\tau$ est une repr\'esentation non n\'ecessairement irr\'eductible d'un groupe de m\^eme type que $G$.
  En particulier si
$$
\sigma\vert.\vert_{F}^s\times \pi' \hookrightarrow (\times_{x\in [(a-1)/2,-(a-1)/2]}\check{\rho}\vert.\vert_{F}^{x-s})\times \tau$$
avec $\tau$  comme ci-dessus, l'induite $\sigma\vert.\vert_{F}^s\times \pi'$ est irr\'educ\-tible.}

Preuve.   Il est clair que le quotient de Langlands de l'induite $\sigma\vert.\vert_{F}^{s}\times \pi'$ a la propri\'et\'e requise car il est inclus dans $\check{\sigma}\vert.\vert_{F}^{-s}\times \pi'$ et  cette induite est elle-m\^eme incluse dans $(\times_{x\in [(a-1)/2,-(a-1)/2]}\check{\rho}\vert.\vert_{F}^{x-s})\times\pi'$. Soit donc $\pi$ un sous-quotient irr\'eductible de l'induite ayant une inclusion comme dans l'\'enonc\'e. Comme $\pi$ est irr\'eductible, une telle inclusion en donne une de m\^eme type mais avec $\tau$ irr\'eductible. On suppose donc que $\tau$ est irr\'educ\-tible.
Par r\'eciprocit\'e de Frobenius un module de Jacquet convenable de $\pi'$ admet $(\otimes_{x\in [(a-1)/2,-(a-1)/2]}\check{\rho}\vert.\vert_{F}^{x-s})\otimes \tau$ comme quotient.

On sait calculer tous les termes du module de Jacquet de l'induite $\sigma\vert.\vert_{F}^{s}\times \pi'$. Les termes   d'un module de Jacquet cuspidal de cette induite peuvent se regrouper en sous-ensembles param\'etr\'es par un demi-entier entier $x\in [(a+1)/2,-(a-1)/2]$ et un terme cuspidal, $(\otimes_{j=1,...,k}\rho_{j}\vert.\vert_{F}^{y_{j}})\times \pi'_{cusp}$ du module de Jacquet de $\pi'$. Les termes correspondants, \'ecrits sous la forme $$
(\otimes_{(\rho',z')\in {\mathcal E} }\rho'\vert.\vert_{F}^{z'})\otimes \pi'_{cusp}
$$
sont tels que le l'ensemble ordonn\'e ${\mathcal E}$
s'obtient en m\'elangeant les ensembles ordonn\'es \'ecrits ci-dessous sans permuter l'ordre interne \`a chacun de ces sous-ensembles:
$$
\{\rho\vert.\vert_{F}^{s+i}; i\in [(a-1)/2,x]\} \hbox{(c'est l'ensemble vide si $x=(a+1)/2$ )};$$
$$ \{ \check{\rho}\vert.\vert_{F}^{i'-s};i'\in [(a-1)/2,-x[\} \hbox{ (c'est l'ensemble vide si $x=-(a-1)/2$) };$$
$$ \{\rho_{j}\vert.\vert_{F}^{y_{j}};j=1,...,k\}.$$
Ainsi les $a$ premiers exposants d'un tel terme, \'ecrits $z'_{i}$ pour $i=1,...,a$ se d\'ecomposent en $a_{1}$ termes du premier ensemble, $a_{2}$ termes du deuxi\`eme et $a_{3}$ termes du troisi\`eme. Donc on a:
$$
\sum_{i=1,...,a}z'_{i}=\sum_{\ell\in [(a-1)/2,(a+1)/2-a_{1}]} (s+\ell) +\sum_{\ell\in [(a-1)/2,(a+1)/2-a_{2}]}(-s+\ell) +\sum_{j=1,...,a_3}y_{j}.
$$
Puisque $\pi'$ est une repr\'esentation temp\'er\'ee, on a s\^urement $\sum_{j=1,...,a_3}y_{j}\geq 0$ et donc
$$\sum_{i=1,...,a}z'_{i}
\geq a_{1}s-a_{2}s \geq -a_{2}s.
$$
Pour le terme consid\'er\'e du module de Jacquet de $\pi$, cette somme vaut $-as$. Comme $a_{2}\leq a$, on doit avoir \'egalit\'e $a_{2}=v$, d'o\`u $a_{1}=a_{3}=0$. Notons $P$ le sous-groupe parabolique qui sert \`a d\'efinir l'induite $\sigma\vert .\vert _{F}^s\times \pi'$. Les termes  du module de Jacquet de l'induite $\sigma\vert .\vert _{F}^s\times \pi'$ qui v\'erifient les conditions pr\'ec\'edentes  sont ceux qui proviennent du sous-quotient $\check{\sigma}\vert .\vert _{F}^{-s}\otimes \pi'$  du module de Jacquet $(\sigma\vert .\vert _{F}^s\times \pi')_{P}$. Ils n'interviennent que dans le quotient de Langlands de notre induite.  Donc $\pi$ est ce quotient de Langlands. Cela d\'emontre la premi\`ere assertion de l'\'enonc\'e. Sous l'hypoth\`ese de la seconde assertion, l'unique sous-module irr\'eductible de $\sigma\vert .\vert _{F}^s\times \pi'$ est le quotient de Langlands d'apr\`es ce que l'on vient de prouver. Puisque l'unique sous-module irr\'eductible et aussi l'unique quotient irr\'eductible, la repr\'esentation est irr\'eductible. $\square$

\bigskip

\subsection{Induction et $L$-paquets temp\'er\'es}
 
  Ci-dessous, on utilise \`a plusieurs reprises la remarque \'el\'ementaire suivante. Soit $x$ un nombre r\'eel non nul et $\rho$ une repr\'esentation  irr\'eductible cuspidale et unitaire d'un groupe lin\'eaire. Soient aussi $m\geq1$ un entier et $\sigma$ une repr\'esentation de $G(F)$. On suppose qu'il n'existe pas d'inclusion de la forme $\sigma\hookrightarrow \rho\vert.\vert_{F}^x\times \sigma'$ o\`u $\sigma'$ est une repr\'esentation quelconque. Alors l'induite
$
\rho\vert.\vert_{F}^{x}\times \cdots \rho\vert.\vert_{F}^{x}\times \sigma,
$, o\`u il y a $m$ copies de $\rho\vert.\vert_{F}^{x}$,
a  un unique sous-module irr\'eductible. C'est un calcul de module de Jacquet et de r\'eciprocit\'e de Frobenius pour le parabolique standard de sous-groupe de Levi $GL(md_{\rho})\times G$. En effet, les sous-quotients irr\'eductibles de ce module de Jacquet de la forme $\tau \otimes \tau'$ v\'erifient soit que le support cuspidal de $\tau$ n'est pas $m$ copies de $\rho\vert.\vert_{F}^x$, soit que $\tau$ est l'induite irr\'eductible $\rho\vert.\vert_{F}^x\times \cdots \times \rho\vert.\vert_{F}^x$. Dans ce dernier cas, on a n\'ecessairement $\tau'\simeq \sigma$ et un tel terme n'intervient qu'avec multiplicit\'e au plus 1 comme sous-quotient irr\'eductible du module de Jacquet. Par r\'eciprocit\'e de Frobenius, tout sous-module irr\'eductible de l'induite \'ecrite a son module de Jacquet qui admet $(\rho\vert.\vert_{F}^x\times \cdots\times \rho\vert.\vert_{F}^{x})\otimes \sigma$ comme quotient irr\'eductible. Comme le foncteur de Jacquet est exact, cela force  l'unicit\'e d'un tel sous-module irr\'eductible. Remarquons que, d'apr\`es 2.1(2),  l'hypoth\`ese sur $\sigma$ est v\'erifi\'ee si les deux conditions suivantes le sont:

$\bullet$ $\sigma$ est une repr\'esentation irr\'eductible temp\'er\'ee;

$\bullet$ $x$ n'est pas un demi-entier positif,  ou $x$ est un tel demi-entier mais $Jord(\sigma)$ ne contient pas $(\rho,2x+1)$.

 \ass{Lemme}{Soit $\Pi$ un paquet de repr\'esentations temp\'er\'ees. Soit $x>0$ un demi-entier tel que $(\rho,2x+1)\in Jord(\Pi)$. On note $m$ la multiplicit\'e de $(\rho,2x+1)$ dans $Jord(\Pi)$ et $\Pi^-$ le paquet de repr\'esentations temp\'er\'ees qui se d\'eduit de $\Pi$ en rempla\c{c}ant les $m$ copies de $(\rho,2x+1)$ par $m$ copies de $(\rho,2x-1)$. Pour tout $\pi^-\in \Pi^-$, on note $\pi$ l'unique sous-module irr\'eductible de l'induite $\rho\vert.\vert_{F}^x\times \cdots \rho\vert.\vert_{F}^{x}\times \pi^-$. Alors $\pi$ est une repr\'esentation temp\'er\'ee de $\Pi$ et l'application ainsi d\'efinie de $\Pi^-$ dans $\Pi$ est une injection.}

Preuve. Soit $\pi\in \Pi$. Supposons qu'il existe une inclusion de la forme $\pi\hookrightarrow \rho\vert.\vert_{F}^{x}\times \cdots \times \rho\vert.\vert_{F}^{x}\times \sigma$, avec $m$ copies de $\rho\vert .\vert _{F}^x$. On v\'erifie que l'on n'a certainement aucune inclusion de la forme $\sigma\hookrightarrow \rho\vert.\vert_{F}^{x}\times \sigma'$. Sinon, on aurait une inclusion $\pi\hookrightarrow \rho\vert.\vert_{F}^{x}\times \cdots \times \rho\vert.\vert_{F}^{x}\times \sigma$, avec $m+1$ copies de $\rho\vert .\vert _{F}^x$ et $Jord(\pi)$ contiendrait  $m+1$ copies de $(\rho,a)$, cf. 2.1(3).  Les calculs de modules de Jacquet expliqu\'es ci-dessus montrent alors que si une telle inclusion se produit avec $\sigma$ que l'on suppose irr\'eductible,  tout sous-quotient irr\'eductible du module de Jacquet de $\pi$ de la forme $\rho\vert.\vert_{F}^{x}\times \cdots \rho\vert.\vert_{F}^{x}\otimes \sigma''$, o\`u il y a $m$ copies de $\rho\vert.\vert_{F}^{x}$, v\'erifie $\sigma''\simeq \sigma$.    Un peu plus g\'en\'eralement supposons que $\pi\in \Pi$ et que le module de Jacquet de $\pi$ contient un sous-quotient irr\'eductible de la forme $\rho\vert.\vert_{F}^{x}\times \cdots \rho\vert.\vert_{F}^{x}\otimes \sigma''$; alors quitte \`a changer $\sigma''$ on peut supposer que ce sous-quotient est en fait un quotient du module de Jacquet de $\pi$ et, par r\'eciprocit\'e de Frobenius que $\pi$ est un sous-module comme ci-dessus.

Quand on applique le module de Jacquet pour le parabolique $GL(md_{\rho},F)\otimes G'$ (pour $G'$ convenable) \`a la distribution stable form\'ee de la somme des \'el\'ements de $\pi$ et que l'on projette sur le caract\`ere de la repr\'esentation $\rho\vert.\vert_{F}^x\times \cdots \times \rho\vert.\vert_{F}^x$ (avec $m$-copies), on obtient une distribution stable associ\'ee \`a $\Pi^-$. Ainsi d'apr\`es ce que l'on a vu ci-dessus, une repr\'esentation $\pi\in \Pi$ soit dispara\^{\i}t dans cette proc\'edure, soit contribue par le caract\`ere de $\sigma$ (avec les notations ci-dessus). Ainsi n\'ecessairement $\sigma\in \Pi^-$ et   toute repr\'esentation de $\Pi^-$ est obtenue par cette proc\'edure.  On a donc d\'efini sur un sous-ensemble de $\Pi$ un inverse surjectif de l'application de l'\'enonc\'e. Cela montre que cette derni\`ere application    est bien une injection de $\Pi^-$ dans $\Pi$. $\square$

\bigskip

\subsection{ Propri\'et\'e des repr\'esentations temp\'er\'ees ayant un mod\`ele de Whittaker} 

  Supposons $G$ quasi-d\'eploy\'e. On d\'efinit de la fa\c{c}on habituelle la notion de mod\`ele de Whittaker. Il y a plusieurs types de tels mod\`eles, autant que de classes de conjugaison d'\'el\'ements unipotents r\'eguliers dans $G(F)$. On utilise la propri\'et\'e suivante:
\bigskip

{\it soit $\pi$ une repr\'esentation irr\'eductible temp\'er\'ee de $G(F)$; il existe une unique repr\'esentation temp\'er\'ee irr\'eductible $\pi_{0}$ de $G(F)$ ayant un mod\`ele de Whittaker d'un type fix\'e et telle que $Jord(\pi_{0})=Jord(\pi)$.}

Cf. [K] th\'eor\`eme 3.4, [W1] th\'eor\`eme 4.9.

 \ass{Lemme}{Supposons $G$ quasi-d\'eploy\'e, soit $\pi$ une repr\'esentation irr\'eductible de $G(F)$ ayant un mod\`ele de Whittaker. On suppose que le support cuspidal \'etendu de $\pi$ est celui d'une repr\'esentation temp\'er\'ee. Alors $\pi$ est temp\'er\'ee.}

Preuve. On d\'emontre ce lemme par r\'ecurrence sur le rang de $G$ et pour cela on a besoin d'une cons\'equence du lemme. Pour fixer les notations, on note $\Pi$ un paquet de repr\'esentations temp\'er\'ees et  on suppose que le support cuspidal \'etendu de $\pi$ est le m\^eme que celui des \'el\'ements de $\Pi$. Le lemme dit alors que $\pi$ appartient \`a $\Pi$.  Ceci montre que le support cuspidal ordinaire de $\pi$ est  bien d\'etermin\'e. On peut d\'ecrire  ce support cuspidal. En effet, notons $\varphi$ le morphisme de $W_{DF}$ dans le $L$-groupe de $G$ param\'etrisant $\Pi$. Notons  $\varphi_{L}$ le morphisme de $W_{DF}$ dans ce $L$-groupe qui est trivial sur $SL(2,{\mathbb C})$ et qui, sur $W_{F}$, est le compos\'e de l'inclusion:
$
W_{F}\hookrightarrow W_{DF}=W_{F}\times SL(2,{\mathbb C}); w\mapsto \biggl(w, \begin{pmatrix}\vert w\vert_{F}^{1/2} &0\\0 &\vert w\vert_{F}^{-1/2}\end{pmatrix}\biggr)$
avec $\varphi$. Ainsi $\varphi_{L}$ d\'etermine un paquet de Langlands assez particulier. Le support cuspidal \'etendu de tout \'el\'ement de ce paquet est le m\^eme que celui des \'el\'ements de $\Pi$.
Il existe exactement une repr\'esentation de ce paquet qui soit le quotient de Langlands d'une induite ayant un mod\`ele de Whittaker. Notons $\pi_{0}$ le sous-quotient irr\'eductible de cette induite ayant un mod\`ele de Whittaker. Le lemme dit que $\pi_{0}$ est temp\'er\'ee et isomorphe \`a $\pi$. On pourrait montrer que $\pi=\pi_{0}$ est la repr\'esentation duale au sens d'Aubert, Schneider-Stuhler du quotient de Langlands. Ce qui nous importe est que si $\pi$ est cuspidal, alors $\varphi^L=\varphi$, donc $Jord(\pi)$ est un ensemble  de couples $(\rho,1)$.

Prouvons le lemme. On fixe $\pi$  comme dans l'\'enonc\'e et    $\Pi$ comme ci-dessus. Si $\pi$ est cuspidale, elle est a fortiori temp\'er\'ee, il n'y a rien \`a d\'emontrer. Sinon, on \'ecrit $\pi$ comme sous-module d'une induite de repr\'esentations cuspidales
$$
\pi\hookrightarrow (\times_{i=1,...,v}\rho_{i}\vert.\vert_{F}^{s_{i}})\times \pi_{cusp}, \eqno(1)
$$
o\`u $\pi_{cusp}$ est une repr\'esentation irr\'eductible cuspidale d'un groupe de m\^eme type que $G$, $v\geq1$ est un entier  et, pour tout $i=1,...,v$, $\rho_{i}$ est une repr\'esentation cuspidale unitaire irr\'eductible d'un  groupe lin\'eaire et $s_{i}$ est un nombre r\'eel. Comme $\pi_{cusp}$ a n\'ecessairement un mod\`ele de Whittaker et est une repr\'esentation d'un groupe de rang plus petit que $G$, on sait, par r\'ecurrence, que $Jord(\pi_{cusp})$ est un ensemble de couples $(\rho,1)$. Puisque $\pi_{cusp}$ est de la s\'erie discr\`ete, cet ensemble est sans multiplicit\'es et pour tout $(\rho,1)$ y intervenant, $(\rho,1)$ est de bonne parit\'e,  en particulier $\rho$ est autodual.

 Soit $(\rho,a)\in Jord(\Pi)$. On suppose d'abord que $\rho$ n'est pas autoduale. Ainsi $(\check{\rho},a)\in Jord(\Pi)$ et il existe un sous-ensemble ${\mathcal E}$ de $\{1,...,v\}$ tel que $\cup_{i\in {\mathcal E}}\{\rho_{i}\vert.\vert_{F}^{s_{i}},\check{\rho}_{i}\vert.\vert_{F}^{-s_{i}}\}=\cup_{x\in [(a-1)/2,-(a-1)/2]}\{\rho\vert.\vert_{F}^{x},\check{\rho}\vert.\vert_{F}^{x}\}$. Un tel ensemble n'est pas unique, en g\'en\'eral, et on en fixe un. Ainsi $\pi$ est certainement un sous-quotient irr\'eductible de l'induite:
 $$
( \times_{x\in [(a-1)/2,-(a-1)/2]}\rho\vert.\vert_{F}^x )\times(\times_{i=1,...,v; i\notin {\mathcal E}}\rho_{i}\vert.\vert_{F}^{s_{i}})\times \pi_{cusp}.  
 $$
Il existe donc   un sous-quotient irr\'eductible $\sigma$ de l'induite $ \times_{x\in [(a-1)/2,-(a-1)/2]}\rho\vert.\vert_{F}^x $ et  un sous-quotient irr\'eductible $\pi'$ de l'induite $(\times_{i=1,...,v; i\notin {\mathcal E}}\rho_{i}\vert.\vert_{F}^{s_{i}})\times \pi_{cusp}$ tel que $\pi$ soit un sous-quotient irr\'eductible de l'induite $\sigma\times \pi'$. N\'ecessairement $\sigma$ a un mod\`ele de Whittaker au sens usuel et $\pi'$ en a un du m\^eme type que $\pi$. Ainsi $\sigma\simeq <(a-1)/2,-(a-1)/2>_{\rho}$. En fait, on conna\^{\i}t aussi $\pi'$. En effet
 le support cuspidal \'etendu de l'induite $( \times_{i=1,...,v; i\notin {\mathcal E}}\rho_{i}\vert.\vert_{F}^{s_{i}})\times \pi_{cusp}$ est celui des repr\'esentations temp\'er\'ees dans le paquet $\Pi'$, qui se d\'eduit de $\Pi$ en enlevant  $(\rho,a)$ et $(\check{\rho},a)$. En appliquant le lemme par r\'ecurrence, on sait que $\pi'$ est une repr\'esentation temp\'er\'ee. Il en est donc de m\^eme de   tout sous-quotient de $\sigma\times \pi'$, donc $\pi$  est temp\'er\'ee. Cela prouve le lemme dans ce cas.
 
Supposons  maintenant que  $Jord(\Pi)$ contient un \'el\'ement $(\rho,a)$ tel que $\rho$ soit autoduale,  avec une multiplicit\'e not\'ee $m$ v\'erifiant l'une des conditions suivantes:

$\bullet$  $m>2$;

$\bullet$ $m=2$ et  $a$ est pair;

$\bullet$  $m=2$ et $a$ est impair mais $(\rho,1)$ n'intervient pas dans $Jord(\pi_{cusp})$.

Ces conditions sont par exemple automatiques si   $(\rho,a)$ est de mauvaise parit\'e. Alors ${\mathcal E}$ comme ci-dessus existe encore et $\pi'$ d\'efini comme ci-dessus a encore son support cuspidal \'etendu qui se d\'eduit de celui de $\pi$ en enlevant  deux copies de $(\rho,a)$. On conclut comme ci-dessus.

 D'autre part, si $a=1$ pour tout $(\rho,a)\in Jord(\Pi)$, les $s_{i}$ sont tous nuls et $\sigma$ est certainement une repr\'esentation temp\'er\'ee d'apr\`es l'inclusion (1). 
Il nous suffit donc maintenant de d\'emontrer le lemme dans le cas o\`u $Jord(\Pi)$ contient un \'el\'ement $(\rho,a)$ de bonne parit\'e avec $a\geq2$. On  peut se limiter au cas o\`u la multiplicit\'e de $(\rho,a)$ est inf\'erieure ou \'egale \`a $2$. On fixe un tel $(\rho,a)$ et on suppose que $a$ est maximal avec cette propri\'et\'e.
 
 On fait d'abord la d\'emonstration dans le cas o\`u la multiplicit\'e de $(\rho,a)$ dans $Jord(\pi)$ est $1$ pour clarifier la m\'ethode.
 On fixe une inclusion (1) et on sait qu'il existe $i_{0}\in \{1,...,v\}$ tel que $\rho_{i_{0}}\vert.\vert_{F}^{s_{i_{0}}}\simeq \rho\vert.\vert_{F}^{\zeta(a-1)/2}$ pour  un signe convenable $\zeta$. On fixe une telle inclusion avec l'hypoth\`ese suppl\'ementaire que pour tous les choix possibles, $i_{0}$ est minimum. On suppose d'abord que $i_{0}=1$ et on conclut: il existe un sous-quotient $\pi'$ de  l'induite $(\times_{i=2,...,v}\rho_{i}\vert.\vert_{F}^{s_{i}})\times \pi_{cusp}$ tel que $\pi$ soit un sous-module irr\'eductible de l'induite $\rho\vert.\vert_{F}^{\zeta(a-1)/2}\times \pi'$. Le support cuspidal \'etendu de $\pi'$ est celui des repr\'esentations temp\'er\'ees dans le paquet qui se d\'eduit de $ \Pi$ en rempla\c{c}ant $(\rho,a)$ par $(\rho,a-2)$. De plus $\pi'$ admet un mod\`ele de Whittaker. Par l'hypoth\`ese de r\'ecurrence, on sait que $\pi'$ est temp\'er\'ee. Ainsi si $\zeta=+$, $\pi$ est un sous-module irr\'eductible de l'induite $\rho\vert.\vert_{F}^{(a-1)/2}\times \pi'$. On applique le lemme 2.6 avec ici $m=1$, $\pi^-=\pi'$. Ce lemme nous dit que $\pi$ appartient \`a $\Pi$, donc est une repr\'esentation temp\'er\'ee. Consid\'erons le cas o\`u $\zeta=-$. Sous cette hypoth\`ese, $\pi$ est le sous-module de Langlands de l'induite $\rho\vert.\vert_{F}^{-(a-1)/2}\times \pi'$. Comme $\pi$ a un mod\`ele de Whittaker, le th\'eor\`eme 1.1 de [Mu] montre que  cette induite est irr\'eductible et $\pi$ est donc aussi un sous-module irr\'eductible de l'induite $\rho\vert.\vert_{F}^{(a-1)/2}\times \pi'$. On conclut comme pr\'ec\'edemment. Au passage d'ailleurs cette conclusion montre que l'induite ne peut pas \^etre irr\'eductible et que ce cas ne peut se produire.

 Il suffit donc  de d\'emontrer que $i_{0}=1$. Supposons  qu'il n'en soit pas ainsi mais que $\zeta=+$. Dans ce cas $\rho_{i_{0}-1}\vert.\vert_{F}^{s_{i_{0}-1}}\times \rho\vert.\vert_{F}^{(a-1)/2}$ est soit irr\'eductible soit de longueur deux. Dans le premier cas on peut commuter les facteurs ce qui contredit la minimalit\'e de $i_{0}$. Dans le deuxi\`eme cas, $\rho_{i_{0}-1}\vert.\vert_{F}^{s_{i_{0}-1}}\simeq \rho\vert.\vert_{F}^{(a-3)/2}$ et  on peut remplacer $\rho\vert.\vert_{F}^{(a-3)/2}\times \rho\vert.\vert_{F}^{(a-1)/2}$ par son unique sous-quotient ayant un mod\`ele de Whittaker. Ce sous-quotient est  $<(a-1)/2,(a-3)/2>_{\rho}$. Or  $<(a-1)/2,(a-3)/2>_{\rho}\hookrightarrow \rho\vert.\vert_{F}^{(a-1)/2}\times \rho\vert.\vert_{F}^{(a-3)/2}$ et on peut donc encore \'echanger $i_{0}$ et $i_{0}-1$ ce qui contredit la minimalit\'e de $i_{0}$. Donc si $\zeta=+$, $i_{0}=1$.
 
 On suppose que $\zeta=-$. Ici la deuxi\`eme partie de l'argument ci-dessus est diff\'erente. On peut remplacer $\rho\vert.\vert_{F}^{-(a-3)/2}\times \rho\vert.\vert_{F}^{-(a-1)/2}$ par $<-(a-3)/2,-(a-1)/2>_{\rho}$. En proc\'edant ainsi de proche ne proche, on montre qu'il existe un segment d\'ecroissant $[e,-(a-1)/2]$ et une inclusion 
 $$
 \pi\hookrightarrow <e,-(a-1)/2>_{\rho}\times \pi',\eqno(2)
 $$
 o\`u $\pi'$ est une repr\'esentation  irr\'eductible convenable. On va montrer que n\'ecessairement $e=(a-3)/2$. De (1) on tire  que $\pi$ est un sous-quotient irr\'eductible de l'induite $\rho\vert.\vert_{F}^{(a-1)/2}\times \pi_{1}$ o\`u $\pi_{1}$ a pour support cuspidal \'etendu celui de $\pi$ o\`u on a remplac\'e $(\rho,a)$ par $(\rho,a-2)$.  Comme $\pi_{1}$ a n\'ecessairement un mod\`ele de Whittaker, $\pi_{1}$ est temp\'er\'ee. On sait calculer les modules de Jacquet de l'induite $\rho\vert.\vert_{F}^{(a-1)/2}\times \pi_{1}$.  On voit que (2) force l'existence d'une inclusion:
 $$
 \pi_{1}\hookrightarrow \times_{j\in [e,-(a-3)/2]}\rho\vert.\vert_{F}^j\times \pi_{2},
 $$
 o\`u $\pi_{2}$ ne nous int\'eresse pas. Comme $\pi_{1}$ est une repr\'esentation temp\'er\'ee, n\'ecessairement $e=(a-3)/2$ comme annonc\'e.  
Le support cuspidal \'etendu de $\pi_{1}$ contient donc  
$$\cup_{x\in [(a-3)/2,-(a-3)/2}\{\rho\vert.\vert_{F}^{x}\}$$
 avec multiplicit\'e au moins $2$. Ceci est donc aussi vrai pour le support cuspidal de $\pi$ et  il existe des entiers $b,b'$ sup\'erieurs ou \'egaux \`a $a-2$ tel que $(\rho,b)$ et $(\rho,b')$ soient dans $Jord(\Pi)$. Par maximalit\'e de $a$, on a $b,b'\leq a$ et comme $(\rho,a)$ n'intervient qu'avec multiplicit\'e $1$ dans $Jord(\Pi)$, l'un des deux vaut $a-2$. Ainsi $Jord(\Pi)$ contient $(\rho,a)$ et $(\rho,a-2)$ et le support cuspidal \'etendu de $\pi'$ se d\'eduit de celui de $Jord(\pi)$ en enlevant $(\rho,a)$ et $(\rho,a-2)$. On applique l'hypoth\`ese de r\'ecurrence \`a $\pi'$ qui a n\'ecessairement un mod\`ele de Whittaker et on sait que $\pi'$ est temp\'er\'ee. Ainsi $\pi$ est le sous-module de Langlands de l'induite $<(a-3)/2,-(a-1)/2>_{\rho}\times \pi'$. D'apr\`es [Mu] théorème 1.1, cette induite doit \^etre irr\'eductible, elle est donc isomorphe \`a $<(a-1)/2,-(a-3)/2>_{\rho}\times \pi'$  on trouve encore une inclusion avec $i_{0}=1$ ce qui contredit la minimalit\'e de $i_{0}$.

 Il nous reste donc \`a voir le cas o\`u la multiplicit\'e de $(\rho,a)$ dans $Jord(\pi)$ est   \'egale \`a $2$. On proc\`ede comme ci-dessus mais ici, $i_{0}$ et remplac\'e par $i_{1}< i_{2}$ o\`u pour $j=1,2 $, il existe un signe $\zeta_{j}$ tel que $\rho_{i_{j}}\vert.\vert_{F}^{s_{i_{j}}}=\rho\vert.\vert_{F}^{\zeta_{j}(a-1)/2}$. On fixe une telle  inclusion. En proc\'edant comme ci-dessus, on pousse d'abord vers la gauche les $\rho\vert.\vert_{F}^{\zeta_{i_{j}}(a-1)/2}$ o\`u $\zeta_{j}=+$. Puis on pousse les autres. On montre que pour tout $j=1,2$ tel que $\zeta_{j}=-$, il existe   un demi-entier $e_{j}$ tel que $[e_{j},\zeta_{j}(a-1)/2]$ soit un segment d\'ecroissant et une inclusion
 $$
 \pi\hookrightarrow (\times_{j=1,2; \zeta_{j}=+}\rho\vert.\vert_{F}^{(a-1)/2}) \times(\times_{j=1,2; \zeta_{j}=-}<e_{j},\zeta_{j}(a-1)/2>_{\rho})\times \pi',
 $$
 pour $\pi'$ irr\'eductible convenable. Comme ci-dessus on v\'erifie que ou bien $\zeta_{j}=+$, ou bien $e_{j}\geq (a-3)/2$ et que l'\'egalit\'e est n\'ecessaire par maximalit\'e de $a$ (ici on utilise le fait que l'on a d\'ej\`a pouss\'e tous les $\rho\vert.\vert_{F}^{(a-1)/2}$ en premi\`ere position). On note $m_{+}$ le nombre de $j=1,2$ tels que $\zeta_{j}=+$ et $m_{-}=2-m_{+}$.  On pose
 $$\tau=\rho\vert .\vert _{F}^{(a-1)/2}\times...\times\rho\vert .\vert _{F}^{(a-1)/2}\times<(a-3)/2,-(a-1)/2>_{\rho}\times...\times <(a-3)/2,-(a-1)/2>_{\rho},$$
 avec $m_{+}$ copies de la premi\`ere repr\'esentation et $m_{-}$ copies de la seconde. Ainsi $\pi'$ est une repr\'esentation irr\'eductible ayant un mod\`ele de Whittaker et dont le support cuspidal \'etendu s'obtient \`a partir $Jord(\Pi)$ en enlevant le support cuspidal de $\tau$ et celui de $\check{\tau}$. Autrement dit, le support cuspidal \'etendu de $\pi'$ s'obtient en rempla\c{c}ant d'abord dans $Jord(\Pi)$ les $2$ copies de $(\rho,a)$ par $2$ copies de $(\rho,a-2)$ puis en enlevant $2m_{-}$ copies de $(\rho,a-2)$; il se peut que $2m_{-}>2$ mais comme dans la d\'emonstration du cas $m=1$ cela ne g\^ene pas. Ainsi le support cuspidal \'etendu de $\pi'$ est le support cuspidal d'une repr\'esentation temp\'er\'ee et $\pi'$ est donc une repr\'esentation temp\'er\'ee. On pose $m=inf(m_{+},m_{-})$.

On suppose que $m\neq 0$, donc $m=m_{+}=m_{-}=1$, et on remarque que l'inclusion $\pi\hookrightarrow \tau\times \pi'$ se factorise necessairement par le sous-module irr\'eductble de $\tau$ ayant un mod\`ele de Whittaker. Celui-ci est   $<(a-1)/2,  -(a-1)/2>_{\rho}$. La repr\'esentation $\pi'$ a pour support cuspidal \'etendu celui de $\pi$ dont on a enlev\'e les deux  copies de $(\rho,a)$. C'est donc bien le support cuspidal \'etendu d'une repr\'esentation temp\'er\'ee et par hypoth\`ese de r\'ecurrence, on sait que $\pi'$ est temp\'er\'ee. Ainsi $\pi$ est un sous-module d'une repr\'esentation induite temp\'er\'ee et est donc temp\'er\'ee.

On suppose donc que $m=0$. On a donc soit $m_{+}=2$ et $m_{-}=0$, soit $m_{-}=2$ et $m_{+}=0$. Dans le cas o\`u $m_{+}=2$, on sait que $\pi$ est une repr\'esentation temp\'er\'ee en appliquant le lemme 2.6. Si $m_{-}=2$, on sait que $\pi$ est le sous-module de Langlands de l'induite $\tau \times \pi'$. Cette induite est donc irr\'eductible d'apr\`es [Mu] théorème 1.1. D'o\`u encore $$\pi\simeq \check{\tau}\times \pi'$$
$$\hookrightarrow  \rho\vert .\vert _{F}^{(a-1)/2}\times \rho\vert .\vert _{F}^{(a-1)/2}\times <(a-3)/2, -(a-3)/2>_{\rho}\times <(a-3)/2,-(a-3)/2>_{\rho}\times \pi'.
$$
On trouve encore une inclusion
$$
\pi\hookrightarrow \rho\vert.\vert_{F}^{(a-1)/2}\times   \rho\vert.\vert_{F}^{(a-1)/2}\times \pi'',
$$
o\`u $\pi''$ est une repr\'esentation temp\'er\'ee. L'ensemble $Jord(\pi'')$ se d\'eduit de $Jord(\Pi)$ en rempla\c{c}ant les $2$ copies de $(\rho,a)$ par $2$ copies de $(\rho,a-2)$. En particulier $Jord(\pi'')$ ne contient pas $(\rho,a)$ et contient au moins $2$ copies de $(\rho,a-2)$ (on peut avoir $a=2$ bien que le cas o\`u $a$ est pair ait d\'ej\`a \'et\'e d\'emontr\'e). Il suffit d'appliquer le lemme  2.6 pour conclure que $\sigma$ est une repr\'esentation temp\'er\'ee. Cela termine la preuve. $\square$

\bigskip

\subsection{ D\'efinition des points de r\'eductibilit\'e possible pour une repr\'esentation temp\'er\'ee}

Soit $\pi$ une repr\'esentation irr\'eductible temp\'er\'ee de $G(F)$. On note $RP(\pi)$ (pour ''r\'eductibilit\'e possible'' ) l'ensemble des couples $(\rho,x)$, o\`u $\rho$ une repr\'esentation irr\'eductible cuspidale unitaire d'un groupe lin\'eaire et $x$ un nombre r\'eel, tel que:

$\bullet$ $\rho\vert.\vert_{F}^{x}\times \pi$ ait pour support cuspidal \'et\'endu le support cuspidal d'une repr\'esentation temp\'er\'ee de $GL(\hat{d}_{G}+2d_{\rho},F)$;

$\bullet$ de plus,  si $x=0$,  $(\rho,1)$ est de bonne parit\'e mais n'appartient pas \`a $Jord(\pi)$.

 Ce sont exactement les couples d\'ecrits dans la remarque 2.3 sauf quand $x=0$ o\`u on a  restreint les possibilit\'es. Remarquons que, pour $(\rho,x)\in RP(\pi)$, on a $\check{\rho}=\rho$ et  $(\rho,-x)\in RP(\pi)$.

 \ass{Proposition}{Soient $\pi$ une repr\'esentation irr\'eductible temp\'er\'ee de $G(F)$, $\rho$ une repr\'esentation cuspidale unitaire irr\'eductible d'un groupe lin\'eaire et $x$ un r\'eel non nul.  
On suppose que $(\rho,x)\notin RP(\pi)$. Alors l'induite $\rho\vert.\vert_{F}^{x}\times \pi$ est irr\'eductible.}

Preuve. Le cas $x=0$ r\'esulte de la classification des repr\'esentations temp\'er\'ees.
Dans tout ce qui suit, on suppose $x\not=0$ et, par sym\'etrie, on peut supposer $x>0$.

On appelle s\'erie discr\`ete strictement positive une s\'erie discr\`ete telle que tous les termes de son module de Jacquet cuspidal soient de la forme $(\otimes_{i=1,...,v}\rho_{i}\vert.\vert_{F}^{s_{i}})\otimes \pi_{cusp}$ o\`u   $\pi_{cusp}$ est une repr\'esentation irr\'eductible cuspidale d'un groupe de m\^eme type que $G$,  $v$ est un entier convenable et, pour $i=1,...,v$, $\rho_{i}$ est une repr\'esentation irr\'eductiible cuspidale unitaire d'un groupe lin\'eaire et $s_{i}$ est un r\'eel  strictement positif. Remarquons que les $s_{i}$ sont forc\'ement des demi-entiers.

 On suppose d'abord que $\pi$ est une s\'erie discr\`ete strictement positive. On va alors montrer  la propri\'et\'e suivante:   l'induite $\rho\vert.\vert_{F}^{x}\times \pi$ est irr\'eductible ou elle contient un sous-quotient irr\'eductible qui est une repr\'esentation temp\'er\'ee. En effet, soit $\sigma$ un sous-quotient irr\'eductible de l'induite $\rho\vert.\vert_{F}^{x}\times \pi$.

On regarde les termes d'un module de Jacquet cuspidal de $\sigma$ c'est-\`a-dire les inclusions 
$$
\sigma\hookrightarrow (\times_{i=1,...,v+1}\rho_{i}\vert.\vert_{F}^{s_{i}})\times \pi_{cusp}. \eqno (1)
$$
D'apr\`es la propri\'et\'e de positivit\'e de $\pi$ tous les $s_{i}$ qui interviennent sont strictement positifs sauf \'eventuellement un, pour $i=i_{0}$ disons, tel que $\rho_{i_{0}}\simeq \check{\rho}$ et $s_{i_{0}}=-x$. Si $i_{0}=1$, alors $\sigma$ est le quotient de Langlands de l'induite $\rho\vert.\vert_{F}^{x}\times \pi$ d'apr\`es le lemme 2.5. Supposons que $i_{0}>1$ et supposons d'abord que $x\neq 1/2$. Cette hypoth\`ese assure que, pour tout $i=1...i_{0}-1$, les induites $\rho_{i}\vert.\vert_{F}^{s_{i}}\times \check{\rho}\vert.\vert_{F}^{-x}$ sont irr\'eductibles car $s_{i}\in 1/2 {\mathbb N}$ et $s_{i}+x $ est un nombre r\'eel positif qui n'est pas $1$. Ainsi l'inclusion ci-dessus donne une inclusion analogue mais o\`u $\check{\rho}\vert.\vert_{F}^{-x}$ a \'et\'e pouss\'e \`a la premi\`ere place. Dans ce cas $\sigma$ est  le quotient de Langlands. Donc si $x\neq 1/2$, $\sigma$ est soit une s\'erie discr\`ete positive soit est le quotient de Langlands. Ainsi si l'induite $\rho\vert.\vert_{F}^{x}\times \pi$ est r\'eductible, elle contient un sous-quotient qui est une s\'erie discr\`ete et $x\in RP(\pi)$ comme annonc\'e.

Le cas o\`u $x=1/2$ est de m\^eme nature; on ne peut pas ''pousser'' $\check{\rho}\vert.\vert_{F}^{-1/2}$ \`a la premi\`ere place s'il existe $i=1,...,i_{0}-1$ tel que $\rho_{i}\simeq \check{\rho}$ et $s_{i}=1/2$. Mais dans ce cas, le terme du module de Jacquet que l'on consid\`ere v\'erifie la propri\'et\'e de positivit\'e large qui caract\'erise les repr\'esentations temp\'er\'ees. Donc ici, soit $\sigma$ est une repr\'esentation temp\'er\'ee soit $\sigma$ est le quotient de Langlands de l'induite $\rho\vert.\vert_{F}^{x}\times \pi$. Et on conclut comme ci-dessus.

On suppose maintenant que $\pi$ est une repr\'esentation temp\'er\'ee quelconque et on prouve la proposition par r\'ecurrence sur le rang de $G$. On suppose que $x\notin RP(\pi)$ et on montre que l'induite $\rho\vert.\vert_{F}^{x}\times \pi$ est irr\'eductible.

On suppose que $\pi$ n'est pas une s\'erie discr\`ete strictement positive puisque ce cas a d\'ej\`a \'et\'e vu. Alors il existe:

  $\bullet$ une repr\'esentation irr\'eductible cuspidale unitaire $\rho_{0}$ d'un groupe lin\'eaire; 
  
  $\bullet$  un segment $[e_{0},f_{0}]$ form\'e de demi-entiers tel que $e_{0}\geq0\geq f_{0}$; 
  
  $\bullet$ une repr\'esentation  irr\'eductible temp\'er\'ee $\pi'$ d'un groupe de m\^eme type que $G$;
  
  \noindent de sorte que l'on ait une inclusion:
$$
\pi\hookrightarrow <e_{0},f_{0}>_{\rho_{0}}\times \pi'.
$$
Si $\pi$ n'est pas de la s\'erie discr\`ete, cela r\'esulte de 2.2(1). Si $\pi$ est de la s\'erie discr\`ete, c'est le lemme 3.1 de [M2]. Si $\rho_{0}\not\simeq \check{\rho}_{0}$ n\'ecessairement $e_{0}=-f_{0}$ et d'apr\`es le lemme 1.4:
$$
RP(\pi)= \begin{cases}RP(\pi')\cup(\rho_{0},e_{0}+1) \cup (\rho_{0},-f_{0}+1) \hbox{ si } \rho_{0}\simeq \check{\rho}_{0},\\
RP(\pi')\cup(\rho_{0},e_{0}+1) \cup (\check{\rho}_{0},e_{0}+1) \hbox{ si } \rho_{0}\not\simeq\check{\rho}_{0}.
\end{cases}
$$
Ainsi si $(\rho,x)\notin RP(\pi)$, l'induite
$
\rho\vert.\vert_{F}^{x}\times <e_{0},f_{0}>_{\rho_{0}}
$
est irr\'eductible car soit $\rho\not\simeq \rho_{0}$, soit $\rho\simeq \rho_{0}$ et $x\neq d_{0}+1$, (on rappelle que $x>0$ ce qui force $x\neq f_{0}-1$). De m\^eme soit $\rho\not\simeq \check{\rho}_{0}$, soit $\rho\simeq \check{\rho}_{0}$ et $x\neq -f_{0}+1$, d'o\`u l'irr\'eductibilit\'e de l'induite $<-f_{0},-e_{0}>_{\check{\rho}_{0}}\times \rho\vert.\vert_{F}^{x}$ et par dualit\'e celle de $<e_{0},f_{0}>_{\rho_{0}}\times \check{\rho}\vert.\vert_{F}^{-x}$. Comme $x\notin RP(\pi')$, par l'hypoth\`ese de r\'ecurrence, on sait aussi que $\rho\vert.\vert_{F}^{x}\times \pi'$ est irr\'eductible. On a donc une suite de morphismes:
$$
\rho\vert.\vert_{F}^{x}\times \pi \hookrightarrow \rho\vert.\vert_{F}^{x}\times <e_{0},f_{0}>_{\rho_{0}}\times \pi' \simeq  <e_{0},f_{0}>_{\rho_{0}}\times\rho\vert .\vert _{F}^x\times \pi'$$
$$\simeq <e_{0},f_{0}>_{\rho_{0}}\times \check{\rho}\vert .\vert _{F}^{-x}\times \pi' \simeq \check{\rho}\vert.\vert_{F}^{-x}\times <e_{0},f_{0}>_{\rho_{0}}\times \pi'.
$$
L'irr\'eductibilit\'e cherch\'ee r\'esulte alors du lemme 2.5 et cela termine la preuve. $\square$

\bigskip

\subsection{R\'eductibilit\'e et g\'en\'ericit\'e}

 \ass{Lemme}{ Soit $\rho$ une repr\'esentation irr\'eductible d'un groupe lin\'eaire, que l'on suppose cuspidale, unitaire et autoduale. Soient $e,f$ des demi-entiers tels que $e\geq0$ et $e+f\not=0$.  On suppose que  $(\rho,2e+1)$ et $(\rho,2\vert f\vert +1)$ ont la bonne parit\'e. Supposons $G$ quasi-d\'eploy\'e, soit $\pi$ une repr\'esentation temp\'er\'ee de $G(F)$, irr\'eductible, ayant un mod\`ele de Whittaker. Alors

(i) si $f\leq0$, l'induite $<e, f>_{\rho}\times \pi$ est r\'eductible;

(ii)  supposons que $f>0$ et que $(\rho,f)\in RP(\pi)$; alors l'induite $<e, f>_{\rho}\times \pi$ est r\'eductible.}

Preuve. Les hypoth\`eses de (i) comme de (ii) assurent que tout sous-quotient irr\'eductible de l'induite \'ecrite a m\^eme support cuspidal \'etendu qu'une repr\'esentation temp\'er\'ee $\pi'$: dans le cas (i), $Jord(\pi')=Jord(\pi)\cup\{(\rho,2d+1),(\rho,-2f+1)\}$ et dans le cas (ii), $Jord(\pi')$ se d\'eduit de $Jord(\pi)$ en rempla\c{c}ant $(\rho,2f-1)$ par $(\rho,2d+1)$. Or ces induites ont un sous-quotient irr\'eductible ayant un mod\`ele de Whittaker. On applique alors le lemme 2.7 pour montrer que ce sous-quotient irr\'eductible est temp\'er\'e. Puisque l'induite n'est pas temp\'er\'ee, elle est r\'eductible. $\square$

\bigskip

\subsection{La notion de liaison}

On fixe  une repr\'esentation irr\'eductible cuspidale et unitaire $\rho$ d'un groupe lin\'eaire et un segment $[e,f]$ de nombres r\'eels.  Soit $\pi$ une repr\'esentation irr\'eductible temp\'er\'ee  de $G(F)$.  On dit que $(\rho,e,f)$ et $Jord(\pi)$  sont  li\'es si les conditions suivantes sont satisfaites:

$e$ est un demi-entier et

(1) si  $(\rho,2\vert e\vert +1)$ n'est pas  de bonne parit\'e, il existe un entier $a\geq1$ de m\^eme parit\'e que $2e+1$ tel que    $(\rho,a)\in Jord(\pi)$ et  les segments $[e,f]$ et $[(a-1)/2,-(a-1)/2]$ sont li\'es au sens de Zelevinsky;

(2) si $(\rho,2\vert e\vert +1)$ est de bonne parit\'e, alors soit  $e\geq -1/2$ et $f\leq1/2$, soit il existe $(\rho,a)\in Jord(\pi)$ avec $(a+1)/2\in [d,f]\cup [-f,-d]$.

La notion de liaison ne d\'epend que de $Jord(\pi)$ et non de $\pi$. De plus, $(\rho,e,f)$ et $Jord(\pi)$ sont li\'es si et seulement si $(\check{\rho},-f,-e)$ et $Jord(\pi)$ le sont.

 \ass{Proposition} {On suppose $G$ quasi-d\'eploy\'e. Soit $(\rho,e,f)$ comme ci-dessus et soit $\pi$ une repr\'esentation irr\'eductible temp\'er\'ee de $G(F)$ ayant un mod\`ele de Whittaker. On suppose que $(\rho,e,f)$ et $Jord(\pi)$ sont li\'es et que $e+f\not=0$. Alors l'induite $<e, f>_{\rho}\times \pi$ est r\'eductible.}
 
 Preuve. Par sym\'etrie, on suppose $e+f>0$, a fortiori $e>0$. On suppose d'abord que $(\rho,e,f)$ satisfait la propri\'et\'e (1) de la condition de liaison. On fixe $a$ comme dans cette propri\'et\'e; comme $(\rho,a)$ n'a pas bonne parit\'e (puisque $(\rho,2e+1)$ ne l'a pas) $\pi$ est une induite de la forme
$<(a-1)/2,  -(a-1)/2>_{\rho}\times \pi'$, o\`u $\pi'$ est une repr\'esentation irr\'eductible temp\'er\'ee convenable.
On a donc un isomorphisme:
$$
<e,f>_{\rho}\times \pi\simeq <e,f>_{\rho}\times <(a-1)/2,-(a-1)/2>_{\rho}\times \pi'.
$$
Or l'induite $<e,f>_{\rho}\times <(a-1)/2,-(a-1)/2>_{\rho}$ dans le $GL$ convenable n'est pas irr\'eductible d'o\`u la r\'eductibilit\'e de l'induite de gauche comme annonc\'e. Ici on n'a d'ailleurs pas utilis\'e le fait que $\pi$ a un mod\`ele de Whittaker.

On suppose que c'est (2) de la d\'efinition de la liaison qui est satisfait. En particulier $(\rho,2e+1)$ et $(\rho,2\vert f\vert+1)$ sont de bonne parit\'e. Si $f\leq 1/2$ le lemme 2.9 montre que l'induite est r\'eductible. On suppose donc que $f>1/2$ et qu'il existe un entier $a\geq1$ tel que $(\rho,a)\in Jord(\pi)$ et $(a+1)/2\in [e,f]$.  On  note $\sigma$ le sous quotient irr\'eductible de l'induite $<e, f>_{\rho}\times \pi$ ayant un mod\`ele de Whittaker. Ainsi $\sigma$ est aussi un sous-quotient irr\'eductible de l'induite
$$
<(a-1)/2,f>_{\rho}\times <e,(a+1)/2>_{\rho}\times \pi.
$$
On note $\pi'$ le sous-quotient irr\'eductible de l'induite $<e,(a+1)/2>_{\rho}\times \pi$ ayant un mod\`ele de Whittaker et on sait que $\pi'$ est une repr\'esentation temp\'er\'ee, cf. lemme 2.7. Ainsi $\sigma$ est un sous-quotient de l'induite $<(a-1)/2,f>_{\rho}\times \pi'$. Si l'induite $<e,f>_{\rho}\times \pi$ est irr\'eductible, elle co\"{\i}ncide avec $\sigma$ et $\sigma$ est un sous-module de l'induite $<-f,-e>_{\rho}\times \pi$. Donc un module de Jacquet cuspidal de $\sigma$ contient un terme $\rho\vert .\vert_{F} ^{-f}\otimes...\otimes \rho\vert .\vert _{F}^{-e}\otimes...$. On sait calculer les modules de Jacquet cuspidaux de l'induite $<(a-1)/2,f>_{\rho}\times \pi'$. Pour qu'ils contiennent le terme pr\'ec\'edent,
il faut n\'ecessairement qu'il existe $x\in [-f,-e]$ tel que le module de Jacquet de $\pi'$ contienne un terme de la forme $\rho\vert.\vert_{F}^{-x}\otimes \cdots$ car le facteur $<(a-1)/2,f>_{\rho}$ ne peut contribuer au mieux qu'au sous-segment $[-f,(a-1)/2]$ de $[-f,-e]$. Ceci contredit le fait que $\pi'$ est temp\'er\'ee et termine la preuve. $\square$

\bigskip

\subsection{Un r\'esultat d'irr\'eductibilit\'e}

\ass{Lemme} {  Soient $(\rho,e,f)$ comme en 2.10 et $\pi$ une repr\'esentation irr\'eductible et temp\'er\'ee de $G(F)$. On suppose que $(\rho,e,f)$ et $Jord(\pi)$ ne sont pas li\'es, que $e+f\not=0$ et  que, ou bien $e$ n'est pas un demi-entier, ou bien $e$ est demi-entier et  $(\rho,2\vert e\vert +1)$ n'est pas de  bonne parit\'e. Alors l'induite $<e, f>_{\rho} \times \pi$ est irr\'eductible.}

Preuve. On suppose, par sym\'etrie, que $e+f>0$. On d\'emontre d'abord le lemme dans le cas o\`u  soit $\rho\not\simeq \check{\rho}$, soit $0\notin [e,f]$. On gagne le fait que pour tout $x\in [e,f[$, l'induite (dans un groupe lin\'eaire convenable)
$$
<e,  x>_{\rho}\times \check{\rho}\vert.\vert_{F}^{-x+1}$$
est irr\'eductible. Supposons de plus pour le moment que $\pi$ est une s\'erie discr\`ete. Alors pour tout $x\in[e,f]$, $(\rho,x)\notin RP(\pi)$  (puisque $(\rho,x)$ n'a pas la bonne parit\'e)  et l'induite $\rho\vert.\vert_{F}^{x}\times \pi$ est irr\'eductible donc isomorphe \`a $\check{\rho}\vert.\vert_{F}^{-x}\times \pi$.  De proche en proche, on montre alors que l'on a une inclusion
$$
<e, f>_{\rho} \times \pi \hookrightarrow \check{\rho}\vert.\vert_{F}^{-f}\times \cdots \check{\rho}\vert.\vert_{F}^{-x+1}\times <e,  x>_{\rho}\times \pi
$$
$$
\hookrightarrow \check{\rho}\vert.\vert_{F}^{-f}\times \cdots \times \check{\rho}\vert.\vert_{F}^{-e}\times \pi.
$$
Et cela donne l'irr\'eductibilit\'e annonc\'ee dans l'\'enonc\'e d'apr\`es  le lemme 2.5.

On enl\`eve l'hypoth\`ese que $\pi$ est une s\'erie discr\`ete. S'il n'en est pas ainsi, $\pi$ est une sous-repr\'esentation d'une induite de la forme $<(a-1)/2, -(a-1)/2>_{\rho'}\times \pi'$, o\`u $\pi'$ est une repr\'esentation irr\'eductible temp\'er\'ee. Il est imm\'ediat que $(\rho,e,f)$ et $Jord(\pi')$ ne sont pas li\'es car $Jord(\pi')$ est un sous-ensemble de $Jord(\pi)$. De plus puisque $(\rho,e,f)$ et $Jord(\pi)$ ne sont pas li\'es, on sait que si $\rho\simeq \rho'$, les segment $[(a-1)/2,-(a-1)/2]$ et $[e,f]$ ne sont pas li\'es. Il en est de m\^eme si $\rho\simeq \check{\rho}'$ car si $(\rho,a)\in Jord(\pi)$ alors $(\check{\rho},a)$ aussi.   Par sym\'etrie par rapport \`a $0$, les segments $[(a-1)/2,-(a-1)/2]$ et $[-f,-e]$ ne sont li\'es que si les segments $[(a-1)/2,-(a-1)/2]$ et $[e,f]$ sont li\'es. Ainsi l'induite $<e, f>_{\rho}\times <(a-1)/2, -(a-1)/2>_{\rho'}$ est irr\'eductible et il en est de m\^eme de l'induite $<(a-1)/2, -(a-1)/2>_{\rho'}\times <-f,-e>_{\check{\rho}}$. On admet par r\'ecurrence,  que l'induite $<e, f>_{\rho}\times \pi'$ est irr\'eductible donc isomorphe \`a l'induite $<-f,-d>_{\check{\rho}}\times \pi'$ et on a alors une s\'erie d'isomorphismes:
$$
<e, f>_{\rho} \times \pi \hookrightarrow <e, f>_{\rho}\times <(a-1)/2, -(a-1)/2>_{\rho'}\times \pi'$$
$$
\simeq  <(a-1)/2, -(a-1)/2>_{\rho'}\times <e, f>_{\rho}\times \pi'$$
$$
\simeq <(a-1)/2, -(a-1)/2>_{\rho'}\times <-f, -e>_{\check{\rho}}\times \pi'$$
$$\simeq
<-f, -e>_{\check{\rho}}\times <(a-1)/2, -(a-1)/2>_{\rho'}\times \pi'.
$$
Et l'irr\'eductibilit\'e annonc\'ee r\'esulte du lemme 2.5.

Le cas restant est plus d\'elicat. On suppose donc que $\rho\simeq \check{\rho}$ et que $0\in [e,f]$, c'est-\`a-dire que $e$ et $f$ sont des entiers tels que $e\geq0\geq f $. Pour \'eviter les confusions de signes, on pose $f^+=-f\geq 0$.

On traite d'abord le cas o\`u $\pi$ est une s\'erie discr\`ete strictement positive. Soit $\sigma$ un sous-quotient irr\'eductible de l'induite $<e,f>_{\rho}\times \pi$. On va montrer que soit $\sigma$ est une repr\'esentation temp\'er\'ee, soit $\sigma$ est le quotient de Langlands de l'induite $<e,f>_{\rho}\times \pi$.

On consid\`ere les termes constants cuspidaux de $\sigma$ c'est-\`a-dire les inclusions 
$$
\sigma\hookrightarrow (\times_{i=1,...,v'}\rho'_{i}\vert.\vert_{F}^{s'_{i}})\times \pi_{cusp},\eqno(1)
$$o\`u $\pi_{cusp}$ est irr\'eductible et cuspidale et, pour tout $i=1,...,v'$, $\rho'_{i}$ est une repr\'esentation irr\'eductible cuspidale et unitaire d'un groupe lin\'eaire et $s'_{i}\in {\mathbb R}$. Les termes constants cuspidaux de toute l'induite sont index\'es par le choix d'un entier $x\in [e+1,-f^+]$ 
et d'un terme constant cuspidal pour $\pi$, c'est-\`a-dire d'une inclusion:
$$
\pi\hookrightarrow (\times_{i=1,...,v}\rho_{i}\vert.\vert_{F}^{s_{i}})\times \pi_{cusp}
$$
et s'obtiennent en m\'elangeant les 3 ensembles ordonn\'ees suivants (en gardant l'ordre dans chaque ensemble)
$$
\cup_{i\in [e,x]}\{\rho\vert.\vert_{F}^{i}\}; \cup_{i\in [f^+,-x[}\{\rho\vert.\vert_{F}^{i}\}; \cup_{i=1,...,v}\{\rho_{i}\vert.\vert_{F}^{s_{i}}\}.
$$
Un tel terme v\'erifie la condition de positivit\'e des repr\'esentations temp\'er\'ees si $x\leq f^++1$. Le point est donc de d\'emontrer que si $\sigma$ contient un terme comme ci-dessus avec $x>f^++1$ alors $\sigma$ est le quotient de Langlands de l'induite $<e, -f^+>_{\rho}\times \pi$.

On fixe donc $\sigma$ et un terme comme ci-dessus avec $x> f^++1$. On utilise le fait que pour tout $i=1,...,v$, si $\rho_{i}\simeq \rho$ alors $s_{i}$ est un demi-entier non entier: $(\rho,2s_{i}+1)$ doit avoir bonne parit\'e, or $e$ est entier et $(\rho,2e+1)$ n'a pas bonne parit\'e. Donc pour tout entier $y\in [e,-f^+]\cup [f^+,-e]$ et pour tout $i=1,...,v$, l'induite $\rho\vert.\vert_{F}^{y}\times \rho_{i}\vert.\vert_{F}^{s_{i}}$ est irr\'eductible. On peut donc commuter ces facteurs. De m\^eme pour tout $y\in [e,x]$ et tout $y'\in [f^+,-x[$ l'induite $\rho\vert.\vert_{F}^{y}\times \rho\vert.\vert_{F}^{y'}$ est irr\'eductible car $y-y'\geq x-f^+>1$. On peut donc aussi commuter de tels facteurs et finalement l'assertion sur $\sigma$ se traduit par l'existence d'une inclusion:
$$
\sigma\hookrightarrow (\times_{j\in [f^+,-x]}\rho\vert.\vert_{F}^{j})\times(\times_{i=1,...,v}\rho_{i}\vert.\vert_{F}^{s_{i}})\times(\times_{y\in [e,x[}\rho\vert.\vert_{F}^{y})\times \pi_{cusp}.
$$
Pour tout $y\in [e,x[$, $\rho\vert.\vert_{F}^{y}\times \pi_{cusp}$ est irr\'eductible (proposition 2.8) et donc isomorphe \`a $\rho\vert.\vert_{F}^{-y}\times \pi_{cusp}$. D'o\`u finalement une inclusion
$$
\sigma\hookrightarrow (\times_{j\in [f^+,-x]}\rho\vert.\vert_{F}^{j})\times(\times_{i=1,...,v}\rho_{i}\vert.\vert_{F}^{s_{i}})\times(\times_{y\in [x,e]}\rho\vert.\vert_{F}^{-y})\times \pi_{cusp} 
$$
$$
\simeq (\times_{j\in [f^+,-e]}\rho\vert.\vert_{F}^{j})\times(\times_{i=1,...,v}\rho_{i}\vert.\vert_{F}^{s_{i}})\times \pi_{cusp}.
$$
Et $\sigma$ est le quotient de Langlands de l'induite comme annonc\'e, cf. lemme 2.5. Ainsi  $<e, f>_{\rho}\times \pi$ a un seul sous-quotient irr\'eductible qui n'est pas temp\'er\'e, c'est le sous-quotient de Langlands. Mais une telle induite ne peut avoir de sous-quotient irr\'eductible qui sont des repr\'esentations temp\'er\'ees car son  support cuspidal \'etendu  est l'union de celui de $\pi$ avec les 2 segments $[e,-e]\cup [f^+,-f^+]$ bas\'es sur $\rho$.  Ceci n'est pas le support cuspidal \'etendu d'une repr\'esentation temp\'er\'ee car $(\rho,2e+1)$ n'a pas bonne parit\'e, ni d'ailleurs $(\rho,2f^++1)$, et que $e\neq f^+$ par hypoth\`ese. D'o\`u l'irr\'eductibilit\'e.

On consid\`ere maintenant   une repr\'esentation irr\'eductible temp\'er\'ee $\pi$ quelconque et on prouve  l'irr\'eductibilit\'e par r\'ecurrence. Si $\pi$ n'est pas une s\'erie discr\`ete strictement positive, alors,  comme on l'a dit en 2.8, il existe une repr\'esentation cuspidale $\rho'$ (unitaire irr\'eductible), un segment $[e',-f']$ de demi-entiers avec $e',f'\geq0$, et une repr\'esentation irr\'eductible temp\'er\'ee $\pi'$ avec une inclusion:
$$
\pi\hookrightarrow <e',-f'>_{\rho'}\times \pi'.
$$
De plus si $\rho'\simeq \rho$ soit $(\rho,2e'+1)$ et $(\rho,2f'+1)$ sont de bonne parit\'e ce qui entra\^{\i}ne que $e'-e$ et $f'+f$ sont des demi-entiers non entiers, soit $(\rho,2e'+1)$ et $(\rho,2f'+1)$ ne sont pas de bonne parit\'e ce qui force $e'=f'$  et, par l'hypoth\`ese de non liaison, que les segments $[e,f]$ et $[e',-f']$ ne sont pas li\'es. Ainsi l'induite $<e, f>_{\rho}\times <e',-f'>_{\rho'}$ est irr\'eductible. On v\'erifie de la m\^eme fa\c{c}on que l'induite et $<e',-f'>_{\rho'}\times <-f,-e>_{\rho}$ est irr\'eductible. On a donc:
$$
<e,f>_{\rho}\times \pi\hookrightarrow <e,f>_{\rho}\times <e',-f'>_{\rho'}\times \pi'\simeq
$$
$$
<e,-f'>_{\rho'}\times <e,f>_{\rho}\times \pi'.
$$On sait d'apr\`es le lemme 2.4 que $Jord(\pi)= Jord(\pi')\cup\{(\rho,2e'+1),(\rho,2f'+1)\}$ et ainsi $(\rho,e,f)$ n'est pas li\'e \`a $Jord(\pi')$. 
Par l'hypoth\`ese de r\'ecurrence $<e,f>_{\rho}\times \pi'$ est donc aussi irr\'eductible et on a encore une inclusion
$$
<e,f>_{\rho}\times \pi\hookrightarrow <e',-f'>_{\rho'}\times <-f,-e>_{\rho}\times \pi'$$
$$\simeq
<-f,-e>_{\rho}\times <d',-f'>_{\rho'}\times \pi'.
$$
Avec le lemme  2.5, cette inclusion montre l'irr\'eductibilit\'e de l'induite $<e,f>_{\rho}\times \pi$. Cela termine la preuve. $\square$

 \bigskip
 
 \subsection{Un second r\'esultat d'irr\'eductibilit\'e}
 
\ass{Proposition}{Soient $(\rho,e,f)$ comme en 2.10 et  $\pi$ une repr\'esentation  irr\'eductible temp\'er\'ee de $G(F)$. On suppose que $(\rho,e,f)$ et $Jord(\pi)$ ne sont pas li\'es. Alors l'induite $<e,f>_{\rho}\times \pi$ est irr\'eductible.}

Preuve. En tenant compte du lemme pr\'ec\'edent, il nous reste \`a voir le cas o\`u $e$ est demi-entier et $(\rho,2\vert e\vert +1)$ est de bonne parit\'e. On peut encore supposer $e+f>0$. L'hypoth\`ese de non liaison assure alors que $f>1/2$ et que pour tout $x\in [e,f]$, $\rho\vert.\vert_{F}^x\notin RP(\pi)$. En particulier, pour un tel $x$, on sait que $\rho\vert.\vert_{F}^x\times \pi$ est irr\'eductible (proposition 2.8), donc isomorphe \`a $\rho\vert.\vert_{F}^{-x}\times \pi$.  De plus pour tous $y,y'\in [d,f]$, $y+y'\geq 2f>1$ et les induites $\rho\vert.\vert_{F}^{y}\times \rho\vert.\vert_{F}^{-y'}$ sont donc aussi irr\'eductibles. On a donc un isomorphisme:
$$
(\times_{x\in [e,f]}\rho\vert.\vert_{F}^{x})\times \pi\simeq (\times_{ x\in [-f,-e]}\rho\vert.\vert_{F}^{x})\times \pi.
$$
D'o\`u
$$
<e,f>_{\rho}\times \pi\hookrightarrow( \times_{x\in [e,f]}\rho\vert.\vert_{F}^{x})\times \pi\simeq (\times_{ x\in [-f,-e]}\rho\vert.\vert_{F}^{x})\times \pi.
$$
Et on conclut, avec le lemme 2.5, \`a l'irr\'eductibilit\'e de l'induite. $\square$

\bigskip

\subsection{Crit\`ere d'irr\'eductibilit\'e}

 G\'en\'eralisons la notion de liaison. Soit $k\in{\mathbb N}$ et pour tout $i=1,...,k$ soit $(\rho_{i},e_{i},f_{i})$ un triplet form\'e d'une repr\'esentation irr\'eductible $\rho_{i}$ cuspidale et unitaire d'un groupe lin\'eaire et d'un segment   $[e_{i},f_{i}]$ de nombres r\'eels. Soit aussi $\pi$ une repr\'esentation irr\'eductible temp\'er\'ee de $G(F)$. On dit que $\{(\rho_{i},e_{i},f_{i})_{i=1,...,k}\}$ et $Jord(\pi)$ ne sont pas li\'es si pour tout $i=1,...,k$, $(\rho_{i},e_{i},f_{i})$ et $Jord(\pi)$ ne sont pas li\'es et si, pour tout $ j=1,...,k$, $j\not=i$,  les induites $<e_{i},f_{i}>_{\rho_{i}}\times <e_{j},f_{j}>_{\rho_{j}}$ et $<e_{i},f_{i}>_{\rho_{i}}\times <-f_{j},-e_{j}>_{\check{\rho}_{j}}$ sont irr\'eductibles.

 \ass{Th\'eor\`eme}{Soient $\{(\rho_{i},e_{i},f_{i})_{i=1,...,k}\}$ comme ci-dessus et $\pi$ une repr\'esentation irr\'eductible temp\'er\'ee de $G(F)$. On suppose que pour tout $i=1,...,k$, $e_{i}+f_{i}\neq 0$.

(i) On suppose que  $G$ est quasi-d\'eploy\'e et que $\pi$ a un mod\`ele de Whittaker. Alors l'induite $(\times_{i=1,...,k}<e_{i},f_{i}>_{\rho_{i}})\times \pi$ est irr\'eductible  seulement si $\{(\rho_{i},e_{i},f_{i})_{i=1,...,k}\}$ et $Jord(\pi)$ ne sont pas li\'es.

(ii) On suppose que $\{(\rho_{i},e_{i},f_{i})_{i=1,...,k}\}$ et $Jord(\pi)$ ne sont pas li\'es. Alors l'induite $(\times_{i=1,...,k}<e_{i},f_{i}>_{\rho_{i}})\times \pi$ est irr\'eductible.}

En d'autres termes, la condition de non liaison est n\'ecessaire et suffisante pour avoir l'irr\'eductibilit\'e de l'induite si $\pi$ a un mod\`ele de Whittaker et est seulement suffisante sans cette hypoth\`ese de mod\`ele de Whittaker.

Preuve. Prouvons (i). On suppose que l'induite figurant dans cette assertion est irr\'eductible. Il faut certainement que pour tout $i, j=1,...,k$, $i\not=j$, les induites $<e_{i},f_{i}>_{\rho_{i}}\times <e_{j},f_{j}>_{\rho_{j}}$ soient irr\'eductibles. Pour un choix de $j=1,...,k$, l'irr\'eductibilit\'e de l'induite de l'assertion est aussi \'equivalente \`a l'irr\'eductibilit\'e de la m\^eme induite mais o\`u on remplace $<e_{j},f_{j}>_{\rho_{j}}$ par sa contragr\'ediente $<-f_{j},-e_{j}>_{\check{\rho}_{j}}$; on doit donc avoir aussi l'irr\'eductibilit\'e pour tout $ i,j=1,...,k$, $i\not=j$, des induites $<e_{i},f_{i}>_{\rho_{i}}\times <-f_{j},-e_{j}>_{\check{\rho}_{j}}$. De plus  pour tout $i=1,...,k$, l'induite $<e_{i},f_{i}>_{\rho_{i}}\times \pi$ doit aussi \^etre irr\'eductible et donc, d'apr\`es la proposition 2.10,  $(\rho_{i},e_{i},f_{i})$ et $Jord(\pi)$ ne sont pas li\'es. Cela montre que les conditions de non liaison sont bien n\'ecessaires pour avoir l'irr\'eductibilit\'e quand $\pi$ a un mod\`ele de Whittaker.

Prouvons (ii).  On suppose que, pour tout $i=1,...,k$, $e_{i}+f_{i}> 0$, sinon on remplace $<e_{i},f_{i}>_{\rho_{i}}$ par $<-f_{i},-e_{i}>_{\check{\rho}_{i}}$.

En permutant, comme on en a le droit, on suppose aussi que $e_{1}+f_{1}\geq \cdots \geq e_{k}+f_{k}$.
Avec ces choix, l'induite
$$
(\times_{i=1,...,k}<e_{i},f_{i}>_{\rho_{i}})\times \pi  \eqno(1)
$$
a un unique quotient irr\'eductible, $\sigma$. De plus $\sigma$  intervient avec multiplicit\'e 1 comme sous-quotient irr\'eductible et l'induite
$$
(\times_{i=1,...,k}<-f_{i},-e_{i}>_{\check{\rho}_{i}})\times \pi \simeq(\times_{i=k,...,1}<-f_{i},-e_{i}>_{\check{\rho}_{i}})\times \pi  \eqno(2)
$$
a un unique sous-module irr\'eductible et il est isomorphe \`a $\sigma$. Pour d\'emontrer l'irr\'eductibilit\'e  de l'induite (1), il suffit donc de montrer que (1) est un sous-module de (2) et on va m\^eme v\'erifier que (1) et (2) sont isomorphes.
On a les isomorphismes
$$
(\times_{i=1,...,k}<e_{i},f_{i}>_{\rho_{i}})\times \pi \simeq (\times_{i=1,...,k-1} <e_{i},f_{i}>_{\rho_{i}})\times <-f_{k},-e_{k}>_{\check{\rho}_{k}}\times \pi $$
$$\simeq
<-f_{k},-e_{k}>_{\check{\rho}_{k}}\times (\times_{i=1,...,k-1} <e_{i},f_{i}>_{\rho_{i}})\times \pi.$$
Et de proche en proche on construit donc un isomorphisme de (1) et (2). Ceci termine la preuve. $\square$

\bigskip

\subsection{Irr\'eductibilit\'e et $L$-paquets g\'en\'eriques}

Soit $\underline{G}$ la forme int\'erieure quasi-d\'eploy\'ee de $G$. Pour un $L$-paquet $\Pi$ de repr\'esentations temp\'er\'ees de $G(F)$, on note $\underline{\Pi}$ le $L$-paquet de repr\'esentations temp\'er\'ees de $\underline{G}(F)$ qui correspond \`a $\Pi$, autrement dit dont le param\`etre de Langlands est le m\^eme que celui de $\Pi$. On a $Jord(\Pi)=Jord(\underline{\Pi})$.  Dans l'\'enonc\'e ci-dessous, l'induite d\'esigne une repr\'esentation de $G(F)$ ou de $\underline{G}(F)$ selon que $\pi$ appartient \`a $\Pi$ ou \`a $\underline{\Pi}$.

\ass{Corollaire}{Soit $k\in {\mathbb N}$. Pour tout $i=1,...,k$, soit $\sigma_{i}$ une repr\'esentation irr\'eductible d'un groupe lin\'eaire, unitaire et de la s\'erie discr\`ete, et soit $s_{i}$ un nombre r\'eel  non nul. Soit $\Pi$ un $L$-paquet de repr\'esentations irr\'eductibles temp\'er\'ees de $G(F)$. Alors l'induite
$$
(\times_{i=1,...,k}\sigma_{i}\vert.\vert_{F}^{s_{i}}) \times \pi
$$
est irr\'eductible pour tout $\pi\in \Pi\sqcup \underline{\Pi}$ si et seulement s'il existe une repr\'esentation $\pi_{0}\in \underline{\Pi}$ ayant un mod\`ele de Whittaker et tel que l'induite pr\'ec\'edente soit irr\'eductible pour $\pi=\pi_{0}$.
}

Preuve. Il est \'evident que si l'induite de l'\'enonc\'e est irr\'eductible pour tout $\pi\in \Pi\sqcup \underline{\Pi}$, elle l'est en particulier pour $\pi=\pi_{0}$ ayant un mod\`ele de Whittaker. R\'eciproquement on suppose que cette induite est irr\'eductible pour une repr\'esentation $\pi=\pi_{0}$ ayant un mod\`ele de Whittaker. On \'ecrit, pour tout $i=1,...,k$, $\sigma_{i}\vert.\vert_{F}^{s_i}$ sour la forme $<e_{i},f_{i}>_{\rho_{i}}$ et on a $s_{i}=e_{i}+f_{i}\neq 0$. Le (i) du th\'eor\`eme pr\'ec\'edent dit que les $\{(\rho_{i},e_{i},f_{i})_{i=1,...,k}\}$ et $Jord(\Pi)$ ne sont pas li\'es. Le (ii) du m\^eme th\'eor\`eme donne alors l'irr\'eductibilit\'e de l'induite pour tout $\pi\in \Pi\sqcup\underline{ \Pi}$. Cela termine la preuve. $\square$

\bigskip

\subsection{Le cas des groupes sp\'eciaux orthogonaux pairs}

On suppose que $G$ est sp\'ecial orthogonal pair.   La principale diff\'erence avec les cas symplectique et sp\'ecial orthogonal impair est que notre notation simplifi\'ee $\sigma_{1}\times...\times \sigma_{m}\times \pi$ pour les induites est trop impr\'ecise car d'une part, il peut y avoir deux  sous-groupes paraboliques non conjugu\'es de m\^eme L\'evi $GL(d_{1})\times ...\times GL(d_{t})\times G'$, d'autre part, l'identification du groupe central $G'$ n'est canonique qu'\`a conjugaison pr\`es par un \'el\'ement du groupe orthogonal tout entier (ce qui fait que selon l'identification choisie, $\pi$ est remplac\'e par $\pi^w$). On peut rem\'edier \`a cela de la fa\c{c}on suivante. Fixons une d\'ecomposition
$$V=Fv_{1}\oplus...\oplus Fv_{t}\oplus V_{an}\oplus Fv_{-t}\oplus...\oplus Fv_{-1}$$
de sorte que

$\bullet$ les vecteurs $v_{i}$ sont isotropes et $q(v_{i},v_{-i})=1$ pour $i=\pm 1,...,\pm t$;

$\bullet$ les espaces $Fv_{1}\oplus Fv_{-1}$,..., $Fv_{t}\oplus Fv_{-t}$ et $V_{an}$ sont deux \`a deux orthogonaux;

$\bullet$ la restriction de $q$ \`a $V_{an}$ est anisotrope.

Appelons famille parabolique une famille $(I_{j})_{j=1,...,m}$ v\'erifiant les conditions suivantes:

$\bullet$ chaque $I_{j}$ est un sous-ensemble non vide de $\{\pm 1,...,\pm t\}$;

$\bullet$ en posant $\check{I}_{j}=\{-i; i\in I_{j}\}$, les ensembles $I_{1}$,...,$I_{m}$, $\check{I}_{1}$,...,$\check{I}_{m}$ sont deux \`a deux disjoints.

Pour une telle famille et pour $j=1,...,t$, notons $X_{j}$ le sous-espace de $V$ engendr\'e par les vecteurs $v_{i}$ pour $i\in \cup_{k=1,...,j}I_{k}$. Notons $P_{I_{1},...,I_{m}}$ le sous-groupe parabolique de $G$ form\'e des \'el\'ements qui conservent le drapeau de sous-espaces
$$X_{1}\subset X_{2}\subset...\subset X_{m}.$$
La composante de L\'evi de $P_{I_{1},...,I_{m}}$ s'identifie \`a $GL(d_{1})\times...\times GL_{d_{m}}\times G'$ o\`u, pour tout $j$, $d_{j}$ est le nombre d'\'el\'ements de $I_{j}$ et $G'$ est le groupe sp\'ecial orthogonal du sous-espace de $V$ engendr\'e par $V_{an}$ et les $v_{i}$ pour $i\not\in \bigcup_{j=1,....,m}(I_{j}\cup\check{I}_{j})$. Cette identification est canonique (\`a automorphismes int\'erieurs pr\`es pour les blocs $GL(d_{j})$, mais cela n'a pas d'importance). Supposons d'abord $V_{an}\not=\{0\}$. Alors l'application qui, \`a une famille parabolique $(I_{j})_{j=1,...,m}$, associe $P_{I_{1},...,I_{m}}$ est injective. La classe de conjugaison de ce parabolique est d\'etermin\'ee par la suite d'entiers $d_{1},...,d_{m}$. Supposons maintenant $V_{an}=\{0\}$. L'application n'est plus injective. Le d\'efaut d'injectivit\'e est le suivant. Consid\'erons une famille parabolique $(I_{j})_{j=1,...,m}$ telle qu'il existe $i_{0}\in \{1,...,t\}$ de sorte que
$$(\bigcup_{j=1,...,m}(I_{j}\cup\check{I}_{j}))=\{\pm 1,...,\pm t\}\setminus\{i_{0},-i_{0}\}.$$
On peut compl\'eter cette famille en deux autres, en lui ajoutant un \'el\'ement $I_{m+1}$ \'egal soit \`a $\{i_{0}\}$, soit \`a $\{-i_{0}\}$. On a
$$P_{I_{1},...,I_{m}}=P_{I_{1},...,I_{m},\{i_{0}\}}=P_{I_{1},...,I_{m},\{-i_{0}\}}.$$
Remarquons que les L\'evi de ces paraboliques ne sont pas exactement les m\^emes. Le premier a un bloc central $G'$ qui est un groupe $SO(2)$ d\'eploy\'e, les deux autres n'ont pas de bloc central, mais ont un facteur $GL(1)$ suppl\'ementaire. Cela correspond aux deux identifications possibles de $SO(2)$ d\'eploy\'e avec $GL(1)$. Par ailleurs, la suite $d_{1},...,d_{m}$ ne d\'etermine plus toujours la classe de conjugaison du parabolique.

Dans les constructions des paragraphes pr\'ec\'edents, il convient maintenant, \`a chaque fois que l'on se donne une induite $\sigma_{1}\times...\times\sigma_{m}\times \pi$, de fixer auparavant une famille parabolique $(I_{j})_{j=1,...,m}$ et de pr\'eciser que l'induite est $Ind_{P_{I_{1},...,I_{m}}}^G(\sigma_{1}\otimes...\otimes\sigma_{m}\otimes \pi)$. Au cours d'une d\'emonstration, la famille parabolique peut changer. Par exemple, dans le cas o\`u $m=2$, dire que l'on peut permuter $\sigma_{1}$ et $\sigma_{2}$ signifie que l'on a 
$$Ind_{P_{I_{1},I_{2}}}^G(\sigma_{1}\otimes\sigma_{2}\otimes \pi)=Ind_{P_{I_{2},I_{1}}}^G(\sigma_{2}\otimes \sigma_{1}\otimes \pi).$$
De m\^eme, quand $Ind_{P_{I_{1}}}^G(\sigma\otimes \pi)$ admet un quotient de Langlands, celui-ci est une sous-repr\'esentation de $Ind_{P_{\check{I}_{1}}}^G(\check{\sigma}\otimes \pi)$.

Comme on le voit, le groupe $SO(2)$ d\'eploy\'e peut intervenir dans nos raisonnements. Il convient de consid\'erer que, pour ce groupe, aucune repr\'esentation n'est de la s\'erie discr\`ete, a fortiori aucune n'est cuspidale. Les repr\'esentations temp\'er\'ees $\pi$ sont les caract\`eres unitaires. Pour une telle repr\'esentation, on a $Jord(\pi)=\{\rho,\rho^{-1}\}$, o\`u $\rho$ est un caract\`ere unitaire de $GL(1,F)$. L'\'el\'ement $\rho$ a bonne parit\'e si et seulement si $\rho=\rho^{-1}$, c'est-\`a-dire $\rho$ est un caract\`ere quadratique.

Certains raisonnements utilisent la caract\'erisation des repr\'esentations temp\'er\'ees par la positivit\'e de leurs exposants cuspidaux. La notion de positivit\'e est a priori diff\'erente pour un groupe sp\'ecial orthogonal pair de ce qu'elle est pour un groupe symplectique ou  sp\'ecial orthogonal impair. Par exemple, pour un exposant   cuspidal $(b_{1},...,b_{t})$ relatif \`a un sous-groupe parabolique minimal, la condition de positivit\'e s'exprime par les in\'egalit\'es
$$b_{1}\geq0,\,b_{1}+b_{2}\geq0,...,b_{1}+...+b_{t}\geq0\eqno(1)$$
dans les cas symplectique ou sp\'ecial orthogonal impair, tandis qu'il faut ajouter la condition
$$b_{1}+...+b_{t-1}-b_{t}\geq0\eqno(2)$$
dans le cas sp\'ecial orthogonal pair. Mais l'utilisation des familles paraboliques ci-dessus \'elimine cette diff\'erence. En effet, posons $I_{j}=\{j\}$ pour tout $j=1,...,t$. Si $(b_{1},...,b_{t})$ est un exposant relatif au parabolique $P_{I_{1},...,I_{t}}$, alors $(b_{1},...,b_{t-1},-b_{t})$ est un exposant relatif \`a $P_{I_{1},...,I_{t-1},\check{I}_{t}}$. Si on impose les conditions (1) pour chacune des familles paraboliques, la condition (2) en r\'esulte.

Dans certains raisonnements, par exemple la preuve du lemme 2.6, il convient de remplacer les paquets $\Pi$ par les paquets plus grossiers $\bar{\Pi}$ d\'efinis en 2.1. Enfin, la premi\`ere propri\'et\'e de 2.7 n'est plus vraie. Il peut y avoir deux repr\'esentations $\pi_{1}$ et $\pi_{2}$ de $G(F)$ ayant un mod\`ele de Whittaker d'un type fix\'e et v\'erifiant $Jord(\pi_{1})=Jord(\pi_{2})$. Mais alors $\pi_{2}=\pi_{1}^w$ et cette unicit\'e \`a conjugaison pr\`es par le groupe orthogonal tout entier suffit pour assurer la validit\'e de nos raisonnements.

Avec ces adaptations, on v\'erifie que les r\'esultats des paragraphes pr\'ec\'edents et leurs preuves restent valables.

{\bf Remarque.} On a pr\'esent\'e en 2.1 et admis la forme fine des conjectures. On peut se contenter d'un forme un peu plus faible o\`u l'on ne param\`etre plus que les paquets $\bar{\Pi}$, cf. [W1] 4.4 pour des \'enonc\'es pr\'ecis. Nos raisonnements restent valables. L'int\'er\^et de cette forme plus faible des conjectures est qu'elle est un r\'esultat annonc\'e par Arthur, ainsi qu'on l'a dit dans l'introduction.

\bigskip

\section{Preuve du th\'eor\`eme}

\bigskip

Soient $G$ et $G'$ comme dans l'introduction. Rappelons le th\'eor\`eme principal que nous allons maintenant d\'emontrer.

\ass{Th\'eor\`eme}{Soient $\varphi\in \Phi(G)$ et $\varphi'\in \Phi(G')$. On suppose $\varphi$ et $\varphi'$ g\'en\'eriques. Alors:

(i) toutes les repr\'esentations induites dont les \'el\'ements de $\Pi^G(\varphi)$ et de $\Pi^{G'}(\varphi')$ sont les quotients de Langlands sont irr\'eductibles;

(ii) si $E(\varphi,\varphi')=-\mu(G,G')$, on a $m(\sigma,\sigma')=0$ pour tous $\sigma\in \Pi^G(\varphi)$, $\sigma'\in \Pi^{G'}(\varphi')$;

(iii) si $E(\varphi,\varphi')=\mu(G,G')$, on a
$$m(\sigma(\varphi,\boldsymbol{\epsilon}),\sigma'(\varphi',\boldsymbol{\epsilon}'))=1$$
et $m(\sigma,\sigma')=0$ pour tous $\sigma\in \Pi^G(\varphi)$, $\sigma'\in \Pi^{G'}(\varphi')$ tels que $(\sigma,\sigma')\not=(\sigma(\varphi,\boldsymbol{\epsilon}), \sigma'(\varphi',\boldsymbol{\epsilon}'))$.}

 Preuve. Introduisons la forme quasi-d\'eploy\'ee $\underline{G}$ de $G$.   Il y a un L\'evi
$$\underline{L}=GL(d_{1})\times...\times GL(d_{m})\times \underline{G}_{0}$$
 de $\underline{G}$, des repr\'esentations irr\'eductibles temp\'er\'ees $\sigma_{j}$ de $GL(d_{j},F)$ et des r\'eels $b_{j}$ pour $j=1,...,m$, et un param\`etre de Langlands temp\'er\'e $\varphi_{0}\in \Phi(\underline{G}_{0})$ de sorte que la suite $(b_{1},...,b_{m})$ soit strictement positive relativement \`a un sous-groupe parabolique de composante de L\'evi $\underline{L}$ et que $\Pi^{\underline{G}}(\varphi)$ soit l'ensemble des quotients de Langlands des induites 
 $$\sigma_{1}\vert .\vert _{F}^{b_{1}}\times...\times \sigma_{m}\vert .\vert _{F}^{b_{m}}\times \pi_{0}\eqno(1)$$
 pour $\pi_{0}\in \Pi^{\underline{G}_{0}}(\varphi_{0})$. 
 
 {\bf Remarque.} Ici, on peut supposer et on suppose que groupe $\underline{G}_{0}$ n'est pas un groupe $SO(2)$ d\'eploy\'e: on peut remplacer un tel groupe par un facteur $GL(1)$.
 
 Si $\underline{L}$ ne correspond \`a aucun L\'evi de $G$ (cela se produit quand $\underline{G}_{0}=\{1\}$ et $G$ n'est pas quasi-d\'eploy\'e), le paquet $\Pi^G(\varphi)$ est vide. Sinon, $\underline{L}$ correspond \`a un L\'evi
  $$L=GL(d_{1})\times...\times GL(d_{m})\times G_{0}$$
   de $G$ et $\Pi^G(\varphi)$ est form\'e des quotients de Langlands des induites (1) pour $\pi_{0}\in \Pi^{G_{0}}(\varphi_{0})$. 
 
 Par hypoth\`ese, il existe $\pi\in \Pi^{\underline{G}}(\varphi)$ qui admet un mod\`ele de Whittaker d'un certain type. C'est le quotient de Langlands de l'induite (1) pour une repr\'esentation $\pi_{0}\in \Pi^{\underline{G}_{0}}(\varphi_{0})$ qui admet elle-aussi un mod\`ele de Whittaker. D'apr\`es [Mu] théorème 1.1, le quotient de Langlands de l'induite  admet un mod\`ele de Whittaker si et seulement si l'induite est irr\'eductible. Donc cette induite est irr\'eductible. En appliquant le corollaire 2.14 au paquet $\Pi^{\underline{G}_{0}}(\varphi_{0})$ et, quand il existe, au paquet $\Pi^{G_{0}}(\varphi_{0})$, on voit que toutes les induites (1) sont irr\'eductibles pour $\pi_{0}\in \Pi^{\underline{G}_{0}}(\varphi_{0})$ ou $\pi_{0}\in \Pi^{G_{0}}(\varphi_{0})$. Cela d\'emontre l'assertion concernant $\Pi^G(\varphi)$ du  (i) du th\'eor\`eme.  
 
 La m\^eme chose vaut du c\^ot\'e de $G'$. On introduit les objets similaires, auquels on ajoute des $'$. L'assertion concernant $\Pi^{G'}(\varphi')$ du (i) du th\'eor\`eme s'obtient comme ci-dessus.  Supposons $E(\varphi,\varphi')=-\mu(G,G')$. Si l'un des paquets $\Pi^G(\varphi)$ ou $\Pi^{G'}(\varphi')$ est vide, l'assertion (ii) l'est aussi. Supposons ces deux paquets non vides, a fortiori, les L\'evi $L$ et $L'$ existent. Soient $\epsilon\in{\cal E}^G(\varphi)$ et $\epsilon'\in {\cal E}^{G'}(\varphi')$. Alors $\sigma(\varphi,\epsilon)$ est l'induite (1) pour $\pi_{0}=\sigma(\varphi_{0},\epsilon)$. On a une assertion analogue pour $\sigma'(\varphi',\epsilon')$. D'apr\`es la proposition 1.3, on a
 $$m(\sigma(\varphi,\epsilon),\sigma'(\varphi',\epsilon'))=m(\sigma(\varphi_{0},\epsilon),\sigma'(\varphi'_{0},\epsilon')).\eqno(2)$$
 On a $E(\varphi,\varphi')=E(\varphi_{0},\varphi'_{0})$ et $\mu(G,G')=\mu(G_{0},G'_{0})$, donc $E(\varphi_{0},\varphi'_{0})=-\mu(G_{0},G'_{0})$.  D'apr\`es [W1] th\'eor\`eme 4.9, on a $m(\sigma(\varphi_{0},\epsilon),\sigma'(\varphi'_{0},\epsilon'))=0$, d'o\`u l'\'egalit\'e cherch\'ee $m(\sigma(\varphi,\epsilon),\sigma'(\varphi',\epsilon'))=0$.
 
 Supposons maintenant $E(\varphi,\varphi')=\mu(G,G')$. Montrons qu'alors, les L\'evi  $L$ et $L'$ existent bel et bien. Si ce n'est pas le cas, l'un des groupes $G_{0}$ ou $G'_{0}$ est r\'eduit \`a $\{1\}$. Le param\`etre $\varphi_{0}$ ou $\varphi'_{0}$ correspondant est  vide et on a  $E(\varphi,\varphi')=E(\varphi_{0},\varphi'_{0})=1$. Alors $\mu(G,G')=1$ et les deux groupes $G$ et $G'$ sont quasi-d\'eploy\'es. Donc les L\'evi $L$ et $L'$ existent, contrairement \`a l'hypoth\`ese. De plus, l'hypoth\`ese $E(\varphi,\varphi')=\mu(G,G')$ entra\^{\i}ne que le couple $(\boldsymbol{\epsilon},\boldsymbol{\epsilon}')$ param\`etre un \'el\'ement du produit $\Pi^G(\varphi)\times \Pi^{G'}(\varphi')$. Ces paquets sont donc non vides. La preuve du (iii) du th\'eor\`eme est alors la m\^eme que celle de (ii): pour $\epsilon\in {\cal E}^G(\varphi)$ et $\epsilon'\in {\cal E}^{G'}(\varphi')$, on a l'\'egalit\'e (2); en appliquant le th\'eor\`eme 4.9 de [W1]  au membre de droite de cette \'egalit\'e, on obtient l'assertion cherch\'ee. $\square$
 
 {\bf Remarque.} Supposons que $F$ soit le compl\'et\'e d'un corps de nombres $k$ en une place $v$ de ce corps et que $G$ soit la composante en cette place $v$ d'un groupe sp\'ecial orthogonal ${\bf G}$ d\'efini sur $k$. Notons ${\mathbb A}$ l'anneau des ad\`eles de $k$. Consid\'erons une repr\'esentation automorphe irr\'eductible $\boldsymbol{\pi}$ de ${\bf G}({\mathbb A})$ intervenant dans le spectre discret. Conjecturalement, il lui correspond une famille $(\boldsymbol{\pi}_{i})_{i=1,...,k}$ de repr\'esentations automorphes irr\'eductibles de groupes lin\'eaires intervenant dans le spectre discret. En notant $d_{i}$ l'entier tel que $\boldsymbol{\pi}_{i}$ soit une repr\'esentation automorphe de $GL(d_{i},{\mathbb A})$, on a $\hat{d}_{G}=\sum_{i=1,...,k}d_{i}$. Ce r\'esultat conjectural est annonc\'e par Arthur, sous les m\^emes r\'eserves que dans le cas  local, cf. l'introduction. Supposons que toutes les $\boldsymbol{\pi}_{i}$ soient cuspidales. Notons $\pi=\boldsymbol{\pi}_{v}$ la composante de $\boldsymbol{\pi}$ en la place $v$. Alors $\pi$ appartient \`a un paquet $\Pi$ auquel on peut appliquer nos r\'esultats, c'est-\`a-dire que c'est un paquet g\'en\'erique. Cela r\'esulte de [M3] proposition 5.1.

\bigskip

{\bf Bibliographie}

[AGRS] A. Aizenbud, D. Gourevitch, S. Rallis, G. Schiffmann: {\it Multiplicity one theorems}, pr\'epublication 2007

[A1] J. Arthur: {\it An introduction to the trace formula}, Clay Math. Proc. 4 (2005), p.1-253

[A2] -----------: {\it A local trace formula}, Publ. Sc. IHES 73 (1991), p.5-96

[BZ] I. Bernstein, A. Zelevinsky: {\it Induced representations of reductive $p$-adic groups I}, Ann. Sc. ENS 10 (1977), p.441-472

[GGP] W. Gan, B. Gross, D. Prasad: {Symplectic local root numbers, central critical $L$-values and restriction problems in the representation theory of classical groups}, pr\'epublication 2009

[K]  T. Konno: {\it Twisted endoscopy implies the generic packet conjecture}, Isra\" el J. of Math. 129 (2002), p.253-289

[M1] C. Moeglin: {\it Repr\'esentations quadratiques unipotentes des groupes classiques $p$-adiques}, Duke Math. J. 84 (1996), p.267-332

[M2] --------------: {\it Sur la classification des s\'eries discr\`etes des groupes classiques $p$-adiques: param\`etres de Langlands et exhaustivit\'e}, J. Eur. Math. Soc. 4 (2002), p.143-200

[M3] --------------: {\it Image des op\'erateurs d'entrelacements normalis\'es et p\^oles des s\'eries d'Eisenstein}, pr\'epublication 2009

[MW] C. Moeglin, J.-L. Waldspurger: {\it Sur le transfert des traces d'un groupe classique $p$-adique \`a un groupe lin\'eaire tordu}, Selecta math. 12 (2006), p.433-515

[Mu] G. Muic: {\it A proof of Casselman-Shahidi's conjecture for quasi-split classical groups}, Can. Math. Bull. 43 (2000), p.90-99

[W1] J.-L. Waldspurger: {\it La conjecture locale de Gross-Prasad pour les repr\'esentations temp\'er\'ees des groupes sp\'eciaux orthogonaux}, pr\'epublication 2009

[W2] -----------------------: {\it Calcul d'une valeur d'un facteur $\epsilon$ par une formule int\'egrale}, pr\'epublication 2009

[W3] -----------------------: {\it Une formule int\'egrale reli\'ee \`a la conjecture locale de Gross-Prasad, $2^{\text{\`eme}}$ partie: extension aux repr\'esentations temp\'er\'ees}, pr\'epublication 2009

[Z] A. Zelevinsky: {\it Induced representations of reductive $p$-adic groups II. On irreducible representations of $GL(n)$}, Ann. Sc. ENS 13 (1980), p.165-210

 \bigskip
 
 CNRS Institut de math\'ematiques de Jussieu
 
 175 rue du Chevaleret
 
 75013 Paris
 
 moeglin@math.jussieu.fr, waldspur@math.jussieu.fr

\end{document}